\documentclass[11pt,a4paper,reqno]{amsart}
\pdfoutput=1
\usepackage{amsmath,amsthm,latexsym,amssymb,wasysym}
\usepackage[margin=1.0in,tmargin=0.9in,bmargin=0.80in]{geometry}
\usepackage{hyperref}
\usepackage{bbm}
\usepackage{xargs}
\usepackage{enumerate}
\usepackage{enumitem}
\setlist[enumerate]{label={\upshape(\roman*)}}
\usepackage[textsize=footnotesize]{todonotes}
\usepackage[flushleft]{threeparttable}
\usepackage{multirow}
\usepackage{longtable}
\usepackage{times}
\usepackage{microtype}
\usepackage{tabularx}
\usepackage{ltablex}
\keepXColumns
\usepackage{caption} 
\usepackage[backend=bibtex, style=alphabetic, maxbibnames=10]{biblatex}
\usepackage{listings}

\usepackage{xcolor}
\definecolor{codegreen}{rgb}{0,0.6,0}
\definecolor{codegray}{rgb}{0.5,0.5,0.5}
\definecolor{codepurple}{rgb}{0.58,0,0.82}
\definecolor{backcolour}{rgb}{0.95,0.95,0.92}

\lstdefinestyle{myStyle}{
	language=Python,
	backgroundcolor=\color{white},   
	commentstyle=\color{codegreen},
	keywordstyle=\color{magenta},
	numberstyle=\tiny\color{codegray},
	stringstyle=\color{codegreen},
	basicstyle=\ttfamily\footnotesize,  
	breakatwhitespace=false,         
	breaklines=true,                 
	captionpos=b,                    
	keepspaces=true,                 
	numbers=left,                    
	numbersep=5pt,
	showspaces=false,                
	showstringspaces=false,
	showtabs=false,                  
	tabsize=2,
	frame=single
}

\usepackage{float}
\definecolor{bred}{rgb}{0.8,0,0}
\hypersetup{colorlinks,linkcolor={blue},citecolor={bred},urlcolor={blue}}

\newtheorem{theorem}{Theorem}[section]
\newtheorem{proposition}[theorem]{Proposition}
\newtheorem{lemma}[theorem]{Lemma}
\newtheorem{corollary}[theorem]{Corollary}
\newtheorem{remark}[theorem]{Remark}
\newtheorem{definition}[theorem]{Definition}

\newtheorem{assumption}{Assumption}

\newcommand{\E}{\mathbb{E}}
\newcommand{\N}{\mathbb{N}}
\newcommand{\R}{\mathbb{R}}
\newcommand{\1}{\mathbbm{1}}
\newcommand{\dee}{\mathrm{d}}
\newcommand{\mc}{\mathcal}
\newcommand{\tr}{\operatorname{tr}}
\newcommand{\lfrf}[1]{\left\lfloor #1 \right\rfloor}
\newcommand{\lcrc}[1]{\left\lceil #1 \right\rceil}
\newcommand{\lara}[1]{\left\langle #1 \right\rangle}
\newcommand{\lbrb}[1]{\left| #1 \right|}
\newcommand{\lsrs}[1]{\left[ #1 \right]}
\newcommand{\ce}[1]{\left. #1 \right|}
\allowdisplaybreaks

\begin{document}

\title[]{Non-asymptotic convergence bounds for modified tamed unadjusted Langevin algorithm in non-convex setting}

\author[A. Neufeld]{Ariel Neufeld}
\author[M. Ng]{Matthew Ng Cheng En}
\author[Y. Zhang]{Ying Zhang} 

\address{Division of Mathematical Sciences, Nanyang Technological University, 21 Nanyang Link, 637371 Singapore}
\email{ariel.neufeld@ntu.edu.sg}

\address{Division of Mathematical Sciences, Nanyang Technological University, 21 Nanyang Link, 637371 Singapore}
\email{matt0037@e.ntu.edu.sg}

\address{Financial Technology Thrust, The Hong Kong University of Science and Technology~(Guangzhou)}
\email{yingzhang@hkust-gz.edu.cn}

\date{}
\thanks
{
}

\begin{abstract}
We consider the problem of sampling from a high-dimensional target distribution $\pi_\beta$ on $\R^d$ with density proportional to $\theta\mapsto e^{-\beta U(\theta)}$ using explicit numerical schemes based on discretising the Langevin stochastic differential equation (SDE). In recent literature, taming has been proposed and studied as a method for ensuring stability of Langevin-based numerical schemes in the case of super-linearly growing drift coefficients for the Langevin SDE. In particular, the Tamed Unadjusted Langevin Algorithm (TULA) was proposed in \cite{brosse2019tamed} to sample from such target distributions with the gradient of the potential $U$ being super-linearly growing. However, theoretical guarantees in Wasserstein distances for Langevin-based algorithms have traditionally been derived assuming strong convexity of the potential $U$. In this paper, we propose a novel taming factor and derive, under a setting with possibly non-convex potential $U$ and super-linearly growing gradient of $U$, non-asymptotic theoretical bounds in Wasserstein-1 and Wasserstein-2 distances between the law of our algorithm, which we name the modified Tamed Unadjusted Langevin Algorithm (mTULA), and the target distribution $\pi_\beta$. We obtain respective rates of convergence $\mc{O}(\lambda)$ and $\mc{O}(\lambda^{1/2})$ in Wasserstein-1 and Wasserstein-2 distances for the discretisation error of mTULA in step size $\lambda$. High-dimensional numerical simulations which support our theoretical findings are presented to showcase the applicability of our algorithm.
\end{abstract}

\maketitle

\section{Introduction}
Sampling from a given high-dimensional distribution is a problem of integral importance for applications in fields such as Bayesian statistics \cite{cotter2013mcmc}, machine learning \cite{andrieu2003introduction}, and molecular dynamics \cite{lelievre2016partial}. To this end, a class of algorithms, collectively described as Langevin Monte Carlo (LMC) algorithms, has been developed and studied extensively in the literature. This class of algorithms involves explicit numerical schemes based on the Langevin stochastic differential equation (SDE)
\begin{align}
\dee Z_t = -h(Z_t)\ \dee t + \sqrt{2\beta^{-1}}\ \dee B_t\label{intro:langevin}
\end{align}

for sampling from a target distribution $\pi_\beta$ with Lebesgue density proportional to $\R^d\ni\theta \mapsto e^{-\beta U(\theta)}$, where $d\in\N$, $\beta>0$ is the so-called inverse temperature parameter, $(B_t)_{t\geq 0}$ is a standard Brownian motion in $\R^d$, and $U:\R^d\to\R$ is a continuously differentiable function with $\nabla U = h^T$ \cite{dalalyan2019user}. We note that, for sufficiently large $\beta$, $\pi_\beta$ concentrates around the minimisers of $U$ \cite{hwang}. This connects the optimization problem with the problem of sampling, namely, that minimising $U$ is equivalent to sampling from $\pi_\beta$ (when $\beta$ takes large values). To solve the sampling problem, a typical approach is to consider the Unadjusted Langevin Algorithm (ULA), given as the first-order Euler discretisation of the Langevin SDE (\ref{intro:langevin}). However, while computationally efficient, it has been well-established, e.g., in \cite{mattingly2002ergodicity} and \cite{hutzenthaler2011strong}, that ULA may diverge to infinity in certain sense if the gradient $h$ of the potential $U$ is growing super-linearly, that is, $\liminf_{|\theta|\to\infty}|h(\theta)|/|\theta|=\infty$; we refer to Remark \ref{rmk:ulaunstable} below for a detailed discussion. To address the issue of divergence in the case of super-linearly growing gradient, a technique known as \textit{taming}, where the gradient $h$ is divided by an appropriately chosen \textit{taming factor} in the numerical scheme to control the super-linear growth of $h$, has been introduced in recent literature, see \cite{hutzenthaler2012strong} and \cite{sabanis2013note}. \\

While non-asymptotic theoretical guarantees of ULA and its tamed variants have been widely established in existing literature, such results have typically been obtained under at least one of the following assumptions: (a) -- global Lipschitz continuity on the gradient $h$; (b) -- strong convexity of the potential $U$. More precisely, we say that $h$ is \textit{globally Lipschitz continuous} if there exists a constant $L>0$ such that
\begin{align}
|h(\theta)-h(\theta')|\leq L|\theta-\theta'|,\quad \forall\theta,\theta\in\R^d,\label{intro:lipschitz}
\end{align}
and that the continuously differentiable potential $U$ is \textit{strongly convex} if there exists a constant $m>0$ for which its gradient $h$ satisfies
\begin{align}
\lara{\theta-\theta', h(\theta)-h(\theta')}\geq m|\theta-\theta'|^2,\quad \forall\theta,\theta'\in\R^d.\label{intro:strong_convex}
\end{align}

Under both assumptions of global Lipschitz continuity of the gradient $h$ and strong convexity of the potential $U$, various non-asymptotic theoretical bounds on the total variation distance and Wasserstein distance between the law of the algorithm and the sampling target distribution have been derived for ULA in \cite{dalalyan2017theoretical}, \cite{durmus2017nonasymptotic}, and \cite{durmus2019high}. In \cite{cheng2018sharp} and \cite{majka2020nonasymptotic}, the assumption of strong convexity is relaxed by only requiring the condition (\ref{intro:strong_convex}) to hold on a set where $|\theta-\theta'|$ is sufficiently large. Such a condition can be satisfied by non-convex $U$ and is comparable to the `convex at infinity' condition which we impose in our work, stated explicitly in Assumption \ref{assumption:3}. Under this relaxed condition, which the authors of \cite{majka2020nonasymptotic} term as `contractivity at infinity', non-asymptotic upper bounds on Wasserstein distances of orders 1 and 2 with respective orders of convergence $1/2$ and $1/4$ in the step size of the algorithm were obtained for ULA in \cite{majka2020nonasymptotic}. Similarly, a non-asymptotic upper bound on Wasserstein-1 distance with order of convergence $1/2$ in the step size was obtained for the stochastic gradient Langevin dynamics (SGLD) algorithm in \cite{chau2021stochastic}, where a stochastic gradient $H(\theta, X_n)$ dependent on a data sequence $(X_n)_{n\in\N}$ is used in the numerical scheme in place of the exact deterministic gradient $h$. This result was also derived under a non-convex setting, where the typical strong convexity condition was relaxed to a dissipativity condition which is equivalent to Remark \ref{remark:2.4} in our work. Non-asymptotic analyses of Langevin-based algorithms under such a dissipativity condition were also presented in \cite{raginsky2017non} and \cite{xu2018global}. In addition, there are convergence results that are established under certain functional inequalities, see, e.g., \cite{vempala2019rapid}, \cite{erdogdu2022convergence}, \cite{mousavi2023towards}, and references therein. More precisely, \cite{vempala2019rapid} provides a non-asymptotic convergence estimate in Kullback Leibler (KL) divergence under the condition that the target distribution satisfies a Log-Sobolev inequality. \cite{chewi2021analysis} extends the results in \cite{vempala2019rapid} to distributions satisfying a Poincar\'e inequality and provides convergence results for ULA in Chi-squared and R\'enyi divergence.  \cite{mousavi2023towards} then further extends the results to handle the case of weak Poincar\'e inequalities. In addition, \cite{erdogdu2022convergence} provides theoretical guarantees in Chi-squared and R\'enyi divergence for the ULA algorithm to sample from target distributions with potentials satisfying a dissipativity condition, which can also be translated to convergence results in KL divergence, total variation, and Wasserstein-2 distance. However, the results in \cite{cheng2018sharp}, \cite{majka2020nonasymptotic}, \cite{chau2021stochastic}, \cite{raginsky2017non}, \cite{xu2018global}, \cite{vempala2019rapid}, \cite{chewi2021analysis}, \cite{erdogdu2022convergence}, and \cite{mousavi2023towards} do not directly apply to the case of super-linearly growing gradient, as they require the gradient to satisfy global Lipschitz (or H\"older) continuity assumptions.\\

In \cite{brosse2019tamed}, the authors proposed the Tamed Unadjusted Langevin Algorithm (TULA) to address the case of super-linearly growing gradient. The typical global Lipschitz assumption on the gradient was relaxed to a polynomial Lipschitz condition, and non-asymptotic bounds on the total variation distance and Wasserstein-2 distance between the law of the algorithm and the target distribution were derived. A convergence order of up to $1$ in the step size for the Wasserstein-2 distance was established, however, the bounds on Wasserstein-2 distance required imposing additionally the strong convexity assumption on $U$. In \cite{lytras2023taming}, the authors propose a Langevin-dynamics based algorithm called the splitted Tamed Unadjusted Langevin Algorithm (sTULA) to sample from distributions with super-linearly growing potentials satisfying a Log-Sobolev inequality, and provide non-asymptotic convergence guarantees for the performance of sTULA in KL divergence, total variation, and Wasserstein-2 distance. A crucial result in \cite{lytras2023taming} is that, under certain convexity at infinity condition of the target distribution $\pi_\beta$, one can obtain a Log-Sobolev inequality with constant that is independent of the dimension and has exponential dependence on $\beta$. In \cite{johnston2023kinetic}, the authors propose two algorithms called  tKLMC1 and tKLMC2 to sample from log-concave distributions with super-linearly growing potentials in Hamiltonian setting, and provide non-asymptotic error estimates in Wasserstein-2 distance between the laws of the algorithms and the target distribution. In \cite{lovas2020taming} and \cite{lim2023non}, the Tamed Unadjusted Stochastic Langevin Algorithm (TUSLA), a tamed variant of the SGLD algorithm, was proposed and studied under a polynomial Lipschitz assumption on the gradient and in a non-convex setting. Non-asymptotic upper bounds for Wasserstein distances of orders 1 and 2 with respective convergence orders $1/2$ and $1/4$ in the step size were obtained for TUSLA.\\

\textbf{Main Contribution.} In this paper, we propose a new tamed variant of ULA, which we name the modified Tamed Unadjusted Langevin Algorithm (mTULA), and derive theoretical non-asymptotic upper bounds on the Wasserstein distances of orders 1 and 2 between the law of the algorithm and the sampling target distribution. In comparison with \cite{brosse2019tamed}, we propose a different taming factor for the gradient, and our theoretical guarantees in Wasserstein distances are obtained under an equivalent polynomial Lipschitz assumption on the gradient, but in a non-convex setting, as opposed to the typical strong convexity condition imposed in \cite{brosse2019tamed} in order to obtain the bounds in Wasserstein-2 distance. In addition, we impose a comparable `convex at infinity' assumption on the potential $U$ as in \cite{lim2023non}, but achieve markedly improved convergence orders of $1$ and $1/2$ in step size for Wasserstein distances of orders 1 and 2, respectively, compared to $1/2$ and $1/4$, respectively, obtained in \cite{lim2023non}.\\

An overview of the remaining of this paper is as follows. Section \hyperref[subsec:assumptions]{2.1} contains the precise statements of the assumptions we impose to obtain our main results. In particular, the polynomial Lipschitz condition and the `convex at infinity' assumption are respectively stated in Assumptions \ref{assumption:2} and \ref{assumption:3}. Sections \hyperref[subsec:main_results]{2.2} and \hyperref[subsec:comparison]{2.3} present our main results of the non-asymptotic theoretical upper bounds which we derive for mTULA and compare them with that of related works, while Section \hyperref[sec:simulations]{3} contains results of some high-dimensional numerical simulations which support our theoretical findings. Section \hyperref[sec:proofs]{4} presents an overview of the methodology we adopted to establish our theoretical bounds. Finally, Section \hyperref[sec:proof_aux]{5} contains the proofs of auxiliary results and remarks in the paper, while Appendix \hyperref[appendix:A]{A} presents a summary table of the explicit analytic expressions of all constants which appear in our intermediate results and main theorems.\\

We conclude this section by introducing some notation. Let $(\Omega,\mathcal{F},\mathbb{P})$ be a probability space. We denote by $\E[Z]$  the expectation of a random variable $Z$. For $p\in[1,\infty)$, $L^p$ is used to denote the  space of $p$-integrable real-valued random variables. Fix integers $d, m \geq 1$. For an $\R^d$-valued random variable $Z$, its law on $\mathcal{B}(\R^d)$, i.e. the Borel sigma-algebra of $\R^d$, is denoted by $\mathcal{L}(Z)$. For a positive real number $a$, we denote by $\lfrf{a}$ its integer part, and $\lcrc{a} = \lfrf{a}+1 $. The notation $\1_\cdot$ is used to denote indicator functions. Given a normed space $(X, \lVert\cdot\rVert_X)$ and an element $x\in X$, we denote the norm of $x$ by $\lVert x\rVert_X$. In the particular case $X=\R^d$ and $\lVert\cdot\rVert$ is the Euclidean norm, we understand the notation $|x|$ as referring to $|x| = \lVert x\rVert_{\R^d}$ for $x\in\R^d$. Similarly, for a real-valued $m\times d$ matrix $A\in\R^{m\times d}$, we understand $|A|$ as referring to the operator norm $|A| = \sup\{|Ax|: |x|\leq 1, x\in\R^d\}$. The Euclidean scalar product is denoted by $\langle \cdot,\cdot\rangle$. For any integer $q \geq 1$, let $\mathcal{P}(\R^q)$ denote the set of probability measures on $\mathcal{B}(\R^q)$. For $\mu,\mu'\in\mathcal{P}(\R^d)$, let $\mathcal{C}(\mu,\mu')$ denote the set of couplings of $\mu,\mu'$, that is, probability measures $\zeta$ on $\mathcal{B}(\R^{2d})$ such that its respective marginals are $\mu,\mu'$. For two Borel probability measures $\mu$ and $\mu'$ defined on $\R^d$ with finite $p$-th moments, the Wasserstein distance of order $p \geq 1$ is defined as
\begin{align}
{W}_p(\mu,\mu'):=
\left(\inf_{\zeta\in\mathcal{C}(\mu,\mu')}\int_{\R^d\times \R^d}\lbrb{\theta-\theta'}^p\ \dee\zeta(\theta,\theta')\right)^{1/p}.\label{eqn:wassertein_p}
\end{align}

\section{Assumptions and main results}
Let $U:\R^d\to \R$ be a twice continuously differentiable function, and denote by $h:=(\nabla U)^T$ its gradient transposed. With a slight abuse of notation, denote by $\nabla h:=\nabla\nabla^T U$ the Hessian of $U$. Furthermore, for any $\beta>0$, define the sampling target distribution
\begin{align}
\pi_\beta(A):=\frac{\int_A e^{-\beta U(\theta)} \dee\theta}{\int_{\R^d} e^{-\beta U(\theta)} \dee\theta},\quad A\in\mathcal{B}(\R^d),\label{eqn:target_distribution}
\end{align}
where we assume $\int_{\R^d} e^{-\beta U(\theta)} \dee\theta<\infty$.\\

In this paper, we propose the modified Tamed Unadjusted Lagenvin Algorithm (mTULA), which we recursively define by
\begin{align}
\theta_0^\lambda:=\theta_0,\quad \theta_{n+1}^\lambda:=\theta_n^\lambda-\lambda h_\lambda(\theta_n^\lambda) + \sqrt{2\lambda\beta^{-1}}\xi_{n+1},\quad n\in\N_0,\label{eqn:tula}
\end{align}
where $\lambda >0$ is the step size, $(\xi_n)_{n\in \N}$ is a sequence of independent standard $d$-dimensional Gaussian random variables independent of the $\R^d$-valued random variable $\theta_0$, and where for all $\theta\in\R^d$, the tamed gradient $h_\lambda$ takes the form
\begin{align}
h_\lambda(\theta) := \frac{h(\theta)}{(1+\lambda\lbrb{\theta}^{2r})^{1/2}}.\label{eqn:tamed_gradient}
\end{align}

\begin{remark}\label{rmk:ulaunstable}
We highlight that ULA is unstable when sampling from a target distribution $\pi_{\beta}$ whose potential has a super-linearly growing gradient. More precisely, in \cite[Lemma 6.3]{mattingly2002ergodicity}, the authors consider an example where $h(\theta) =\theta^3$, $\theta\in\R$, and show that the Euler–Maruyama discretization of \eqref{intro:langevin} (denoted by $(\theta^{\mathsf{EM}}_n)_{n\in\N_0}$) is not ergodic in the following sense: first, if $\E[\theta^{\mathsf{EM}}_0]\geq 2/\lambda$ with $\lambda>0$ being the step size, then $\lim_{n \to \infty} \E[(\theta^{\mathsf{EM}}_n)^2]=\infty$; second, for any $\theta^{\mathsf{EM}}_0\in \R$ and $\lambda>0$, the sample path $\theta^{\mathsf{EM}}_n$ diverges to infinity with positive probability. In addition, it is shown in \cite[Theorem 1 and Equation (2.7)]{hutzenthaler2011strong} that if $|h(\theta)|\geq |\theta|^l/\mathsf{C}$ and $\sqrt{2\beta^{-1}}\leq \mathsf{C}|\theta|^{l'}$ for all $|\theta|\geq \mathsf{C}$ with respect to constants $l>1, l>l'\geq 0, \mathsf{C}>0$, then the absolute moments of ULA diverge towards infinity in finite time. One example among many super-linearly growing functions that satisfy the conditions in \cite{hutzenthaler2011strong} is again $h(\theta) =\theta^3$, $\theta\in\R$, where we choose $\mathsf{C} = \max\{\sqrt{2\beta^{-1}}, \sqrt{\beta/2}\}$, $l = 3, l'=0$ with some $\beta>0$. Thus, to tackle the aforementioned sampling problem, in this paper, we consider the mTULA algorithm \eqref{eqn:tula}-\eqref{eqn:tamed_gradient}, which can be viewed as a tamed version of ULA. 
\end{remark}

\subsection{Assumptions} \label{subsec:assumptions} In this section, we present the conditions required to establish the main results. Let $r,\nu\in \N_0$ be fixed. Denote $r_*:=\max\{8r+8,4\nu+4,2\nu+2r+4\}$. We impose the following assumptions.

\begin{assumption}
\label{assumption:1}
The initial condition $\theta_0$ has a finite $r_*$-th moment, i.e., $\E
\lsrs{\lbrb{\theta_0}^{r_*}}<\infty$.
\end{assumption}

\begin{assumption}
\label{assumption:2} We state this assumption in two parts:
\begin{itemize}
\item[(a)] There exists a constant $L>0$ such that, for all $\theta,\theta'\in\R^d$,
\begin{align*}
\lbrb{h(\theta)-h(\theta')}\leq L(1+\lbrb{\theta}+\lbrb{\theta'})^{r}\lbrb{\theta-\theta'}.
\end{align*}
\item[(b)] In addition, there exists a constant $K>0$ such that, for all $\theta\in\R^d$, 
\begin{align*}
\lbrb{h(\theta)}\leq K(1+\lbrb{\theta}^{r+1}).
\end{align*}
\end{itemize}
\end{assumption}

\begin{remark}
Note that it is possible that $r=0$. In particular, our results also include, but do not require, the case where the gradient $h$ is globally Lipschitz continuous.
\end{remark}

\begin{remark}
Note that from the (polynomial) Lipschitz condition of Assumption \hyperref[assumption:2]{2-(a)}, one can deduce a polynomial growth condition on $h$ in the form of Assumption \hyperref[assumption:2]{2-(b)} with $K=\max\left\{L2^{r-1}+|h(0)|, L2^r\right\}$. However, such a choice of $K$ can be unnecessarily large which, in view of equation (\ref{rf:max_step}), leads to a maximum step size restriction $\tilde{\lambda}_{\max}$ of the algorithm which is unnecessarily small. As an example, consider $h(\theta)=(|\theta|^2-1)\theta$, i.e. the gradient of the double-well potential . From Proposition \ref{proposition:3.1}, Assumption \hyperref[assumption:2]{2-(a)} is satisfied for this choice of gradient with $r=2, L=1$. Thus, the constant from the polynomial growth condition derived from Assumption \hyperref[assumption:2]{2-(a)} is given by $\max\left\{L2^{r-1}+|h(0)|, L2^r\right\}=4$. However, by direct computation, one sees that Assumption \hyperref[assumption:2]{2-(b)} holds for this choice of $h$ with $K=2$. This is the reason we choose to impose separately a polynomial growth condition on $h$.
\end{remark}

\begin{assumption}
\label{assumption:3}
We state this assumption in the following two cases:
\begin{enumerate}
\item If $r>0$, then there exist constants $a,b>0$ and $\bar{r}\in[0,r)$ such that, for all $\theta,\theta'\in \R^d$, 
\begin{align}
\lara{\theta-\theta',h(\theta)-h(\theta')}\geq a\lbrb{\theta-\theta'}^2(\lbrb{\theta}^{r}+\lbrb{\theta'}^{r})-b\lbrb{\theta-\theta'}^2(\lbrb{\theta}^{\bar{r}}+\lbrb{\theta'}^{\bar{r}}).\label{eqn:convex_at_infinity}
\end{align}
\item If $r=0$, then there exist constants $\tilde{a},\tilde{b}>0$ such that, for all $\theta\in\R^d$,
\begin{align}
\lara{\theta,h(\theta)} \geq \tilde{a}\lbrb{\theta}^{2}-\tilde{b}.\label{eqn:dissipativity}
\end{align}
\end{enumerate}
\end{assumption}

\begin{remark}
\label{remark:2.3}
In the case $r>0$, Assumption \ref{assumption:3} is a `convex at infinity' condition and can be understood as follows. When the norms of $\theta,\theta'$ are large, the first term of the RHS of condition (\ref{eqn:convex_at_infinity}) dominates and it essentially behaves like the condition
\begin{align}
\lara{\theta-\theta',h(\theta)-h(\theta')}\geq a\lbrb{\theta-\theta'}^2(\lbrb{\theta}^{r}+\lbrb{\theta'}^{r}).\label{eqn:local_convex}
\end{align}
Condition (\ref{eqn:local_convex}) is comparable to the usual strong convexity condition (\ref{intro:strong_convex}), but with the constant $m$ in (\ref{intro:strong_convex}) being dependent on $\theta,\theta'$. When the norms of $\theta,\theta'$ are large, condition (\ref{eqn:local_convex}) implies condition (\ref{intro:strong_convex}) for a suitable choice of the constant $m$. However, on domains where the norms of $\theta,\theta'$ are sufficiently small, the second term of the RHS of condition (\ref{eqn:convex_at_infinity}) dominates and $\lara{\theta-\theta',h(\theta)-h(\theta')}$ can be negative. That is, $U$ may not satisfy (strong) convexity on domains where $\theta,\theta'$ are sufficiently small. This means that our convex at infinity assumption accommodates for possible non-convexity of the potential function $U$, but it comes at the expense of imposing a stronger assumption on $U$ than strong convexity when the norms of $\theta,\theta'$ are large. \\

As an illustration, consider again $h(\theta)=(|\theta|^2-1)\theta$, i.e. the gradient of the double-well potential. Then one can verify that $h$ does not satisfy condition (\ref{eqn:local_convex}), nor is it strongly convex. However, by Proposition \ref{proposition:3.1}, Assumption \ref{assumption:3} is satisfied with $r=2$, $\bar{r}=0$, $a=1/2$, and $b=1$. Moreover, since, for this choice of $h$,
\begin{align}
\lara{\theta-\theta', h(\theta)-h(\theta')}\geq \tfrac{1}{2}\left(|\theta|^2+|\theta'|^2-2\right)|\theta-\theta'|^2\nonumber
\end{align}
for all $\theta,\theta'\in\R^d$, this implies that 
\begin{align}
\lara{\theta-\theta', h(\theta)-h(\theta')}\geq |\theta-\theta'|^2\nonumber
\end{align}
for all $\theta,\theta'\in\R^d$ satisfying $|\theta|,|\theta'|\geq \sqrt{2}$. That is, the double-well potential is strongly convex only on subdomains where $\theta,\theta'$ have sufficiently large norms, but it is not (strongly) convex on $\R^d$ due to its behaviour when $\theta,\theta'$ are near $0$. This gives intuition to the choice of terminology `at infinity' and illustrates how our assumption can be satisfied by potentials $U$ which are not strongly convex.\\
\end{remark}

\begin{remark}
\label{rmk:conatinffunineq}
Let us comment on the connection between the convexity at infinity condition in Assumption~\ref{assumption:3} and functional inequalities. By \cite[Theorem 5.3]{lytras2023taming}, under Assumptions~\ref{assumption:2} and~\ref{assumption:3}, it holds that $\pi_\beta$ defined in \eqref{eqn:target_distribution} satisfies a Log-Sobelev inequality with a constant that is independent of the dimension and has exponential dependence on $\beta$.
\end{remark}

Under Assumptions \ref{assumption:1}, \ref{assumption:2}, \ref{assumption:3}, one can obtain dissipativity conditions for the gradient $h$. The explicit statement is given in the following remark which proof is given in Section \hyperref[proof:remark2.4]{5.1}.
\begin{remark}
\label{remark:2.4}
By Assumptions \ref{assumption:1}, \ref{assumption:2}, \ref{assumption:3}, there exist constants $\bar{a},\bar{b}>0$ such that, for all $\theta\in\R^d$,
\begin{align*}
\lara{\theta,h(\theta)} \geq \bar{a}\lbrb{\theta}^{r+2}-\bar{b},
\end{align*}
where the $\bar{a}, \bar{b}$ are explicitly given as
\begin{align}
\bar{a} &:= \frac{a}{2}\1_{\{r>0\}}+\tilde{a}\1_{\{r=0\}},\nonumber\\
\bar{b} &:= \left(\left(b+\tfrac{a}{2}\right)R^{\bar{r}+2}+\tfrac{K^2}{2a}\right)\1_{\{r>0\}}+\tilde{b}\1_{\{r=0\}},\nonumber\\
R &:=\max\left\{\left(\tfrac{4b}{a}\right)^{1/(r-\bar{r})},2^{1/r}\right\}.\nonumber
\end{align}
Moreover, it holds that for all $\theta\in \R^d$,
\begin{align*}
\lara{\theta,h(\theta)} \geq \bar{a}\lbrb{\theta}^{2}-\bar{b}',
\end{align*}
where $\bar{b}':=(\bar{b}+2^{2/r}\bar{a})\1_{\{r>0\}} + \tilde{b}\1_{\{r=0\}}$.
\end{remark}
\textit{Proof.} See Section \hyperref[proof:remark2.4]{5.1}.\\

\begin{remark}
\label{remark:2.5}
Under Assumptions \ref{assumption:1}, \ref{assumption:2}, \ref{assumption:3}, one obtains a one-sided Lipschitz continuity condition on $h$ stated as follows: for all $\theta,\theta'\in\R^d$, the inequality
\begin{align}
\lara{\theta-\theta', h(\theta)-h(\theta')} \geq -\bar{L}|\theta-\theta'|^2\nonumber
\end{align}
holds, where $\bar{L}:= L(1+2\bar{R})^{r}$ and $\bar{R}:=(b/a)^{1/(r-\bar{r})}\1_{\{r>0\}}$.
\end{remark}
\textit{Proof.} See Section \hyperref[proof:remark2.5]{5.1}.\\

The final assumption is a polynomial Lipschitz continuity condition imposed on $\nabla h$, i.e. the Hessian of $U$. This additional smoothness condition is similar to Assumption H4 of \cite{brosse2019tamed} and we impose it in order to obtain improved convergence rates compared to those derived in \cite{lovas2020taming} and \cite{lim2023non} for TUSLA. 
\begin{assumption}
\label{assumption:4}
There exist constant $L_\nabla > 0$ such that for all $\theta,\theta'\in \R^d$,
\begin{align}
\lbrb{\nabla h(\theta) - \nabla h(\theta')} \leq L_\nabla (1+|\theta| + |\theta'|)^\nu|\theta-\theta'|.\nonumber
\end{align}
\end{assumption}

\begin{remark}
\label{remark:2.6}
Under Assumption \ref{assumption:4}, the following assertions hold.
\begin{enumerate}
\item[(i)] For all $\theta\in\R^d$, $|\nabla h(\theta)| \leq C_\nabla(1+|\theta|^{\nu+1})$, where $C_\nabla := 2\max\{2^{\nu-1}L_\nabla, |\nabla h(0)|\}$.
\item[(ii)] If $r=\nu+1$, Assumption \hyperref[assumption:2]{2-(a)} holds with the choice of constant $L=C_\nabla$.
\item[(iii)]For all $\theta,\theta'\in\R^d$,
\begin{align}
|h(\theta)-h(\theta')-\nabla h(\theta')(\theta-\theta')|\leq \bar{L}_\nabla(1+|\theta|^\nu + |\theta'|^\nu)|\theta-\theta'|^2,\nonumber
\end{align}
where $\bar{L}_\nabla := 3^{\nu - 1}L_\nabla$.
\end{enumerate}
\end{remark}

\textit{Proof.} See Section \hyperref[proof:remark2.6i]{5.1}.
\subsection{Main Results} \label{subsec:main_results} Define the maximum step size restriction $\tilde{\lambda}_{\max}$ of mTULA as 
\begin{align}
\tilde{\lambda}_{\max}:=\min\left\{1,\tfrac{\bar{a}^2}{8K^4}, \tfrac{1}{\bar{a}^2}\right\}.\label{rf:max_step}
\end{align}
The main results of this paper are non-asymptotic upper bound estimates in Wasserstein-1 and Wasserstein-2 distance between the law of the algorithm $\theta_n^\lambda$ defined in (\ref{eqn:tula}) and the target distribution $\pi_\beta$ defined in (\ref{eqn:target_distribution}), stated precisely in Theorems \ref{theorem:2.7} and \ref{theorem:2.8}.
\begin{theorem}
\label{theorem:2.7}
Let Assumptions \ref{assumption:1}, \ref{assumption:2}, \ref{assumption:3}, \ref{assumption:4} hold. Then, for any $\lambda\in(0,\tilde{\lambda}_{\max})$, $n\in\N_0$, the mTULA algorithm $(\theta_n^\lambda)_{n\in\N_0}$ has the following non-asymptotic upper bound estimate in Wasserstein-1 distance:
\begin{align}
W_1(\mc{L}(\theta_{n}^\lambda), \pi_\beta)
\leq &\ C_1 e^{-C_0\lambda n}\left(1+\E\lsrs{|\theta_0|^{r_*}}\right) + C_2\lambda\nonumber
\end{align}
where the constants $C_0, C_1, C_2$ are given in (\ref{eqn:theorem2.7_constants_C0_C1_C2}).
\end{theorem}


\begin{theorem}
\label{theorem:2.8}
Let Assumptions \ref{assumption:1}, \ref{assumption:2}, \ref{assumption:3}, \ref{assumption:4} hold. Then, for any $\lambda\in(0,\tilde{\lambda}_{\max})$, $n\in\N_0$, the mTULA algorithm $(\theta_n^\lambda)_{n\in\N_0}$ has the following non-asymptotic upper bound estimate in Wasserstein-2 distance:
\begin{align}
W_2(\mc{L}(\theta_n^\lambda),\pi_\beta)\leq C_4e^{-C_3\lambda n}\left(1+\E\lsrs{|\theta_0|^{r_*}}\right)^{1/2} + C_5\lambda^{1/2}\nonumber
\end{align}
where the constants $C_3, C_4, C_5$ are given in (\ref{eqn:theorem2.8_constant_C3_C4_C5}).
\end{theorem}


\begin{remark}
The first and second terms of the RHS of the upper bound estimates of Theorems \ref{theorem:2.7} and \ref{theorem:2.8}, respectively, can be interpreted as the error due to algorithm having finite time horizon $\lambda n$, and the error due to discretisation with step size $\lambda$. In particular, we obtain in Theorems \ref{theorem:2.7} and \ref{theorem:2.8} rates of convergence of $\mc{O}(\lambda)$ in Wasserstein-1 distance and $\mc{O}(\lambda^{1/2})$ in Wasserstein-2 distance for the discretisation error of mTULA.\\

We highlight that Remark~\ref{remark:2.5} and Remark~\ref{remark:2.6} are key inequalities in establishing the aforementioned improved rates of convergence of mTULA in Wasserstein-1 and Wasserstein-2 distances, which hold due to Assumption~\ref{assumption:3} and Assumption~\ref{assumption:4}, respectively. One may also refer to the proof of Lemma~\ref{lemma:4.9} for the detailed arguments. Moreover, we note that the absolute moments of mTULA are finite due to Remark~\ref{remark:2.4} (Assumption~\ref{assumption:3}), see Lemma~\ref{lemma:4.2}, which are used to obtain the convergence in Wasserstein distances.
\end{remark}

\begin{remark}
We note that the constants $C_1, C_2, C_4, C_5$ in Theorem~\ref{theorem:2.7} and~\ref{theorem:2.8} have exponential dependence on the dimension $d$ and on $\beta$ due to \cite[Theorem 2.2]{eberle2019quantitative}. A very recent result \cite[Corollary 5.7]{lytras2023taming} shows that the corresponding constants for sTULA proposed in \cite{lytras2023taming} only depend polynomially on the dimension, which is achieved by applying a Log-Sobolev inequality \cite[Theorem 5.3]{lytras2023taming}. 
However, as mTULA \eqref{eqn:tula}-\eqref{eqn:tamed_gradient} utilises a taming factor which is different from that of sTULA, we cannot apply that result directly and leave it for future research to see whether one can also obtain a polynomial dependence for the constants $C_1, C_2, C_4, C_5$ in Theorem~\ref{theorem:2.7} and~\ref{theorem:2.8} for mTULA by employing  a Log-Sobolev inequality.
\end{remark}

\subsection{Comparison with Related Works}
\label{subsec:comparison}
In this section, we compare our main results and assumptions under which they were obtained with those of \cite{brosse2019tamed} and \cite{lim2023non} where tamed variants of ULA were proposed to deal with the case of super-linearly growing gradient.\\

In \cite{brosse2019tamed}, the Tamed Unadjusted Langevin Algorithm (TULA) was proposed and non-asymptotic bounds in Wasserstein-2 distances were obtained under the two sets of assumptions (H1, H2, H3) and (A1, A2, H2, H3, and H4), in the notation of \cite{brosse2019tamed}. Assumption H1 of \cite{brosse2019tamed}  is a polynomial Lipschitz condition on the gradient which is equivalent to Assumption \ref{assumption:2} in our work. Under Assumption H2 of \cite{brosse2019tamed}, the gradient $\nabla U(\theta)$ becomes radially unbounded as $|\theta|$ diverges to infinity and the expression $\lara{\theta, \nabla U(\theta)}$ can be bounded from below by a negative constant on $\R^d$. In this sense, Assumption H2 is loosely comparable to our dissipativity condition (\ref{eqn:dissipativity}) of Remark \ref{remark:2.4} in our work from which we may deduce the same two properties for the potential $U$. Assumption A1 of \cite{brosse2019tamed} is a condition on closeness in norm between the tamed and actual gradient for sufficiently small step size $\lambda$, which is satisfied for our particular choice of taming factor; see, for example, a derivation of an upper bound on $|h(\theta)-h_\lambda(\theta)|$ for our taming factor embedded in our argument for deriving the bound (\ref{eqn:lemma4.9_I31_ub}). Assumption A2 of \cite{brosse2019tamed} which was imposed to ensure the finiteness of exponential moments of TULA may not be directly comparable to the assumptions imposed in our work where we adopt a different approach to ensure finiteness of polynomial moments of mTULA. Note that the authors of \cite{brosse2019tamed} show that for their choices of taming factors, H1 and H2 together imply A1 and A2. \\

Assumption H4 of \cite{brosse2019tamed} is an additional smoothness condition imposed on $U$ which assumes polynomial H\"{o}lder continuity on the Hessian of $U$ with exponent $\beta\in[0,1]$. When the exponent $\beta$ equals $1$, this assumption is equivalent to Assumption \ref{assumption:4} in our paper. Under this additional smoothness condition, the authors of \cite{brosse2019tamed} were able to improve the order of convergence in Wasserstein-2 distance from $1/2$ to $(1+\beta)/2$. Even though such rates of convergence are stronger than the ones we obtained in Theorem \ref{theorem:2.8}, the key difference is the condition of strong convexity imposed in \cite{brosse2019tamed} as Assumption H3 for both their estimates in Wasserstein-2 distances, whereas our results can be applied to non-convex potential functions $U$ in view of our convex at infinity condition imposed in Assumption \ref{assumption:3} of our work. For a detailed comparison between the assumption of strong convexity and our convex at infinity condition, we refer the reader to  Remark \ref{remark:2.3}.\\

In \cite{lim2023non}, the polynomial Lipschitz condition imposed in Assumptions 2 and 3 of \cite{lim2023non}, and the convex at infinity condition imposed in Assumption 4 of \cite{lim2023non} are directly analogous to Assumptions \ref{assumption:2} and \ref{assumption:3} of this paper in the case of deterministic gradient. With our choice of taming factor as $(1+\lambda|\theta|^{2r})^{1/2}$, see (\ref{eqn:tamed_gradient}), as opposed to $(1+\sqrt{\lambda}|\theta|^{r})$ in \cite{lim2023non}, as well as assuming additional smoothness conditions on $U$ in Assumption \ref{assumption:4} of our work, we were able to improve the orders of convergence of $1/2$ and $1/4$ for Wasserstein distances of orders 1 and 2, respectively, derived in \cite{lim2023non}, to $1$ and $1/2$, respectively, in Theorems \ref{theorem:2.7} and \ref{theorem:2.8}.

\section{Numerical Examples}
\label{sec:simulations}
In this section, we illustrate our theoretical results with some numerical simulations. We first present a motivating example where the ULA algorithm fails to work, highlighting numerically the importance of the taming technique and the tamed algorithms. Then, we use mTULA to sample from several high-dimensional distributions illustrating its wide applicability. Finally, we conclude this section by providing discussions for our numerical results.

\subsection{Motivating Example} 
We consider to sample from a high-dimensional double-well potential distribution denoted by $\pi^{\mathsf{dw}}_\beta$ where $U(\theta)=|\theta|^4/4- |\theta|^2/2$ and $h(\theta) = (|\theta|^2-1)\theta$ for all $\theta \in \R^d$. It can be shown that this example satisfies our Assumptions \ref{assumption:1}-\ref{assumption:4}, see Proposition \ref{proposition:3.1}. We aim to numerically obtain the absolute second moment of $\pi^{\mathsf{dw}}_\beta$ as its explicit value is not available. We note that the absolute second moment of $\pi^{\mathsf{dw}}_\beta$ is finite due to \cite[Lemma A.1]{lim2023non} and \cite[Proposition 1-(ii)]{durmus2019high} under our assumptions. To this end, we run ULA \cite[Eq. (2)]{brosse2019tamed}, TULA \cite[Eq. (3)]{brosse2019tamed}, mTULA \eqref{eqn:tula}-\eqref{eqn:tamed_gradient} with step size $\lambda = 0.001$ to sample from $\pi^{\mathsf{dw}}_\beta$ and calculate the absolute second moment of the approximations. We set $\beta= 1$ and $d=100$, $\theta_0=10\in\R^d$, the number of iterations is $4\times10^5$, and the number of independent Markov chains for each algorithm is $250$. The experiments are conducted for five times and the numerical results are summarised in Table~\ref{table:1}. The reference value, which equals to $10.785$, is obtained by runing mTULA for $n = 4\times10^{7}$ iterations. 
\\

\begin{table}[h]
\centering
\begin{tabular}{|c c c c|} 
 \hline
 &$\E\left[|\theta_n^\lambda|^2\right]$ & $\E\left[|\theta_n^{\mathsf{TULA}}|^2\right]$ &$\E\left[|\theta_n^{\mathsf{ULA}}|^2\right]$\\ [1ex] 
 \hline
Simulation 1 &10.743&10.503&NaN\\
Simulation 2 &10.784&10.784&NaN\\
Simulation 3 &10.703&10.612&NaN\\
Simulation 4 &10.805&10.721&NaN\\
Simulation 5 &10.721&10.614&NaN\\ [1ex] 
 \hline
\end{tabular}
\caption{Absolute second moment of approximations obtained using mTULA, TULA, and ULA, denoted by $\E\left[|\theta_n^\lambda|^2\right]$, $\E\left[|\theta_n^{\mathsf{TULA}}|^2\right]$, and $\E\left[|\theta_n^{\mathsf{ULA}}|^2\right]$, respectively.}
\label{table:1}
\end{table}
Table \ref{table:1} illustrates that the absolute second moment of the ULA algorithm diverges in all the experiments while mTULA and TULA produce finite values which are close to the reference value.

\subsection{Sampling using mTULA}
By running mTULA, we drew samples from three high-dimensional target distributions, where the choices of target distribution $\dee\pi_\beta(\theta)\propto e^{-U(\theta)}\ \dee\theta$ considered were:


\begin{enumerate}
\item[(i)] \textit{Multivariate Standard Gaussian $\mc{N}(0,I_d)$}, with potential 
\begin{align}
U(\theta)=\tfrac{1}{2}|\theta|^2,\quad \theta\in\R^d.\label{eqn:gaussian_u}
\end{align}\\[-20pt]
\item[(ii)] \textit{Multivariate Gaussian Mixture}, with potential
\begin{align}
U(\theta)= \tfrac{1}{2}|\theta-\dot{a}|^2 - \log(1+e^{-2\lara{\dot{a},\theta}}),\quad \theta\in\R^d.\label{eqn:gaussian_mixture_u}
\end{align}
for a given $\dot{a}\in\R^d$.\\
\item[(iii)] \textit{Double-well Potential}, with potential
\begin{align}
U(\theta)=\tfrac{1}{4}|\theta|^4-\tfrac{1}{2}|\theta|^2,\quad \theta\in\R^d.\label{eqn:doublewell_u}
\end{align}
\end{enumerate}

The next proposition states that our assumptions are satisfied for these three choices of target distributions. Therefore, our theoretical results apply to samples from these target distributions drawn with mTULA.

\begin{proposition}
\label{proposition:3.1} 
Assumptions \ref{assumption:2}, \ref{assumption:3}, \ref{assumption:4} are satisfied for the choices of potential functions $U$ given in equations (\ref{eqn:gaussian_u}), (\ref{eqn:gaussian_mixture_u}), and (\ref{eqn:doublewell_u}). For each $U$, choices of constants $r$, $\nu$, $L$, $K$, $a$, $b$, $\bar{r}$, $\tilde{a}$, $\tilde{b}$, $L_\nabla$ for which the assumptions are satisfied are respectively given by
\begin{enumerate}
\item[(i)] \textit{Multivariate Standard Gaussian $\mc{N}(0,I_d)$}: $r=0$, $\nu=0$, $L=1$, $K=1$, $\tilde{a}=1$, $\tilde{b}=1$, $L_\nabla=1$.
\item[(ii)] \textit{Multivariate Gaussian Mixture}: $r=0$, $\nu=0$, $L=1+4|\dot{a}|^2$, $K=\max\{1,|\dot{a}|\}$, $\tilde{a}=1/2$, $\tilde{b}=2$, $L_\nabla=8|\dot{a}|^3$.
\item[(iii)] \textit{Double-well potential}: $r=2$, $\nu=1$, $L=1$, $K=2$, $a=1/2$, $b=1$, $\bar{r}=0$, $L_\nabla=3$.
\end{enumerate}
\end{proposition}

\textit{Proof}. See Section \hyperref[proof:proposition3.1]{5.2}.\\

In particular, we note that the multivariate standard Gaussian and Gaussian mixture distributions have potentials which are globally Lipschitz continuous, with that of the multivariate standard Gaussian distribution also satisfying strong convexity. In contrast, the double-well potential is representative of the typical non-convex and 'super-linear growth' setting under which our main convergence results in Theorems \ref{theorem:2.7} and \ref{theorem:2.8} were derived -- it is both non-convex and has gradient which is locally, but not globally, Lipschitz continuous.\\

In our numerical experiments, we set $\beta=1$ and the dimension as $d=100$. To obtain samples from each target distribution $\pi_\beta$, we used deterministic initialisation $\theta_0=0\in\R^d$, ensuring in particular that Assumption \ref{assumption:1} satisfied, and ran $m=250$ independent mTULA Markov chains for each step size $\lambda\in \{0.001, 0.005, 0.01, 0.025, 0.05, 0.1\}$, holding the time horizon of each chain constant at $\lambda n = 400$, with $n$ denoting the length of the chain. For the Gaussian mixture model, $\dot{a}\in\R^d$ was chosen such that all its components are equal and $|\dot{a}|=2$, so that its potential $U$ is not strongly convex \cite{dalalyan2017theoretical}. Normalised histogram plots were then generated from these samples drawn with mTULA to illustrate the convergence of the algorithm, as detailed in the next section.

\label{subsec:histogram} For each target distribution $\pi_\beta$ and step size $\lambda$, we plotted normalised histograms of the first components of samples obtained from the last iterations of each of the $m=250$ independent mTULA Markov chains. To visualise the closeness of the empirical distributions of the first components of the drawn samples to the marginal distributions of the respective target distributions, we superimposed the theoretical probability densities of the marginal distributions of the first components of the respective target distributions onto the normalised histograms. We note that each component of a sample drawn from a multivariate standard Gaussian distribution is distributed as $\mc{N}(0,1)$, while the densities of the first components of samples drawn from a multivariate Gaussian mixture model and double-well potential are, respectively, given by
\begin{align}
\R\ni\theta^{(1)}\mapsto \frac{1}{2\sqrt{2\pi}}\lsrs{e^{-\left(\theta^{(1)}-a^{(1)}\right)^2/2}+e^{-\left(\theta^{(1)}+a^{(1)}\right)^2/2}}\
\end{align}
and 
\begin{align}
\R\ni \theta^{(1)}\mapsto \frac{\Gamma\left(\tfrac{d}{2}\right)}{\sqrt{\pi}\Gamma\left(\tfrac{d-1}{2}\right)}\cdot\frac{\int^\infty_0 r^{(d-3)/2}\exp\left\{-\tfrac{1}{4}\left(r+\theta^{(1)}\right)^2+\tfrac{1}{2}\left(r+\theta^{(1)}\right)\right\}\ \dee r}{\int^\infty_0 r^{d/2-1}\exp\left\{-\tfrac{1}{4}r^2+\tfrac{1}{2}r\right\}\ \dee r}.
\end{align}

\begin{center}
\begin{figure}[H]
\includegraphics[scale=0.55]{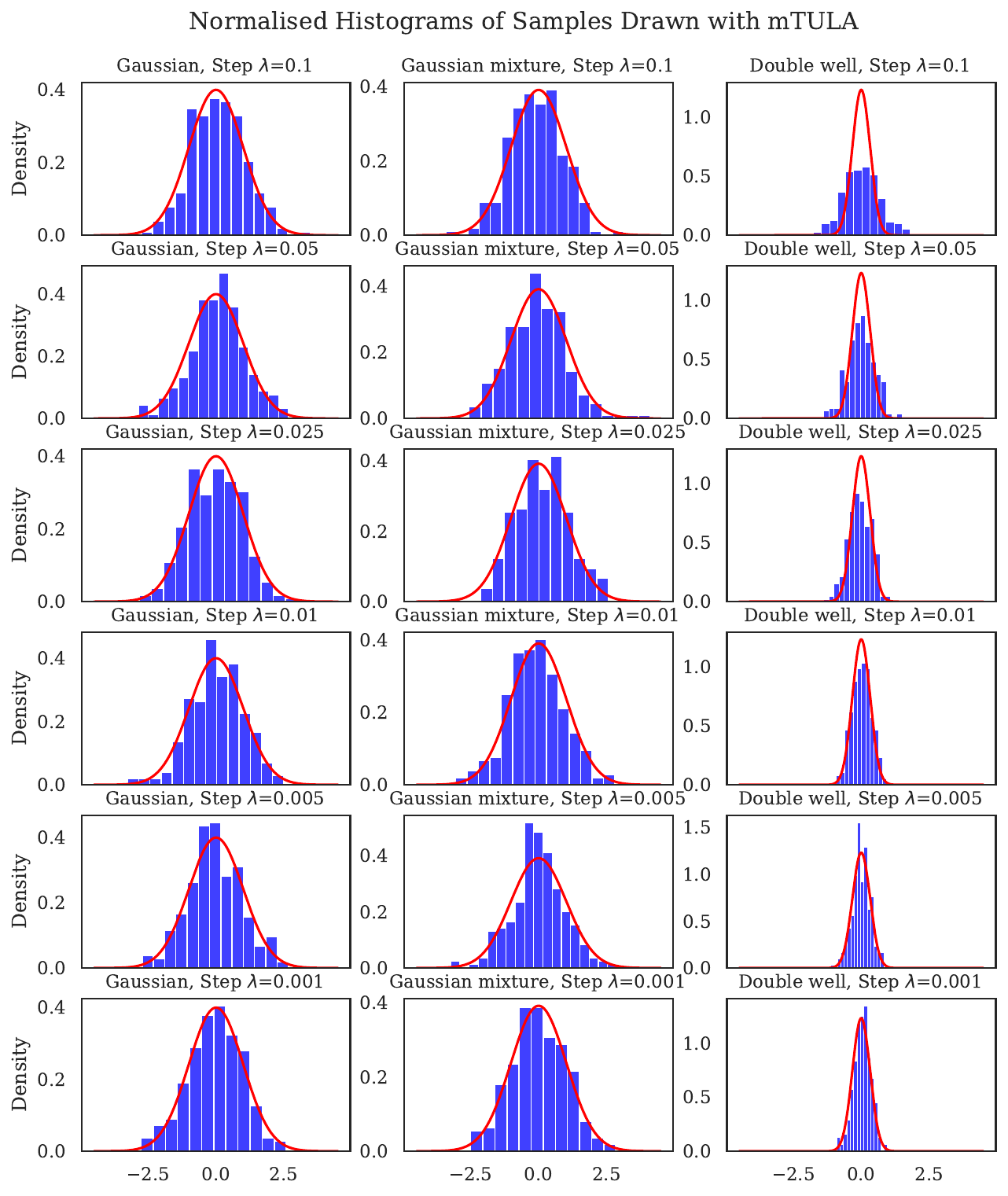}
\caption{Normalised histograms of first components of samples drawn with mTULA}
\label{figure:histogram}
\end{figure}
\end{center}

Here, we denote $\R^d\ni \theta=(\theta^{(1)},\cdots,\theta^{(d)})$, i.e. $\theta^{(1)}$ and $a^{(1)}$ denote, respectively, the first components of $\theta, a\in\R^d$, whereas $\Gamma$ denotes the Euler Gamma function. Figure \hyperref[figure:histogram]{1} displays, for each choice of target distribution $\pi_\beta$ and step size $\lambda$ we considered, the normalised histogram plots of samples generated with mTULA together with the superimposed theoretical marginal probability density curves of the corresponding target distributions. The source code for our numerical experiments are available at the following GitHub repository: \url{https://github.com/tracyyingzhang/mTULA}.

\subsection{Discussion} 
In general, the numerical experiments appear to support our theoretical results. One may observe from the histogram plots that the samples generated by mTULA were close to their corresponding target distributions for sufficiently small step sizes. In particular, what is markedly observable for the double-well potential, which has a super-linearly growing gradient $h$ and non-convex potential $U$, is that the normalised histograms of the generated samples from mTULA approached the exact probability density curve of the target distribution as the step size $\lambda$ decreased. This observation supports our main convergence results in Theorems \ref{theorem:2.7} and \ref{theorem:2.8} which were derived under a super-linear and non-convex setting.

\section{Proof Overview of the Main Results}
\label{sec:proofs}
This section presents our methodical approach for establishing the convergence rate of mTULA in Wasserstein-1 and Wasserstein-2 distances stated in Theorems \ref{theorem:2.7} and \ref{theorem:2.8}. A key ingredient in our proof is the introduction of an auxiliary process $\bar{\zeta}_t^{\lambda, n}$, given explicitly in Definition \ref{definition:4.1}, to ease computations. We obtain some preliminary estimates on moments of functions of the continuous-time interpolated algorithm and of the auxiliary process, which are subsequently used as intermediate results for arriving at the convergence rates of mTULA.

\subsection{Auxiliary Processes} Consider the $\R^d$-valued Langevin SDE $(Z_t)_{t\in\R_+}$ given by
\begin{align}
\dee Z_t = -h(Z_t)\ \dee t + \sqrt{2\beta^{-1}}\ \dee B_t\label{eqn:langevin}
\end{align}
with $Z_0:=\theta_0$, where $(B_t)_{t\geq 0}$ is a standard $d$-dimensional Brownian motion on $(\Omega, \mathcal{F}, \mathbb{P})$. Denote by $(\mathcal{F}_t)_{t\geq 0}$ the $\mathbb{P}$-completion of the natural filtration generated by $(B_t)_{t\geq 0}$, assumed to be independent of $\theta_0$.\\

For each $\lambda>0$, denote by $Z_t^\lambda:=Z_{\lambda t},t\in\R_+$, the time-changed Langevin SDE given by
\begin{align}
\dee Z_t^\lambda = -\lambda h(Z_t^\lambda)\ \dee t + \sqrt{2\lambda\beta^{-1}}\ \dee B_t^\lambda \label{eqn:langevin_time_changed}
\end{align}
with the initial condition $Z_0^\lambda:=\theta_0$, where $B_t^\lambda:=B_{\lambda t}/\sqrt{\lambda}$, $t\geq 0$. Note that $(B_t^\lambda)_{t\geq 0}$ is a $d$-dimensional standard Brownian motion. For each $\lambda >0$, denote by  $(\mathcal{F}_t^\lambda)_{t_\geq 0}$ the natural filtration of $(B_t^\lambda)_{t\geq 0}$ with $\mathcal{F}_t^\lambda:=\mathcal{F}_{\lambda t}$ for each $t\geq 0$, which is also independent of $\theta_0$.\\ 

Then, define the continuous-time interpolation of mTULA (\ref{eqn:tula}), denoted by $(\bar{\theta}_t^\lambda)_{t\in\R_+}$, as
\begin{align}
\dee \bar{\theta}_t^\lambda = -\lambda h_\lambda(\bar{\theta}_{\lfrf{t}}^\lambda)\ \dee t + \sqrt{2\lambda \beta^{-1}} \dee B_t^\lambda\label{eqn:tula_interpolated}
\end{align}
with the initial condition $\bar{\theta}_0^\lambda:=\theta_0$. By construction, $\mc{L}\left(\bar{\theta}_n^\lambda\right)=\mc{L}\left(\theta_n^\lambda\right)$ for $n\in\N_0$. That is, the law of the interpolated process coincides with that of mTULA (\ref{eqn:tula}) at all integer times.\\

Moreover, denote by $\zeta_t^{s,v,\lambda}$, for $0\leq s < t <\infty$, a continuous-time process defined by the SDE
\begin{align}
\dee \zeta_t^{s,v,\lambda} = -\lambda h(\zeta_{t}^{s,v,\lambda})\ \dee t + \sqrt{2\lambda\beta^{-1}}\dee B_t^\lambda\label{eqn:auxiliary_process}
\end{align}

with initial condition $\zeta_s^{s,v,\lambda}=v\in\R^d$.

\begin{definition}
\label{definition:4.1}
For each fixed $\lambda > 0$ and $n\in\N_0$, define $\bar{\zeta}_t^{\lambda,n}:=\zeta_t^{nT, \bar{\theta}_{nT}^\lambda, \lambda}$, $t\geq nT$, where $T\equiv T(\lambda):=\lfrf{1/\lambda}$.
\end{definition}

\subsection{Preliminary Estimates} We first establish some preliminary bounds necessary for establishing our main convergence results. For each $p\in\N_0$, define the Lyapunov function $V_p(\theta):=(1+|\theta|^2)^{p/2}$ for all $\theta\in\R^d$. Similarly, define $v_p(w):=(1+w^2)^{p/2}$ for all $w\geq 0$. Denote by $\mc{P}_{V_p}(\R^d)$ the set of probability measures $\mu\in\mc{P}(\R^d)$ such that $\int_{\R^d}V_p(\theta)\ \dee\mu(\theta)<\infty$.\\

The following lemma provides moment estimates for $(\bar{\theta}_t^\lambda)_{t\geq 0}$ defined in (\ref{eqn:tula_interpolated}).

\begin{lemma}
\label{lemma:4.2}
Let Assumptions \ref{assumption:1}, \ref{assumption:2}, \ref{assumption:3} hold. Then, the following holds:
\begin{enumerate}
\item[(i)] For any $\lambda\in(0,\lambda_{1,\max})$, $n\in\N_0$ and $t\in(n, n+1]$,
\begin{align}
\E\lsrs{|\bar{\theta}^\lambda_t|^2} \leq (1 - \lambda(t-n)\bar{a}\kappa)(1-\lambda\bar{a}\kappa)^n\E\lsrs{|\bar{\theta}_0^\lambda|^2} + c_0\left(1+\frac{1}{\bar{a}\kappa}\right),\nonumber
\end{align}
where $\lambda_{1,\max}:= \min\{1,\frac{\bar{a}^2}{8K^4}\}$ and the constants $c_0,\kappa$ are defined as $c_0 := \bar{a}\kappa+2\bar{b}+2d\beta^{-1}+2K^2$ and $\kappa := \tfrac{1}{\sqrt{2}}$.
\item[(ii)] For any $p\in[2,\infty)\cap\N$, $\lambda \in (0,\tilde{\lambda}_{\max})$, $n\in\N_0$, and $t\in (n,n+1]$,
\begin{align}
\E\lsrs{|\bar{\theta}_t^\lambda|^{2p}}\leq \left(1-\lambda(t-n)\frac{\bar{a}\tilde{\kappa}(2)}{2}\right)\left(1-\lambda\frac{\bar{a}\tilde{\kappa}(2)}{2}\right)^n\E\lsrs{|\bar{\theta}_0^\lambda|^{2p}}+c_3(p)\left(1+\frac{2}{\bar{a}\tilde{\kappa}(2)}\right),\nonumber
\end{align}
where $\tilde{\lambda}_{\max}:= \min\{1,\frac{\bar{a}^2}{8K^4},\frac{1}{\bar{a}^2}\}$ and the constants $c_3(p), \tilde{\kappa}(2)$ are given explicitly in (\ref{eqn:lemma4.2_constant_kappatilde_M1_c3_M2_c2_c1}).
\end{enumerate}
\end{lemma}
\textit{Proof.} See Section \hyperref[proof:lemma4.2]{5.3}.\\

Combining the cases $p=1$ and $p\in [2,\infty)\cap \N$ and applying, for each fixed $\lambda$, the inequalities $1-x \leq e^{-x}$, for all $x\in \R$, and $\frac{1}{2}\leq \lambda T\leq 1$  yields the following corollary from the moment estimates:

\begin{corollary}
\label{corollary:4.3}
Let Assumptions \ref{assumption:1}, \ref{assumption:2}, \ref{assumption:3} hold. Then, for any $p\in \N$, $\lambda \in (0,\tilde{\lambda}_{\max})$, $n\in \N_0$, and $s\geq nT$,
\begin{align}
\E\lsrs{|\bar{\theta}_s^\lambda|^{2p}}\leq e^{-\lambda\bar{a}\kappa_* \lfrf{s}}\E\lsrs{|\bar{\theta}_0^\lambda|^{2p}}+c_*(p)\left(1+\frac{1}{\bar{a}\kappa_*}\right)\leq e^{-\bar{a}\kappa_* n/2}\E\lsrs{|\bar{\theta}_0^\lambda|^{2p}}+c_*(p)\left(1+\frac{1}{\bar{a}\kappa_*}\right),\nonumber
\end{align}
where $\kappa_*:=\min\{\kappa, \frac{\tilde{\kappa}(2)}{2}\}$ and $c_*(p):=\max\{c_0,c_3(p)\1_{\{p\geq 2\}}\}$ for $p\in \N$. The constants $\kappa, c_0, c_3(p), \tilde{\kappa}(2)$ are given in Lemma \ref{lemma:4.2} and (\ref{eqn:lemma4.2_constant_kappatilde_M1_c3_M2_c2_c1}).\\
\end{corollary}

In addition, a drift condition is obtained for the function $V_p$, which precise statement is given below.

\begin{lemma}
\label{lemma:4.4}
Let Assumptions \ref{assumption:1}, \ref{assumption:2}, \ref{assumption:3} hold. Then, for any $p\in [2,\infty)\cap\N$ and $\theta\in\R^d$, one obtains
\begin{align}
\beta^{-1}\Delta V_p(\theta) - \lara{\nabla V_p(\theta), h(\theta)}\leq -c_{V,1}(p)V_p(\theta) + c_{V,2}(p),\nonumber
\end{align}
where, for every $p\in\N_0$, the constants $c_{V,1}(p)$ and $c_{V,2}(p)$ are defined as $c_{V,1}(p):= \frac{\bar{a} p}{2}$ and $c_{V,2}(p):= \frac{\bar{a}p}{2}v_p(M_V(p))$, with $M_V(p):=\left(1+\frac{2\bar{b}'+2\beta^{-1}(d+p-2)}{\bar{a}}\right)^{1/2}$. Here, $\Delta$ denotes the Laplace operator (i.e. the trace of the Hessian matrix).
\end{lemma}

\textit{Proof.} See Section \hyperref[proof:lemma4.4]{5.3}.\\

From this drift condition, one obtains the following relation between the moments of $\bar{\theta}_t^{\lambda}$ and that of the auxiliary process $\bar{\zeta}_t^{\lambda,n}$. The precise statement is given in the following lemma.
\begin{lemma}
\label{lemma:4.5}
Let Assumptions \ref{assumption:1}, \ref{assumption:2}, \ref{assumption:3} hold. Then, for any $p\in [2,\infty)\cap \N$, $\lambda \in (0,\tilde{\lambda}_{\max})$, $n\in\N$, and $t\in (nT, (n+1)T]$, the following relation holds:
\begin{align}
\E\lsrs{V_p(\bar{\zeta}_t^{\lambda,n})}
\leq e^{-\lambda \bar{a}p(t-nT)/2}\E\lsrs{V_p(\bar{\theta}_{nT}^\lambda)}+v_p(M_V(p))(1-e^{-\lambda \bar{a}p(t-nT)/2}),\nonumber
\end{align}
where the constant $M_V(p)$ is defined as in Lemma \ref{lemma:4.4}.
\end{lemma}

\textit{Proof.} See Section \hyperref[proof:lemma4.5]{5.3}.\\

Combining Corollary \ref{corollary:4.3} and Lemma \ref{lemma:4.5} yields then the following moment estimate for the auxiliary process $\bar{\zeta}_t^{\lambda,n}$.

\begin{corollary}
\label{corollary:4.6}
Let Assumptions \ref{assumption:1}, \ref{assumption:2}, \ref{assumption:3} hold. Then, for any $p\in \N$, $\lambda\in (0, \tilde{\lambda}_{\max})$, $n\in\N_0$, and $t\in (nT, (n+1)T]$, 
\begin{align}
\E\lsrs{V_{2p}(\bar{\zeta}_t^{\lambda,n})}\leq 2^{p-1} + 2^{p-1}e^{-\lambda \bar{a}\min\{p,\kappa_*\}t}\E\lsrs{|\bar{\theta}_0^{\lambda}|^{2p}} + 2^{p-1} c_*(p)\left(1+\frac{1}{\bar{a}\kappa_*}\right) + v_{2p}(M_V(2p)).\nonumber
\end{align}
\end{corollary}

\textit{Proof.} See Section \hyperref[proof:corollary4.6]{5.3}.\\
\begin{remark}
\label{remark:4.7}
In the case $p=0$, one sees that the statements of Corollary \ref{corollary:4.3} and Corollary \ref{corollary:4.6} still hold if we define $c_*(0):= 1$.\\
\end{remark}

Our final preliminary estimates required for proving the main results are some moment bounds on the one-step errors $(\bar{\theta}_t^{\lambda}-\bar{\theta}_{\lfrf{t}}^{\lambda})$ and $(\bar{\zeta}_t^{\lambda,n}-\bar{\zeta}_{\lfrf{t}}^{\lambda, n})$ of the processes $(\bar{\theta}_t^\lambda))_{t\geq 0}$ and $(\bar{\zeta}_t^{\lambda,n})_{t\geq 0}$.
\begin{lemma}
\label{lemma:4.8}
Let Assumptions \ref{assumption:1}, \ref{assumption:2}, \ref{assumption:3} hold. Then, for any $p\in[2,\infty)\cap\N$, $\lambda\in (0,\tilde{\lambda}_{\max})$, $n\in\N_0$, and $t\in(nT,(n+1)T]$, one has the estimates
\begin{align}
\E\lsrs{\lbrb{\bar{\theta}_t^\lambda - \bar{\theta}_{\lfrf{t}}^\lambda}^8} &\leq \lambda^4\left(e^{-\bar{a}\left(\tilde{\kappa}(4r+4)\right)n/2}\bar{C}_{1,1}\E\lsrs{|\bar{\theta}_0^{\lambda}|^{8r+8}}+\bar{C}_{2,1}\right),\nonumber\\
\E\lsrs{\lbrb{\bar{\zeta}_t^{\lambda,n}-\bar{\zeta}_{\lfrf{t}}^{\lambda,n}}^4} &\leq \lambda^2\left(e^{- \bar{a}\min\{r+1, \tilde{\kappa}(2r+2)/4\}n}\bar{C}_{1,2}\E\lsrs{|\bar{\theta}_0^{\lambda}|^{4r+4}}+\bar{C}_{2,2}\right),\nonumber
\end{align}
where the constants $\bar{C}_{1,1}, \bar{C}_{2,1}, \bar{C}_{1,2}, \bar{C}_{2,2}$ are given in (\ref{eqn:lemma4.8_constant_C11bar_C21bar}) and (\ref{eqn:lemma4.8_constant_C12bar_C22bar}).
\end{lemma}

\textit{Proof.} See Section \hyperref[proof:lemma4.8]{5.3}.\\

\subsection{Proof of the main theorems}\label{sec:mtproof}
This section presents the main steps of proving our main results using the intermediate results we have derived so far. We provide an overview for the proof of Theorem \ref{theorem:2.7}, and similar arguments can be applied to obtain the result in Theorem \ref{theorem:2.8}. To establish a non-asymptotic estimate in Wasserstein-1 distance between $\mc{L}(\theta_n^\lambda)$ and $\pi_\beta$, we consider the following splitting using the continuous-time interpolation of the mTULA algorithm: for any $n\in\N_0$ and $t\in (nT, (n+1)T]$,
\begin{align}
W_1(\mc{L}(\bar{\theta}_t^\lambda),\pi_\beta)&\leq W_1(\mc{L}(\bar{\theta}_t^\lambda),\mc{L}(\bar{\zeta}_t^{\lambda,n})) + W_1(\mc{L}(\bar{\zeta}_t^{\lambda,n}),\mc{L}(Z_t^\lambda)) + W_1(\mc{L}(Z_t^\lambda),\pi_\beta).\label{eqn:wassertein_1_splitting}
\end{align}
To obtain an estimate for the first term on the RHS of \eqref{eqn:wassertein_1_splitting}, we first consider the corresponding $L^2$-distance $\E\lsrs{\lbrb{\bar{\zeta}_t^{\lambda, n} - \bar{\theta}_t^\lambda}^2}$ and apply the synchronous coupling of $\bar{\zeta}_t^{\lambda, n}$ and $\bar{\theta}_t^\lambda$. As $\bar{\theta}_t^\lambda$ can be viewed as (a continuous-time version of) the Milstein scheme of $\bar{\zeta}_t^{\lambda, n}$, its rate of convergence in $L^2$ is 1, see, e.g., \cite{kloeden2013numerical}. We note that a key result to achieve this is Remark \ref{remark:2.6}, which holds due to Assumption \ref{assumption:4}. Finally, by applying standard techniques in numerical analysis and by applying Gronwall's lemma, we obtain a desired upper estimate.\\

To upper bound the second term on the RHS of \eqref{eqn:wassertein_1_splitting}, we use the definition of $\bar{\zeta}_t^{\lambda, n}$ given in Definition \ref{definition:4.1} and view $\mc{L}(\bar{\zeta}_t^{\lambda,n})$ and $\mc{L}(Z_t^\lambda)$ as the laws of the time-changed Langevin process starting from different initial points. We then apply a contraction result in $w_{1,2}$ (defined in \eqref{eqn:wasserstein_functional} below), which is obtained by applying \cite[Theorem 2.2]{eberle2019quantitative}. Using the fact that $w_{1,2}$ upper bounds $W_1$ yields the desired result.\\

We note that the last term on the RHS of \eqref{eqn:wassertein_1_splitting} can be upper bounded by using the same approach as that for bounding the second term (as described above). This is due to the fact that $\mc{L}(Z_t^\lambda)$ and $\pi_\beta$ can also be viewed as the laws of the Langevin process starting from different initial points, i.e., from $\theta_0$ and a random variable distributed according to the invariant measure $\pi_\beta$, respectively.\\

In the following results, we provide non-asymptotic estimates with explicit constants for each of the terms on the RHS of \eqref{eqn:wassertein_1_splitting}.\\

We start by providing an upper bound for the first term on the RHS of (\ref{eqn:wassertein_1_splitting}).

\begin{lemma}
\label{lemma:4.9}
Let Assumptions \ref{assumption:1}, \ref{assumption:2}, \ref{assumption:3}, \ref{assumption:4} hold. Then, for any $\lambda \in (0,\tilde{\lambda}_{\max})$, $n\in\N_0$, and $t\in (nT, (n+1)T]$, the following inequality holds:
\begin{align}
W_2(\mc{L}(\bar{\theta}_t^\lambda), \mc{L}(\bar{\zeta}_t^{\lambda,n}))\leq \lambda\left(e^{-\bar{a}\min\{r,\kappa_*/2\}n}\bar{C}_0\E\lsrs{|\bar{\theta}_0^\lambda|^{r_*}} +\bar{C}_1\right)^{1/2}\nonumber
\end{align}
where the constant $\kappa_*$ is given in Corollary \ref{corollary:4.3} and the constants $\bar{C}_0$, $\bar{C}_1$ are given explicitly in (\ref{eqn:lemma4.9_constant_C0_C1}).
\end{lemma}

\textit{Proof.} See Section \hyperref[proof:lemma4.9]{5.3}.\\

To obtain an upper bound for the second term of the RHS of (\ref{eqn:wassertein_1_splitting}), we define, for every $p\geq 1$, $\mu,\mu'\in\mc{P}_{V_p}(\R^d)$, the functional
\begin{align}
w_{1,p}(\mu,\mu'):=\inf_{\psi\in\mc{C}(\mu,\mu')}\int_{\R^d\times\R^d}(1\wedge |\theta-\theta'|)(1+V_p(\theta)+V_p(\theta'))\ \dee \psi(\theta,\theta')\label{eqn:wasserstein_functional}
\end{align}

which, in the particular case $p=2$, is related to the Wasserstein-1 and Wasserstein-2 metric through the following inequalities.
\begin{lemma}
\label{lemma:4.10}
For any $\mu,\mu'\in \mc{P}_{V_2}(\R^d)$, the following inequalities hold:
\begin{align}
W_1(\mu,\mu')\leq w_{1,2}(\mu, \mu'),\quad W_2(\mu,\mu')\leq \sqrt{2 w_{1,2}(\mu , \mu')}.\nonumber
\end{align}
\end{lemma}

\textit{Proof.} See Section \hyperref[proof:lemma4.10]{5.3}.\\

In particular, the following contraction property satisfied by $w_{1,2}$ is instrumental for establishing upper bounds for the remaining two terms $W_1(\mc{L}(\bar{\zeta}_t^{\lambda,n}),\mc{L}(Z_t^\lambda))$ and $W_1(\mc{L}(Z_t^\lambda),\pi_\beta)$ of (\ref{eqn:wassertein_1_splitting}). This result is due to Theorem 2.2 of \cite{eberle2019quantitative}.

\begin{proposition}
\label{proposition:4:11}
Let Assumptions \ref{assumption:1}, \ref{assumption:2}, \ref{assumption:3} hold. Let $Z_t$, $Z_t'$, $t\geq 0$, be solutions of the Langevin SDE (\ref{eqn:langevin}) with square-integrable initial conditions $Z_0:=\theta_0$, $Z_0':=\theta_0'$. Then, there exist constants $\hat{c}, \dot{c}>0$ such that for all $t\geq 0$,
\begin{align}
w_{1,2}\left(\mc{L}(Z_t),\mc{L}(Z_t')\right)\leq \hat{c}e^{-\dot{c}t}w_{1,2}(\mc{L}(\theta_0), \mc{L}(\theta_0')).\nonumber
\end{align}
The explicit form of the constants $\hat{c}$ and $\dot{c}$ are given by
\begin{align}
\hat{c}:=&\ 2(1+\overline{R}_2)\exp\{\beta\bar{L}\overline{R}_2^2/8+2\overline{R}_2\}/\epsilon,\nonumber\\
\dot{c}:=&\ \min\left\{\left(\overline{R}_2\sqrt{8\pi\beta/\bar{L}}\exp\left\{\left(\overline{R}_2\sqrt{\beta\bar{L}/8}+\sqrt{8/(\beta\bar{L})}\right)^2\right\}\right)^{-1}, c_{V,1}(2)/2,\ 2c_{V,2}(2)\epsilon c_{V,1}(2)\right\},\nonumber
\end{align}
where $\overline{R}_2:=2(4c_{V,2}(2)(1+c_{V,1}(2))/c_{V,1}(2)-1)^{1/2}$, $\epsilon>0$ is chosen to satisfy the inequality
\begin{align}
\epsilon \leq 1\wedge \left(4c_{V,2}(2)\sqrt{2\pi\beta/\bar{L}}\int^{\overline{R}_1}_0\exp\left\{\left(s\sqrt{\beta\bar{L}/8} +\sqrt{8/(\beta\bar{L})}\right)^2\right\}\ \dee s\right)^{-1},\nonumber
\end{align}
with $\overline{R}_1 := 2(2c_{V,2}(2)/c_{V,1}(2)-1)^{1/2}$, and the constants $c_{V,1}(2), c_{V,2}(2)$ are defined in Lemma \ref{lemma:4.4}.
\end{proposition}

\textit{Proof.} See Section \hyperref[proof:proposition4.11]{5.3}.\\

This contraction property satisfied by $w_{1,2}$ provides upper estimates for the second and third terms of the RHS of (\ref{eqn:wassertein_1_splitting}), stated in the following lemmas.

\begin{lemma}
\label{lemma:4.12}
Let Assumptions \ref{assumption:1}, \ref{assumption:2}, \ref{assumption:3}, \ref{assumption:4} hold. Then, for any $\lambda \in (0,\tilde{\lambda}_{\max})$, $n\in\N_0$, and $t\in(nT, (n+1)T]$, the following inequality holds:
\begin{align}
W_1(\mc{L}(\bar{\zeta}_t^{\lambda,n}),\mc{L}(Z_t^\lambda))\leq \lambda\left(e^{-\min\{\dot{c}/4, \bar{a}/2, \bar{a}r/2,\bar{a}\kappa_*/4\}n}\bar{C}_2\E\lsrs{|\bar{\theta}_0^\lambda|^{r_*}} + \bar{C}_3\right), \nonumber
\end{align}
where the constant $\kappa_*$ is given in Corollary \ref{corollary:4.3}, the constants $\dot{c}$, $\hat{c}$ are given in Proposition \ref{proposition:4:11}, and $\bar{C}_2$, $\bar{C}_3$ are given in (\ref{eqn:lemma4.12_constant_C2_C3}).
\end{lemma}

\textit{Proof.} See Section \hyperref[proof:lemma4.12]{5.3}.\\

Up to this point, we have developed sufficient machinery to establish our main results. Their proofs are given as follows.
\subsection*{\textit{Proof of Theorem \ref{theorem:2.7}}}
\label{proof:theorem2.7}
Let $\lambda\in(0,\tilde{\lambda}_{\max})$, $n\in\N_0$, and let $t\in(nT, (n+1)T]$. Since $\pi_\beta$ is the invariant measure for the time-changed Langevin SDE (\ref{eqn:langevin_time_changed}), direct application of Lemma \ref{lemma:4.10} and Proposition \ref{proposition:4:11} yields 
\begin{align}
W_1(\mc{L}(Z_t^\lambda),\pi_\beta)
\leq &\  w_{1,2}(\mc{L}(Z_t^\lambda),\pi_\beta)\nonumber\\ 
\leq &\ \hat{c}e^{-\dot{c}\lambda t} w_{1,2}(\theta_0,\pi_\beta)\nonumber\\
\leq &\ \hat{c}e^{-\dot{c}\lambda t}\lsrs{1+\E[V_2(\theta_0)] + \int_{\R^d} V_2(\theta)\ \dee\pi_\beta(\theta)}.\label{eqn:theorem2.7_wass1_ub}
\end{align}
Substituting (\ref{eqn:theorem2.7_wass1_ub}) together with the results of Lemmas \ref{lemma:4.9} and \ref{lemma:4.12} into (\ref{eqn:wassertein_1_splitting}) then yields
\begin{align}
W_1(\mc{L}(\bar{\theta}_t^\lambda),\pi_\beta)
\leq&\ \lambda\left(e^{-\bar{a}\min\{r,\kappa_*/2\}n}\bar{C}_0\E\lsrs{|\bar{\theta}_0^\lambda|^{r_*}} +\bar{C}_1\right)^{1/2}\nonumber\\
&\ + \lambda\left(e^{-\min\{\dot{c}/4, \bar{a}/2, \bar{a}r/2,\bar{a}\kappa_*/4\}n}\bar{C}_2\E\lsrs{|\bar{\theta}_0^\lambda|^{r_*}} + \bar{C}_3\right) \nonumber\\
&\ + \hat{c}e^{-\dot{c}\lambda t}\lsrs{1+\E[V_2(\theta_0)] + \int_{\R^d} V_2(\theta)\ \dee\pi_\beta(\theta)}\nonumber\\
\leq&\ \lambda\left(e^{-\min\{\dot{c}/4, \bar{a}/2, \bar{a}r/2,\bar{a}\kappa_*/4\}n}\bar{C}_0^{1/2}\left(1+\E\lsrs{|\bar{\theta}_0^\lambda|^{r_*}}\right) +\bar{C}_1^{1/2}\right)\nonumber\\
&\ + \lambda\left(e^{-\min\{\dot{c}/4, \bar{a}/2, \bar{a}r/2,\bar{a}\kappa_*/4\}n}\bar{C}_2\left(1+\E\lsrs{|\bar{\theta}_0^\lambda|^{r_*}}\right) + \bar{C}_3\right) \nonumber\\
&\ + \hat{c}e^{-\min\{\dot{c}/4, \bar{a}/2, \bar{a}r/2,\bar{a}\kappa_*/4\}n}\Bigg{[}\left(1+\E\lsrs{|\bar{\theta}_0^\lambda|^{r_*}}\right)\nonumber\\
&\ + \left(2+ \int_{\R^d} V_2(\theta)\ \dee\pi_\beta(\theta)\right)\left(1+\E\lsrs{|\bar{\theta}_0^\lambda|^{r_*}}\right)\Bigg{]}\nonumber\\
\leq&\ \left(\bar{C}_0^{1/2}+\bar{C}_2+\hat{c}\left(3+\int_{\R^d} V_2(\theta)\ \dee\pi_\beta(\theta)\right)\right)\times e^{-\min\{\dot{c}/4, \bar{a}/2, \bar{a}r/2,\bar{a}\kappa_*/4\}n}\nonumber\\
&\ \times \left(1+\E\lsrs{|\bar{\theta}_0^\lambda|^{r_*}}\right) + \left(\bar{C}_1^{1/2}+\bar{C}_3\right)\lambda\nonumber\\
=&\ C_1 e^{-C_0(n+1)}\left(1+\E\lsrs{|\bar{\theta}_0^\lambda|^{r_*}}\right) + C_2\lambda,\nonumber
\end{align}
with
\begin{align}
C_0 &:= \min\{\dot{c}/4, \bar{a}/2, \bar{a}r/2,\bar{a}\kappa_*/4\},\nonumber\\
C_1 &:= e^{C_0}\lsrs{\bar{C}_0^{1/2}+\bar{C}_2+\hat{c}\left(3+\int_{\R^d} V_2(\theta)\ \dee\pi_\beta(\theta)\right)},\label{eqn:theorem2.7_constants_C0_C1_C2}\\
C_2 &:= \bar{C}_1^{1/2}+\bar{C}_3,\nonumber
\end{align}

where the constants $\kappa_*$, $\bar{C}_0$ and $\bar{C}_1$, $\dot{c}$ and $\hat{c}$, and $\bar{C}_2$ and $\bar{C}_3$ are given in Corollary \ref{corollary:4.3}, Lemma \ref{lemma:4.9}, Proposition \ref{proposition:4:11}, and Lemma \ref{lemma:4.12} respectively.  By replacing $t$ with $nT\in ((n-1)T, nT]$, this implies that for any $\lambda\in(0,\tilde{\lambda}_{\max})$ and $n\in\N_0$,
\begin{align}
W_1(\mc{L}(\bar{\theta}_{nT}^\lambda),\pi_\beta)\leq C_1 e^{-C_0n}\left(1+\E\lsrs{|\bar{\theta}_0^\lambda|^{r_*}}\right) + C_2\lambda.\nonumber
\end{align}
To obtain a non-asymptotic bound for the mTULA algorithm $(\theta_n^\lambda)_{n\in\N_0}$, we replace $nT$ with $n$ in the above inequality to obtain for any $\lambda\in(0,\tilde{\lambda}_{\max})$, $n\in\N_0$ that
\begin{align}
W_1(\mc{L}(\theta_{n}^\lambda), \pi_\beta)
\leq &\ C_1 e^{-C_0n/T}\left(1+\E\lsrs{|\theta_0|^{r_*}}\right) + C_2\lambda\nonumber\\
\leq &\ C_1 e^{-C_0\lambda n}\left(1+\E\lsrs{|\theta_0|^{r_*}}\right) + C_2\lambda,\nonumber
\end{align}
where the last inequality follows from $T\leq 1/\lambda$. This completes the proof. \qed\\

By using similar arguments as those described in the beginning of Section \ref{sec:mtproof}, we can also obtain an upper estimate in Wasserstein-2 distance. The detailed proof is provided below.\\
\subsection*{\textit{Proof of Theorem \ref{theorem:2.8}}}
\label{proof:theorem2.8}
Let $\lambda\in(0,\tilde{\lambda}_{\max})$, $n\in\N_0$, and let $t\in(nT, (n+1)T]$. Consider the splitting
\begin{align}
W_2(\mc{L}(\bar{\theta}_t^\lambda),\pi_\beta)&\leq W_2(\mc{L}(\bar{\theta}_t^\lambda),\mc{L}(\bar{\zeta}_t^{\lambda,n})) + W_2(\mc{L}(\bar{\zeta}_t^{\lambda,n}),\mc{L}(Z_t^\lambda)) + W_2(\mc{L}(Z_t^\lambda),\pi_\beta).\label{eqn:wasserstein_2_splitting}
\end{align}
Lemma \ref{lemma:4.9} gives an upper bound estimate for the first term of the RHS of (\ref{eqn:wasserstein_2_splitting}). By using $W_2\leq\sqrt{2w_{1,2}}$ in Lemma \ref{lemma:4.10} and modifying the proof of Lemma \ref{lemma:4.12}, one obtains
\begin{align}
W_2(\mc{L}(\bar{\zeta}_t^{\lambda,n}),\mc{L}(Z_t^\lambda))
\leq &\ \lambda^{1/2}\left(e^{-\min\{\dot{c}/8,\bar{a}/4,\bar{a}r/4,\bar{a}\kappa_*/8\}n}\bar{C}_4\left(1+\E\lsrs{|\bar{\theta}_0^\lambda|^{r_*}}\right)^{1/2}+\bar{C}_5\right),\label{eqn:theorem2.8_wass2_aux_langevin_ub}
\end{align}
where
\begin{align}
\bar{C}_4 &:= \sqrt{2\hat{c}}e^{\min\{\dot{c}/8,\bar{a}/4,\bar{a}r/4,\bar{a}\kappa_*/8\}}\left(1+\frac{1}{\min\{\dot{c}/8,\bar{a}/4,\bar{a}r/4,\bar{a}\kappa_*/8\}}\right)\left(\bar{C}_0^{1/2}+\frac{1}{\sqrt{2}}\right),\nonumber\\
\bar{C}_5 &:= \frac{4\sqrt{2\hat{c}}e^{\dot{c}/4}}{\dot{c}}\left(\bar{C}_1^{1/2}+\frac{1+2\sqrt{2}}{4}+\sqrt{\frac{c_*(2)}{2}}\left(1+\frac{1}{\bar{a}\kappa_*}\right)^{1/2}+\frac{(v_4(M_V(4)))^{1/2}}{4}\right),\label{eqn:theorem2.8_constants_C4_C5}
\end{align}
as an upper bound estimate for the second term of the RHS of (\ref{eqn:wasserstein_2_splitting}). The details of the proof of (\ref{eqn:theorem2.8_wass2_aux_langevin_ub}) are omitted as the arguments here follow the same lines as that of Lemma \ref{lemma:4.12}. Similarly, modifying (\ref{eqn:theorem2.7_wass1_ub}) yields
\begin{align}
W_2(\mc{L}(Z_t^\lambda),\pi_\beta)\leq \sqrt{2\hat{c}}e^{-\dot{c}\lambda t/2}\lsrs{1+\E[V_2(\theta_0)] + \int_{\R^d} V_2(\theta)\ \dee\pi_\beta(\theta)}^{1/2}.\label{eqn:theorem2.8_wass2_langevin_target}
\end{align}
Substituting the result of Lemma \ref{lemma:4.9}, (\ref{eqn:theorem2.8_wass2_aux_langevin_ub}), and (\ref{eqn:theorem2.8_wass2_langevin_target}) into (\ref{eqn:wasserstein_2_splitting}) yields
\begin{align}
W_2(\mc{L}(\bar{\theta}_t^\lambda),\pi_\beta)
\leq&\ \lambda\left(e^{-\bar{a}\min\{r,\kappa_*/2\}n}\bar{C}_0\E\lsrs{|\bar{\theta}_0^\lambda|^{r_*}} +\bar{C}_1\right)^{1/2}\nonumber\\
&\ + \lambda^{1/2}\left(e^{-\min\{\dot{c}/8,\bar{a}/4,\bar{a}r/4,\bar{a}\kappa_*/8\}n}\bar{C}_4\left(1+\E\lsrs{|\bar{\theta}_0^\lambda|^{r_*}}\right)^{1/2}+\bar{C}_5\right) \nonumber\\
&\ + \sqrt{2\hat{c}}e^{-\dot{c}\lambda t/2}\lsrs{1+\E[V_2(\theta_0)] + \int_{\R^d} V_2(\theta)\ \dee\pi_\beta(\theta)}^{1/2}\nonumber\\
\leq&\ \lambda\left(e^{-\min\{\dot{c}/8, \bar{a}/4, \bar{a}r/4,\bar{a}\kappa_*/8\}n}\bar{C}_0^{1/2}\left(1+\E\lsrs{|\bar{\theta}_0^\lambda|^{r_*}}\right)^{1/2} +\bar{C}_1^{1/2}\right)\nonumber\\
&\ + \lambda^{1/2}\left(e^{-\min\{\dot{c}/8,\bar{a}/4,\bar{a}r/4,\bar{a}\kappa_*/8\}n}\bar{C}_4\left(1+\E\lsrs{|\bar{\theta}_0^\lambda|^{r_*}}\right)^{1/2}+\bar{C}_5\right) \nonumber\\
&\ + \sqrt{2\hat{c}}e^{-\min\{\dot{c}/8, \bar{a}/4, \bar{a}r/4,\bar{a}\kappa_*/8\}n}\Bigg{[}\left(1+\E\lsrs{|\bar{\theta}_0^\lambda|^{r_*}}\right)^{1/2}\nonumber\\
&\ + \left(2+ \int_{\R^d} V_2(\theta)\ \dee\pi_\beta(\theta)\right)^{1/2}\left(1+\E\lsrs{|\bar{\theta}_0^\lambda|^{r_*}}\right)^{1/2}\Bigg{]}\nonumber\\
\leq&\ \left(\bar{C}_0^{1/2}+\bar{C}_4+\sqrt{2\hat{c}}\left(1+\left(2+\int_{\R^d} V_2(\theta)\ \dee\pi_\beta(\theta)\right)^{1/2}\right)\right)\times e^{-\min\{\dot{c}/8, \bar{a}/4, \bar{a}r/4,\bar{a}\kappa_*/8\}n}\nonumber\\
&\ \times \left(1+\E\lsrs{|\bar{\theta}_0^\lambda|^{r_*}}\right)^{1/2} + \left(\bar{C}_1^{1/2}+\bar{C}_5\right)\lambda^{1/2}\nonumber\\
=&\ C_4 e^{-C_3(n+1)}\left(1+\E\lsrs{|\bar{\theta}_0^\lambda|^{r_*}}\right)^{1/2} + C_5\lambda^{1/2},\label{eqn:theorem2.8_wass2_ub}
\end{align}
with
\begin{align}
C_3 &:= \min\{\dot{c}/8, \bar{a}/4, \bar{a}r/4,\bar{a}\kappa_*/8\},\nonumber\\
C_4 &:= e^{C_3}\lsrs{\bar{C}_0^{1/2}+\bar{C}_4+\sqrt{2\hat{c}}\left(1+\left(2+\int_{\R^d} V_2(\theta)\ \dee\pi_\beta(\theta)\right)^{1/2}\right)},\label{eqn:theorem2.8_constant_C3_C4_C5}\\
C_5 &:= \bar{C}_1^{1/2}+\bar{C}_5,\nonumber
\end{align}

where the constants $\kappa_*$, $\bar{C}_0$ and $\bar{C}_1$, $\dot{c}$ and $\hat{c}$, and $\bar{C}_4$ and $\bar{C}_5$ are given in Corollary \ref{corollary:4.3}, Lemma \ref{lemma:4.9}, Proposition \ref{proposition:4:11}, and (\ref{eqn:theorem2.8_constants_C4_C5}) respectively. Finally, as with the proof of Theorem \ref{theorem:2.7}, replacing $t$ with $nT$ in (\ref{eqn:theorem2.8_wass2_ub}), and then replacing $nT$ with $n$, yields, for every $\lambda\in(0,\tilde{\lambda}_{\max})$ and $n\in\N_0$,
\begin{align}
W_2(\mc{L}(\theta_n^\lambda),\pi_\beta)\leq C_4e^{-C_3\lambda n}\left(1+\E\lsrs{|\theta_0|^{r_*}}\right)^{1/2} + C_5\lambda^{1/2}.
\end{align}
This completes the proof. \qed

\section{Proof of Auxiliary Results}
\label{sec:proof_aux}
\subsection{Proof of Remarks in Section 2}
\subsection*{\textit{Proof of Remark \ref{remark:2.4}}}
\label{proof:remark2.4}
We first consider the case when $r>0$. By Assumptions \ref{assumption:3} and \ref{assumption:2}, we have, for any $\theta\in\R^d$,
\begin{align}
\lara{\theta, h(\theta)}
&\geq a|\theta|^{r+2} - b|\theta|^{\bar{r}+2}+\lara{\theta,h(0)}\nonumber\\
&\geq a|\theta|^{r+2} - b|\theta|^{\bar{r}+2}-|\theta|\cdot|h(0)|\nonumber\\
&\geq a|\theta|^{r+2} - b|\theta|^{\bar{r}+2}-|\theta a^{1/2}|\cdot \frac{K}{a^{1/2}}\nonumber\\
&\geq a|\theta|^{r+2} - b|\theta|^{\bar{r}+2}-\frac{a}{2}|\theta|^2 - \frac{K^2}{2a}\label{eqn:remark2.4_ip_lb1}\\
& = \frac{a}{2}|\theta|^{r+2} + \left(\frac{a}{4}|\theta|^{r+2}-b|\theta|^{\bar{r}+2}\right) + \left(\frac{a}{4}|\theta|^{r+2}-\frac{a}{2}|\theta|^2\right) - \frac{K^2}{2a}.\label{eqn:remark2.4_ip_lb2}
\end{align}
Observe that
\begin{align}
\left(\frac{a}{4}|\theta|^{r+2}-b|\theta|^{\bar{r}+2}\right) > 0 &\iff |\theta| > \left(\frac{4b}{a}\right)^{1/(r-\bar{r})},\nonumber\\
\left(\frac{a}{4}|\theta|^{r+2}-\frac{a}{2}|\theta|^2\right) > 0 &\iff |\theta| > 2^{1/r}.\label{eqn:remark2.4_ip_comp_lb}
\end{align}
Hence, setting $R:=\max\left\{\left(\frac{4b}{a}\right)^{1/(r-\bar{r})},2^{1/r}\right\}>1$, we obtain from (\ref{eqn:remark2.4_ip_lb2}) and (\ref{eqn:remark2.4_ip_comp_lb}),
\begin{align}
\lara{\theta, h(\theta)} \geq \frac{a}{2}|\theta|^{r+2}-\frac{K^2}{2a},\qquad \forall|\theta|>R,\label{eqn:remark2.4_ip_comp1_lb}
\end{align}
and from (\ref{eqn:remark2.4_ip_lb1}),
\begin{align}
\lara{\theta, h(\theta)} 
&\geq \frac{a}{2}|\theta|^{r+2}-b |\theta|^{\bar{r}+2} - \frac{a}{2}|\theta|^2-\frac{K^2}{2a}\nonumber\\
&\geq \frac{a}{2}|\theta|^{r+2} -\left(b+\frac{a}{2}\right)R^{\bar{r}+2}-\frac{K^2}{2a},\qquad \forall|\theta|\leq R.\label{eqn:remark2.4_ip_comp2_lb}
\end{align}
Combining (\ref{eqn:remark2.4_ip_comp1_lb}) and (\ref{eqn:remark2.4_ip_comp2_lb}) yields the inequality
\begin{align}
\lara{\theta, h(\theta)}\geq \frac{a}{2}|\theta|^{r+2}-\left(\left(b+\frac{a}{2}\right)R^{\bar{r}+2}+\frac{K^2}{2a}\right).\label{eqn:remark2.4_ip_lb3}
\end{align}
In the case $r=0$, Assumption \ref{assumption:3} reads
\begin{align}
\lara{\theta,h(\theta)} \geq \tilde{a}\lbrb{\theta}^{2}-\tilde{b}=\tilde{a}\lbrb{\theta}^{r+2}-\tilde{b}\label{eqn:remark2.4_ip_lb4}
\end{align}
for all $\theta\in\R^d$. Hence, combining (\ref{eqn:remark2.4_ip_lb3}) and (\ref{eqn:remark2.4_ip_lb4}) yields, for $r\geq 0$,
\begin{align}
\lara{\theta,h(\theta)} \geq \bar{a}\lbrb{\theta}^{r+2}-\bar{b}\nonumber
\end{align}
for all $\theta\in\R^d$, where
\begin{align}
\bar{a} &:= \tfrac{a}{2}\1_{\{r>0\}}+\tilde{a}\1_{\{r=0\}},\nonumber\\
\bar{b} &:= \left(\left(b+\tfrac{a}{2}\right)R^{\bar{r}+2}+\tfrac{K^2}{2a}\right)\1_{\{r>0\}}+\tilde{b}\1_{\{r=0\}},\label{eqn:remark2.4_constant_a_b_R}\\
R &:=\max\left\{\left(\tfrac{4b}{a}\right)^{1/(r-\bar{r})},2^{1/r}\right\},\nonumber
\end{align}
completing the proof of the first inequality of Remark \ref{remark:2.4}.\\

To establish the second inequality of Remark \ref{remark:2.4}, set 
\begin{align}
\bar{b}':=(\bar{b}+2^{2/r}\bar{a})\1_{\{r>0\}} + \tilde{b}\1_{\{r=0\}}.\label{eqn:remark2.4_constant_bprime}
\end{align}

In the case $r>0$, it suffices to prove that for all $\theta\in\R^d$,
\begin{align}
\bar{a}|\theta|^{r+2}-\bar{b} \geq \bar{a}|\theta|^2 -\left(\bar{b}+2^{2/r}\bar{a}\right)
&\iff \bar{a}|\theta|^{r+2} -\bar{a}|\theta|^2+2^{2/r}\bar{a} \geq 0\nonumber\\
&\iff \frac{\bar{a}}{2}|\theta|^{r+2}+\left(\frac{\bar{a}}{2}|\theta|^{r+2}-\bar{a}|\theta|^2\right) + 2^{2/r}\bar{a}\geq 0.\nonumber
\end{align}
Observe that 
\begin{align}
\left(\frac{\bar{a}}{2}|\theta|^{r+2}-\bar{a}|\theta|^2\right)>0\iff |\theta|>2^{1/r}.\nonumber
\end{align}
Hence, for all $|\theta|>2^{1/r}$,
\begin{align}
\bar{a}|\theta|^{r+2} -\bar{a}|\theta|^2+2^{2/r}\bar{a}
&\geq \frac{\bar{a}}{2}|\theta|^{r+2}+2^{2/r}\bar{a}\geq 0,\label{eqn:remark2.4_poly_comp1_lb}
\end{align}
and for all $|\theta|\leq 2^{1/r}$,
\begin{align}
\bar{a}|\theta|^{r+2} -\bar{a}|\theta|^2+2^{2/r}\bar{a}
&\geq -\bar{a}\cdot 2^{2/r}+2^{2/r}\bar{a}= 0.\label{eqn:remark2.4_poly_comp2_lb}
\end{align}
The veracity of the second inequality of Remark \ref{remark:2.4} is established from (\ref{eqn:remark2.4_poly_comp1_lb}) and (\ref{eqn:remark2.4_poly_comp2_lb}) in the case $r>0$. In the case $r=0$, the second inequality of Remark \ref{remark:2.4} holds trivially. This completes the proof. \qed

\subsection*{\textit{Proof of Remark \ref{remark:2.5}}} 
\label{proof:remark2.5}
We first consider the case where $r>0$. The remark holds trivially when $\theta=\theta'$, hence for the remainder of the proof we address the case $\theta\neq \theta'$. Assumption \ref{assumption:3} states that for all $\theta,\theta'\in\R^d$,
\begin{align}
\lara{\theta-\theta',h(\theta)-h(\theta')}\geq a\lbrb{\theta-\theta'}^2(\lbrb{\theta}^{r}+\lbrb{\theta'}^{r})-b\lbrb{\theta-\theta'}^2(\lbrb{\theta}^{\bar{r}}+\lbrb{\theta'}^{\bar{r}}).\nonumber
\end{align}
Observe that for all $\theta,\theta'\in\R^d$,
\begin{align}
a\lbrb{\theta-\theta'}^2|\theta|^r - b\lbrb{\theta-\theta'}^2|\theta|^{\bar{r}}>0\iff |\theta| > \left(\frac{b}{a}\right)^{1/(r-\bar{r})}.\label{eqn:remark2.5_ip_lb}
\end{align}
Moreover, when $\theta\neq\theta'$, the statement in (\ref{eqn:remark2.5_ip_lb}) still holds if both strict inequalities are replaced with inequalities in (\ref{eqn:remark2.5_ip_lb}). Denote $\bar{R}:=\left(\frac{b}{a}\right)^{1/(r-\bar{r})}$ and $\bar{B}(0,\bar{R})$ (resp. $B(0, \bar{R})$) as the closed (resp. open) ball with radius $\bar{R}$ centred at the zero vector in $\R^d$. It follows from (\ref{eqn:remark2.5_ip_lb}) that for all $\theta,\theta'\notin B(0,\bar{R})$ such that $\theta\neq \theta'$ and either one of $\theta, \theta'$ does not lie in $\bar{B}(0,\bar{R})$,
\begin{align}
\lara{\theta-\theta',h(\theta)-h(\theta')} > 0.\label{eqn:remark2.5_ip_lb2}
\end{align}
Moreover, from the Cauchy-Schwarz inequality and Assumption \ref{assumption:2} that, for all $\theta,\theta'\in \bar{B}(0,\bar{R})$,
\begin{align}
-\lara{\theta-\theta',h(\theta)-h(\theta')}
\leq&\ \lbrb{\theta-\theta'}\cdot \lbrb{h(\theta)-h(\theta')}\nonumber\\
\leq&\ L(1+\lbrb{\theta}+\lbrb{\theta'})^{r}\lbrb{\theta-\theta'}^2\nonumber\\
\leq&\ L(1+2\bar{R})^{r}\lbrb{\theta-\theta'}^2.\label{eqn:remark2.5_ip_lb3}
\end{align}

For the case $|\theta|<\bar{R}$, $|\theta'| > \bar{R}$, consider $\bar{\theta}:= \theta' + c_{\theta,\theta'}(\theta-\theta')$, where $c_{\theta,\theta'}$ is chosen as the unique scalar in the interval $(0,1)$ such that $|\bar{\theta}| = \bar{R}$. Observe that $\theta-\bar{\theta}=(1-c_{\theta,\theta'})(\theta-\theta')$. One then obtains, for all $|\theta|< \bar{R}$, $|\theta'| > \bar{R}$,
\begin{align}
\lara{\theta-\theta',h(\theta)-h(\theta')}
=&\ \lara{\theta-\theta',h(\theta)-h(\bar{\theta})} + \lara{\theta-\theta',h(\bar{\theta})-h(\theta')}\nonumber\\
=&\ \lara{\theta-\theta',h(\theta)-h(\bar{\theta})} + \frac{1}{c_{\theta,\theta'}}\lara{\bar{\theta}-\theta',h(\bar{\theta})-h(\theta')}\nonumber\\
>&\ \lara{\theta-\theta',h(\theta)-h(\bar{\theta})}\label{eqn:remark2.5_ip_lb4}\\
=&\ \frac{1}{1-c_{\theta,\theta'}}\lara{\theta-\bar{\theta}, h(\theta)-h(\bar{\theta})}\nonumber\\
\geq & - \frac{L(1+2\bar{R})^{r}}{1-c_{\theta,\theta'}}\lbrb{\theta-\bar{\theta}}^2\label{eqn:remark2.5_ip_lb5}\\
=&\ -L(1+2\bar{R})^{r}(1-c_{\theta,\theta'})\lbrb{\theta-\theta'}^2\nonumber\\
\geq&\ -L(1+2\bar{R})^{r}\lbrb{\theta-\theta'}^2,\label{eqn:remark2.5_ip_lb6}
\end{align}

where the inequalities (\ref{eqn:remark2.5_ip_lb4}) and (\ref{eqn:remark2.5_ip_lb5}) follow respectively from (\ref{eqn:remark2.5_ip_lb2}) and (\ref{eqn:remark2.5_ip_lb3}), and the last inequality follows from $c_{\theta,\theta'}\in (0,1)$. Interchanging the roles of $\theta$ and $\theta'$ in the above argument shows that (\ref{eqn:remark2.5_ip_lb6}) holds for $|\theta|>\bar{R}$, $|\theta'|< \bar{R}$. Hence, when $r>0$, one concludes that for all $\theta,\theta'\in\R^d$,
\begin{align}
\lara{\theta-\theta',h(\theta)-h(\theta')}\geq -L(1+2\bar{R})^{r}\lbrb{\theta-\theta'}^2.
\end{align}

For the case $r=0$, one observes that the same argument used to derive inequality (\ref{eqn:remark2.5_ip_lb3}) applies ad verbatim and is valid for all $\theta,\theta'\in\R^d$. That is, for all $\theta,\theta'\in\R^d$,
\begin{align}
-\lara{\theta-\theta',h(\theta)-h(\theta')}
\leq&\ \lbrb{\theta-\theta'}\cdot \lbrb{h(\theta)-h(\theta')}\nonumber\\
\leq&\ L(1+\lbrb{\theta}+\lbrb{\theta'})^{r}\lbrb{\theta-\theta'}^2\nonumber\\
=&\ L\lbrb{\theta-\theta'}^2.\nonumber
\end{align}
Therefore, by denoting $\bar{L}:=L(1+2\bar{R})^r$ and $\bar{R}:=(b/a)^{1/(r-\bar{r})}\1_{\{r>0\}}$, one obtains
\begin{align}
\lara{\theta-\theta', h(\theta)-h(\theta')} \geq -\bar{L}|\theta-\theta'|^2,\qquad \forall\theta,\theta'\in\R^d.\nonumber
\end{align} 
This completes the proof. \qed

\subsection*{\textit{Proof of Remark \ref{remark:2.6}(i)}} 
\label{proof:remark2.6i}
From Assumption \ref{assumption:4}, one obtains for all $\theta\in\R^d$,
\begin{align}
|\nabla h(\theta)|
&\leq |\nabla h(\theta) - \nabla h(0)| + |\nabla h(0)|\nonumber\\
&\leq L_\nabla (1+|\theta|)^\nu \cdot |\theta| + |\nabla h(0)|\nonumber\\
&\leq 2^{\nu-1}L_\nabla (1+|\theta|^\nu)\cdot |\theta| + |\nabla h(0)|\nonumber\\
&\leq 2^{\nu-1}L_\nabla (1+|\theta|^{\nu+1}) + 2^{\nu-1}L_\nabla |\theta|^{\nu + 1} + |\nabla h(0)|\nonumber\\
&\leq \frac{C_\nabla}{2} (1+|\theta|^{\nu+1}) + \frac{C_\nabla}{2} |\theta|^{\nu + 1} + \frac{C_\nabla}{2}\nonumber\\
&= C_\nabla (1+|\theta|^{\nu+1}),\nonumber
\end{align}
where $C_\nabla = 2\max\{2^{\nu-1}L_\nabla, |\nabla h(0)|\}$. This completes the proof. \qed

\subsection*{\textit{Proof of Remark \ref{remark:2.6}(ii)}} 
\label{proof:remark2.6ii}
For $\theta,\theta'\in \R^d$, denote the auxiliary function $g_{\theta,\theta'}(t):=h(t\theta + (1-t)\theta')$ such that $g_{\theta,\theta'}(1)=h(\theta)$ and $g_{\theta,\theta'}(0)=h(\theta')$. Observe that $g_{\theta,\theta'}'(t)=\nabla h(t\theta + (1-t)\theta')(\theta - \theta')$. For $\theta,\theta'\in \R^d$, one obtains, from the Fundamental Theorem of Calculus, Remark \ref{remark:2.6}(i) and Jensen's inequality that,
\begin{align}
|h(\theta)-h(\theta')|
&= |g_{\theta,\theta'}(1) - g_{\theta,\theta'}(0)|\nonumber\\
&= \lbrb{\int^1_0 g_{\theta,\theta'}'(t)\ \dee t}\nonumber\\
&\leq \int^1_0 \lbrb{g_{\theta,\theta'}'(t)}\ \dee t\nonumber\\
&\leq \int^1_0 \lbrb{\nabla h(t\theta + (1-t)\theta')}\cdot \lbrb{\theta - \theta'}\ \dee t\nonumber\\
&\leq C_\nabla\int^1_0 \left(1+|t\theta + (1-t)\theta'|^{\nu+1}\right)\cdot \lbrb{\theta - \theta'}\ \dee t\nonumber\\
&\leq C_\nabla\int^1_0 \left(1+|\theta|^{\nu+1} + |\theta'|^{\nu+1} \right)\cdot \lbrb{\theta - \theta'}\ \dee t\nonumber\\
&\leq C_\nabla \left(1+|\theta| + |\theta'|\right)^{\nu+1}|\theta - \theta'|.\nonumber
\end{align}
This completes the proof. \qed

\subsection*{\textit{Proof of Remark \ref{remark:2.6}(iii)}} 
\label{proof:remark2.6iii}
For $\theta,\theta'\in \R^d$, define $g_{\theta,\theta'}(t)$ as in the proof of Remark \ref{remark:2.6}(ii). One obtains by a similar argument
\begin{align}
|h(\theta)-h(\theta')-\nabla h(\theta')(\theta-\theta')|
&= |g_{\theta,\theta'}(1) - g_{\theta,\theta'}(0) - g_{\theta,\theta'}'(0)|\nonumber\\
&= \lbrb{\int^1_0 \left(g_{\theta,\theta'}'(t) - g_{\theta,\theta'}'(0)\right)\ \dee t}\nonumber\\
&\leq \int^1_0 \lbrb{g_{\theta,\theta'}'(t) - g_{\theta,\theta'}'(0)}\ \dee t\nonumber\\
&\leq \int^1_0 \lbrb{\nabla h(t\theta + (1-t)\theta') - \nabla h(\theta')}\cdot \lbrb{\theta - \theta'}\ \dee t\nonumber\\
&\leq L_\nabla\int^1_0 \left(1+|t\theta + (1-t)\theta'| + |\theta'|\right)^\nu\cdot t\lbrb{\theta - \theta'}^2\ \dee t\nonumber\\
&\leq 3^{\nu - 1}L_\nabla\int^1_0 \left(1+|t\theta + (1-t)\theta'|^\nu + |\theta'|^\nu\right)\cdot t\lbrb{\theta - \theta'}^2\ \dee t\nonumber\\
&\leq 3^{\nu - 1}L_\nabla\int^1_0 \left(1+|\theta|^\nu + 2|\theta'|^\nu \right)\cdot t\lbrb{\theta - \theta'}^2\ \dee t\nonumber\\
&\leq 2\cdot 3^{\nu - 1}L_\nabla \left(1+|\theta|^\nu + |\theta'|^\nu \right)\cdot \lbrb{\theta - \theta'}^2\cdot \int^1_0 t\ \dee t\nonumber\\
& = \bar{L}_\nabla \left(1+|\theta|^\nu + |\theta'|^\nu \right)|\theta - \theta'|^2,\nonumber
\end{align}
where $\bar{L}_\nabla = 3^{\nu - 1}L_\nabla$.  This completes the proof. \qed

\subsection{Proof of Auxiliary Results in Section 3}
\subsection*{\textit{Proof of Proposition \ref{proposition:3.1}}}
\label{proof:proposition3.1}
We verify the assumptions for each choice of $U$ considered.
\begin{enumerate}
\item[(i)] For the multivariate standard Gaussian distribution, $h(\theta)=\theta$ and $\nabla h(\theta)=I_d$. Hence, Assumptions \ref{assumption:2}, \ref{assumption:3}, \ref{assumption:4} hold with $r=0$, $\nu=0$, $L=1$, $K=1$, $\tilde{a}=1$, $\tilde{b}=1$, $L_\nabla=1$.\\
\item[(ii)] For the given multivariate Gaussian mixture model, we have, for all $\theta\in\R^d$,
\begin{align}
h(\theta)=\theta-\dot{a}+\frac{2\dot{a}}{1+e^{2\lara{\dot{a},\theta}}},\quad \nabla h(\theta)=I_d -\frac{4e^{2\lara{\dot{a},\theta}}}{(1+e^{2\lara{\dot{a},\theta}})^2}\dot{a}\dot{a}^T.\nonumber
\end{align}
We first verify Assumption \ref{assumption:2}. Fix $\theta,\theta'\in\R^d$ and denote here $x:=2\lara{\dot{a},\theta}$, $y:=2\lara{\dot{a},\theta'}$. Then by the inequality $e^w\geq 1+w$, for all $w\in\R$, and the Cauchy-Schwarz inequality,
\begin{align}
|h(\theta)-h(\theta')|
=&\ \lbrb{(\theta-\theta') + 2\dot{a}\left(\frac{1}{1+e^x}-\frac{1}{1+e^y}\right)}\nonumber\\
\leq&\ |\theta-\theta'| + 2|\dot{a}|\cdot\frac{|e^x-e^y|}{(1+e^x)(1+e^y)}\nonumber\\
=&\ |\theta-\theta'| + 2|\dot{a}|\cdot\frac{e^{\max{\{x,y\}}}\cdot \left(1-e^{-|x-y|}\right)}{(1+e^x)(1+e^y)}\nonumber\\
\leq&\ |\theta-\theta'| + 2|a|\left(1-e^{-|x-y|}\right)\nonumber\\
\leq&\ |\theta-\theta'| + 2|\dot{a}|\cdot |x-y|\nonumber\\
=&\ |\theta-\theta'| + 4|\dot{a}|\cdot \lbrb{\lara{\dot{a},\theta-\theta'}}\nonumber\\
\leq&\ \left(1+4|\dot{a}|^2\right)\cdot|\theta-\theta'|.\nonumber
\end{align}
Furthermore, we have
\begin{align}
|h(\theta)|
\leq&\ |\theta| + \lbrb{1-\frac{2}{1+e^{2\lara{\dot{a},\theta}}}}\cdot |\dot{a}|\nonumber\\
=&\ |\theta| + \lbrb{\frac{e^{2\lara{\dot{a},\theta}}-1}{e^{2\lara{\dot{a},\theta}}+1}}\cdot |\dot{a}|\nonumber\\
\leq&\ |\theta|+|\dot{a}|.\nonumber
\end{align}
Hence, Assumption \ref{assumption:2} holds with $r=0$, $L=1+4|\dot{a}|^2$, and $K=\max\{1,|\dot{a}|\}$.\\

\noindent Next, we verify Assumption \ref{assumption:3}. We have, by the Cauchy-Schwarz inequality,
\begin{align}
\lara{\theta, h(\theta)}
=&\ |\theta|^2-\frac{e^{2\lara{\dot{a},\theta}}-1}{e^{2\lara{\dot{a},\theta}}+1}\cdot \lara{\dot{a},\theta}\nonumber\\
\geq&\ |\theta|^2-\lbrb{\lara{\dot{a},\theta}}\nonumber\\
\geq&\ \tfrac{1}{2}|\theta|^2-\left(|\dot{a}|\cdot |\theta|-\tfrac{1}{2}|\theta|^2\right)\1_{\{|\theta| > 2|\dot{a}|\}}-\left(|\dot{a}|\cdot |\theta|-\tfrac{1}{2}|\theta|^2\right)\1_{\{|\theta|\leq 2|\dot{a}|\}}\nonumber\\
\geq&\ \tfrac{1}{2}|\theta|^2-\left(2|\dot{a}|^2-\tfrac{1}{2}|\theta|^2\right)\1_{\{|\theta|\leq 2|\dot{a}|\}}\nonumber\\
\geq&\ \tfrac{1}{2}|\theta|^2-2|\dot{a}|^2.\nonumber
\end{align}
Hence, Assumption \ref{assumption:3} holds with $\tilde{a}=\frac{1}{2}$ and $\tilde{b}=2$.\\

\noindent Finally, we verify Assumption \ref{assumption:4}. Again, fix $\theta,\theta'\in\R^d$ and denote here $x:=2\lara{\dot{a},\theta}$, $y:=2\lara{\dot{a},\theta'}$. Then,
\begin{align}
\lbrb{\nabla h(\theta) - \nabla h(\theta')}
=&\ 4\lbrb{\frac{e^x}{(1+e^x)^2} - \frac{e^y}{(1+e^y)^2}}\cdot |\dot{a}\dot{a}^T|\nonumber\\
=&\ 4|\dot{a}|^2\cdot \frac{\lbrb{e^x-e^y+e^{x+2y}-e^{2x+y}}}{(1+e^x)^2(1+e^y)^2}\nonumber\\
\leq&\ 4|\dot{a}|^2\cdot \frac{\lbrb{e^x-e^y}+e^{x+y}\lbrb{e^x-e^y}}{(1+e^x)^2(1+e^y)^2}\nonumber\\
\leq&\ 4|\dot{a}|^2\cdot \frac{\lbrb{e^x-e^y}}{(1+e^x)(1+e^y)}\nonumber\\
\leq&\ 8|\dot{a}|^3\cdot |\theta-\theta'|,\nonumber
\end{align}
where, for the last inequality, we recognised the same expression which appeared in the verification of Assumption \ref{assumption:2}. Hence, Assumption \ref{assumption:4} holds with $\nu=0$ and $L_\nabla = 8|\dot{a}|^3$.\\

\item[(iii)] For the double-well potential, we have, for every $\theta\in\R^d$,
\begin{align}
h(\theta) = (|\theta|^2-1)\theta,\quad \nabla h(\theta) = (|\theta|^2-1)I_d+2\theta\theta^T.
\end{align}
For any $\theta,\theta'\in\R^d$, one obtains
\begin{align}
|h(\theta)-h(\theta')|
=&\ \lbrb{\left(|\theta|^2\theta-|\theta'|^2\theta'\right)-(\theta-\theta')}\nonumber\\
\leq&\ \lbrb{|\theta|^2\theta-|\theta|^2\theta'} + \lbrb{|\theta|^2\theta' - |\theta'|^2\theta'} + |\theta-\theta'|\nonumber\\
\leq&\ |\theta|^2\cdot|\theta-\theta'| + |\theta'|\cdot\left(|\theta| + |\theta'|\right)\cdot|\theta-\theta'| + |\theta-\theta'|\nonumber\\
\leq&\ (|\theta| + |\theta'|)^2\cdot |\theta-\theta'| + |\theta-\theta'|\nonumber\\
\leq& \left(1+|\theta|+|\theta'|\right)^2|\theta-\theta'|\nonumber
\end{align}
and
\begin{align}
|h(\theta)|
\leq&\ |\theta|^3+|\theta|\leq 1+2|\theta|^3\nonumber\\
\leq&\ 2\left(1+|\theta|^3\right).\nonumber
\end{align}
Hence, Assumption \ref{assumption:2} holds with $r=2$, $L=1$, and $K=2$. Furthermore,
\begin{align}
\lara{\theta-\theta', |\theta|^2\theta-|\theta'|^2\theta'}
=&\ \lara{\theta-\theta', |\theta|^2\theta-|\theta|^2\theta'+|\theta|^2\theta'-|\theta'|^2\theta'}\nonumber\\
=&\ |\theta|^2\cdot |\theta-\theta'|^2+\left(\lara{\theta,\theta'}-|\theta'|^2\right)\left(|\theta|^2-|\theta'|^2\right).\nonumber
\end{align}
Interchanging $\theta$ and $\theta'$ in the above yields
\begin{align}
\lara{\theta-\theta', |\theta|^2\theta-|\theta'|^2\theta'}=&\ |\theta'|^2\cdot |\theta-\theta'|^2-\left(\lara{\theta,\theta'}-|\theta|^2\right)\left(|\theta|^2-|\theta'|^2\right).\nonumber
\end{align}
Taking the average of these two expressions, 
\begin{align}
\lara{\theta-\theta', |\theta|^2\theta-|\theta'|^2\theta'}
=&\ \tfrac{1}{2}\left(|\theta|^2+|\theta'|^2\right)|\theta-\theta'|^2+\tfrac{1}{2}\left(|\theta|^2-|\theta'|^2\right)^2\nonumber\\
\geq&\ \tfrac{1}{2}\left(|\theta|^2+|\theta'|^2\right)|\theta-\theta'|^2.\nonumber
\end{align}
This implies that
\begin{align}
\lara{\theta-\theta', h(\theta)-h(\theta')}\geq \tfrac{1}{2}\left(|\theta|^2+|\theta'|^2\right)|\theta-\theta'|^2-|\theta-\theta'|^2.\nonumber
\end{align}
Consequently, Assumption \ref{assumption:3} holds with $a=1/2$, $b=1$, and $\bar{r}=0$. Finally,
\begin{align}
&\lbrb{\nabla h(\theta)-\nabla h(\theta')}\nonumber\\
=&\ \lbrb{\left(|\theta|^2-|\theta'|^2\right)I_d+2\left(\theta\theta^T-\theta'\theta'^T\right)}\nonumber\\
\leq&\ \left(|\theta|+|\theta'|\right)|\theta-\theta'| + 2\sup_{w\in\R^d: |w|\leq 1}\lbrb{\lara{w,\theta}\theta - \lara{w,\theta'}\theta'}\nonumber\\
=&\ \left(|\theta|+|\theta'|\right)|\theta-\theta'| + 2\sup_{w\in\R^d: |w|\leq 1}\lbrb{\lara{w,\theta}\theta - \lara{w,\theta}\theta' + \lara{w,\theta}\theta' - \lara{w,\theta'}\theta'}\nonumber\\
\leq&\ \left(|\theta|+|\theta'|\right)|\theta-\theta'| + 2\sup_{w\in\R^d: |w|\leq 1}\lsrs{\lbrb{\lara{w,\theta}}\cdot |\theta-\theta'| + \lbrb{\lara{w,\theta-\theta'}}|\theta'|}\nonumber\\
\leq&\ \left(|\theta|+|\theta'|\right)|\theta-\theta'| + 2\lsrs{|\theta|\cdot |\theta-\theta'| + |\theta-\theta'|\cdot |\theta'|}\nonumber\\
=&\ 3(|\theta|+|\theta'|)\cdot |\theta-\theta'|.\nonumber
\end{align}
Therefore, Assumption \ref{assumption:4} holds with $\nu=1$ and $L_\nabla=3$.\\

\end{enumerate}

This completes the proof. \qed
\subsection{Proof of Auxiliary Results in Section 4}
\subsection*{\textit{Proof of Lemma \ref{lemma:4.2}}(i)}
\label{proof:lemma4.2}
For any $\lambda \in (0,\lambda_{1,\max})$, $t\in (n,n+1]$, $n\in\N_0$, define
\begin{align}
\Delta_{n,t}^\lambda := \bar{\theta}_n^\lambda-\lambda h_\lambda(\bar{\theta}_n^\lambda)(t-n),\quad \Xi_{n,t}^\lambda:=\sqrt{2\lambda \beta^{-1}}(B_t^\lambda-B_n^\lambda),\label{eqn:lemma4.2_delta_xi}
\end{align}
so that 
\begin{align}
\bar{\theta}_{t}^\lambda = \Delta_{n,t}^\lambda + \Xi_{n,t}^\lambda. \label{eqn:lemma4.2_tula_interpolated_splitting}
\end{align}
Taking the norm-squared of (\ref{eqn:lemma4.2_tula_interpolated_splitting}) and applying $\E[\cdot\mid \bar{\theta}_n^\lambda]$ yields
\begin{align}
\E\lsrs{\ce{|\bar{\theta}_t^\lambda|^2} \bar{\theta}_n^\lambda}
&=\E\lsrs{\ce{|\Delta_{n,t}^\lambda|^2} \bar{\theta}_n^\lambda}+ 2 \E\lsrs{\ce{\lara{\Delta_{n,t}^\lambda, \Xi_{n,t}^\lambda}} \bar{\theta}_n^\lambda}+ \E\lsrs{\ce{|\Xi_{n,t}^\lambda|^2} \bar{\theta}_n^\lambda}\nonumber\\
&=|\Delta_{n,t}^\lambda|^2+ 2\lara{\Delta_{n,t}^\lambda, \E\lsrs{\Xi_{n,t}^\lambda}}+ 2\lambda\beta^{-1} \E\lsrs{|B_t^\lambda-B_n^\lambda|^2}\nonumber\\
&=|\Delta_{n,t}^\lambda|^2+ 2\lambda\beta^{-1} \E\lsrs{|B_t^\lambda-B_n^\lambda|^2}\nonumber\\
&=|\Delta_{n,t}^\lambda|^2+ 2\lambda \beta^{-1} (t-n)d.\label{eqn:lemma4.2_tula_interpolated_moment2_condexp_eq}
\end{align}
Furthermore, by Remark \ref{remark:2.4} and Assumption \ref{assumption:2}, the first term of the RHS of (\ref{eqn:lemma4.2_tula_interpolated_moment2_condexp_eq}) can be bounded above by
\begin{align}
&|\Delta_{n,t}^\lambda|^2\nonumber\\
&=|\bar{\theta}_n^\lambda|^2 -2\lambda(t-n)\lara{\bar{\theta}_n^\lambda, h_\lambda(\bar{\theta}_n^\lambda)} + \lambda^2(t-n)^2\lbrb{h_\lambda(\bar{\theta}_n^\lambda)}^2\nonumber\\
&=|\bar{\theta}_n^\lambda|^2 -2\lambda(t-n)\frac{\lara{\bar{\theta}_n^\lambda, h(\bar{\theta}_n^\lambda)}}{(1+\lambda\lbrb{\bar{\theta}_n^\lambda}^{2r})^{1/2}} + \lambda^2(t-n)^2\frac{\lbrb{h(\bar{\theta}_n^\lambda)}^2}{1+\lambda\lbrb{\bar{\theta}_n^\lambda}^{2r}}\nonumber\\
&\leq |\bar{\theta}_n^\lambda|^2 -\frac{2\lambda(t-n)(\bar{a}|\bar{\theta}_n^\lambda|^{r+2}-\bar{b})}{(1+\lambda\lbrb{\bar{\theta}_n^\lambda}^{2r})^{1/2}} + \frac{\lambda^2(t-n)^2 K^2 (1+|\bar{\theta}_n^\lambda|^{r+1})^2}{1+\lambda\lbrb{\bar{\theta}_n^\lambda}^{2r}}\nonumber\\
&\leq |\bar{\theta}_n^\lambda|^2 -\frac{2\lambda(t-n)(\bar{a}|\bar{\theta}_n^\lambda|^{r+2}-\bar{b})}{(1+\lambda\lbrb{\bar{\theta}_n^\lambda}^{2r})^{1/2}} + \frac{2\lambda^2(t-n)^2 K^2 (1+|\bar{\theta}_n^\lambda|^{2r+2})}{1+\lambda\lbrb{\bar{\theta}_n^\lambda}^{2r}}\nonumber\\
&\leq |\bar{\theta}_n^\lambda|^2 -\frac{2\lambda(t-n)\bar{a}|\bar{\theta}_n^\lambda|^{r+2}}{(1+\lambda\lbrb{\bar{\theta}_n^\lambda}^{2r})^{1/2}} + 2\lambda(t-n)\bar{b} + \frac{2\lambda^2(t-n)^2 K^2|\bar{\theta}_n^\lambda|^{2r+2}}{1+\lambda\lbrb{\bar{\theta}_n^\lambda}^{2r}} + 2\lambda^2(t-n)^2 K^2 \nonumber\\
&= |\bar{\theta}_n^\lambda|^2 - \lambda(t-n)\frac{\bar{a}|\bar{\theta}_n^\lambda|^{r+2}}{(1+\lambda\lbrb{\bar{\theta}_n^\lambda}^{2r})^{1/2}}-\lambda(t-n)\lsrs{\frac{\bar{a}|\bar{\theta}_n^\lambda|^{r+2}}{(1+\lambda\lbrb{\bar{\theta}_n^\lambda}^{2r})^{1/2}}-\frac{2\lambda(t-n) K^2|\bar{\theta}_n^\lambda|^{2r+2}}{1+\lambda\lbrb{\bar{\theta}_n^\lambda}^{2r}}}\nonumber\\
&\quad + 2\lambda(t-n)(\bar{b}+\lambda(t-n)K^2)\nonumber\\
&= |\bar{\theta}_n^\lambda|^2 - \lambda(t-n)\bar{a}\cdot T_1^\lambda(\bar{\theta}_n^\lambda)\cdot|\bar{\theta}_n^\lambda|^2-\lambda(t-n)\cdot T_2^\lambda(\bar{\theta}_n^\lambda) + 2\lambda(t-n)(\bar{b}+\lambda(t-n)K^2)\nonumber\\
&\leq |\bar{\theta}_n^\lambda|^2 - \lambda(t-n)\bar{a}\cdot T_1^\lambda(\bar{\theta}_n^\lambda)\cdot|\bar{\theta}_n^\lambda|^2-\lambda(t-n)\cdot T_2^\lambda(\bar{\theta}_n^\lambda) + \lambda(t-n)(2\bar{b}+2K^2),\label{eqn:lemma4.2_delta_moment2_condexp_ub}
\end{align}
where the second and last inequalities follow from $(x+y)^2\leq 2(x^2+y^2)$ and $\lambda(t-n)<1$, respectively, and where for every $\theta\in \R^d$,
\begin{align}
T_1^\lambda(\theta):=\frac{|\theta|^{r}}{(1+\lambda|\theta|^{2r})^{1/2}}\label{eqn:lemma4.2_T1}
\end{align}
and 
\begin{align}
T_2^\lambda(\theta):=\frac{\bar{a}|\theta|^{r+2}}{(1+\lambda\lbrb{\theta}^{2r})^{1/2}}-\frac{2\lambda(t-n) K^2|\theta|^{2r+2}}{1+\lambda\lbrb{\theta}^{2r}}.\label{eqn:lemma4.2_T2}
\end{align}
Observe that, for every $\lambda >0$, the function $s\mapsto \frac{s}{(1+\lambda s^2)^{1/2}}$ is non-decreasing. Thus, since $\lambda < \lambda_{1,\max} \leq1$,
\begin{align}
T_1^\lambda(\theta)\geq \frac{1}{(1+\lambda)^{1/2}}>\frac{1}{\sqrt{2}}:=\kappa,\qquad \forall|\theta|>1.\label{eqn:lemma4.2_T1_lb}
\end{align}
Furthermore, applying the inequality $\sqrt{x+y}\geq \frac{1}{\sqrt{2}}(\sqrt{x}+\sqrt{y})$ for all $x,y\geq 0$ and the fact that $\lambda < \lambda_{1,\max}\leq \frac{\bar{a}^2}{8K^4}$ to the definition of $T_2^\lambda$ gives
\begin{align}
T_2^\lambda(\theta)
&=\frac{\bar{a}|\theta|^{r+2}(1+\lambda |\theta|^{2r})^{1/2}-2\lambda(t-n) K^2|\theta|^{2r+2}}{1+\lambda\lbrb{\theta}^{2r}}\nonumber\\
&\geq \frac{\sqrt{\lambda} (\frac{\bar{a}}{\sqrt{2}}-2\sqrt{\lambda} K^2)|\theta|^{2r+2}}{1+\lambda\lbrb{\theta}^{2r}}\nonumber\\
&\geq \frac{\sqrt{\lambda} (\frac{\bar{a}}{\sqrt{2}}-2\cdot \frac{\bar{a}}{2\sqrt{2}K^2}\cdot K^2)|\theta|^{2r+2}}{1+\lambda\lbrb{\theta}^{2r}}\nonumber\\
&= 0,\qquad \forall\theta\in\R^d.\label{eqn:lemma4.2_T2_lb}
\end{align}
Denote $S_n:=\{\omega\in\Omega: |\bar{\theta}_n^\lambda(\omega)| > 1\}$. Substituting (\ref{eqn:lemma4.2_T1_lb}) and (\ref{eqn:lemma4.2_T2_lb}) into (\ref{eqn:lemma4.2_delta_moment2_condexp_ub}) yields
\begin{align}
&|\Delta_{n,t}^\lambda|^2\nonumber\\
&\leq |\bar{\theta}_n^\lambda|^2 - \lambda(t-n)\bar{a}\cdot T_1^\lambda(\bar{\theta}_n^\lambda)\cdot|\bar{\theta}_n^\lambda|^2-\lambda(t-n)\cdot T_2^\lambda(\bar{\theta}_n^\lambda) + \lambda(t-n)(2\bar{b}+2K^2)\nonumber\\
&\leq |\bar{\theta}_n^\lambda|^2 - \lambda(t-n)\bar{a}\cdot (T_1^\lambda(\bar{\theta}_n^\lambda)\1_{S_n}+T_1^\lambda(\bar{\theta}_n^\lambda)\1_{S_n^c})\cdot|\bar{\theta}_n^\lambda|^2+ \lambda(t-n)(2\bar{b}+2K^2)\nonumber\\
&\leq |\bar{\theta}_n^\lambda|^2 - \lambda(t-n)\bar{a}\cdot (\kappa\1_{S_n}+T_1^\lambda(\bar{\theta}_n^\lambda)\1_{S_n^c})\cdot|\bar{\theta}_n^\lambda|^2+ \lambda(t-n)(2\bar{b}+2K^2)\nonumber\\
&\leq |\bar{\theta}_n^\lambda|^2 - \lambda(t-n)\bar{a}\cdot (\kappa\1_{S_n})\cdot|\bar{\theta}_n^\lambda|^2+ \lambda(t-n)(2\bar{b}+2K^2)\nonumber\\
&= |\bar{\theta}_n^\lambda|^2 - \lambda(t-n)\bar{a}\kappa\cdot|\bar{\theta}_n^\lambda|^2+ \lambda(t-n)\bar{a}\kappa\1_{S_n^c}\cdot|\bar{\theta}_n^\lambda|^2+ \lambda(t-n)(2\bar{b}+2K^2)\nonumber\\
&\leq (1 - \lambda(t-n)\bar{a}\kappa)|\bar{\theta}_n^\lambda|^2+ \lambda(t-n)(\bar{a}\kappa + 2\bar{b}+2K^2),\label{eqn:lemma4.2_delta_moment2_condexp_ub2}
\end{align}
where the fourth inequality follows from $T_1^\lambda(\theta)\geq 0$ for all $\theta\in\R^d$. Therefore, (\ref{eqn:lemma4.2_delta_moment2_condexp_ub2}) and (\ref{eqn:lemma4.2_tula_interpolated_moment2_condexp_eq}) together imply that, for every $t\in(n,n+1], n\in\N_0, \lambda\in (0,\lambda_{1,\max})$, 
\begin{align}
\E\lsrs{\ce{|\bar{\theta}_t^\lambda|^2} \bar{\theta}_n^\lambda}\leq (1 - \lambda(t-n)\bar{a}\kappa)|\bar{\theta}_n^\lambda|^2 + \lambda (t-n)c_0,\label{eqn:lemma4.2_tula_interpolated_moment2_condexp_ub}
\end{align}
where
\begin{align}
c_0 &:= \bar{a}\kappa+2\bar{b}+2d\beta^{-1}+2K^2,\quad \kappa := \tfrac{1}{\sqrt{2}}.\label{eqn:lemma4.2_constants_c0_kappa}
\end{align}
Taking expectations on (\ref{eqn:lemma4.2_tula_interpolated_moment2_condexp_ub}) further yields, for all $t\in(n,n+1], n\in\N_0, \lambda\in (0,\lambda_{1,\max})$,
\begin{align}
\E\lsrs{|\bar{\theta}_t^\lambda|^2} \leq (1 - \lambda(t-n)\bar{a}\kappa)\E\lsrs{|\bar{\theta}_n^\lambda|^2} + \lambda (t-n)c_0.\label{eqn:lemma4.2_tula_interpolated_moment2_ub}
\end{align}
In particular, for all $n\in\N, \lambda\in (0,\lambda_{1,\max})$, 
\begin{align}
\left(\E\lsrs{|\bar{\theta}_n^\lambda|^2} - \frac{c_0}{\bar{a}\kappa}\right) 
&\leq (1 - \lambda\bar{a}\kappa)\left(\E\lsrs{|\bar{\theta}_{n-1}^\lambda|^2} -  \frac{c_0}{\bar{a}\kappa}\right)\nonumber\\
&\leq \vdots\nonumber\\
&\leq (1-\lambda\bar{a}\kappa)^n\left(\E\lsrs{|\bar{\theta}_0^\lambda|^2} -  \frac{c_0}{\bar{a}\kappa}\right).\label{eqn:lemma4.2_tula_moment2_ub}
\end{align} 
Therefore, (\ref{eqn:lemma4.2_tula_interpolated_moment2_condexp_ub}) and (\ref{eqn:lemma4.2_tula_moment2_ub}) together imply
\begin{align}
\E\lsrs{|\bar{\theta}_t^\lambda|^2}
&\leq (1 - \lambda(t-n)\bar{a}\kappa)\lsrs{\frac{c_0}{\bar{a}\kappa}+(1-\lambda\bar{a}\kappa)^n\left(\E\lsrs{|\bar{\theta}_0^\lambda|^2}-\frac{c_0}{\bar{a}\kappa}\right)} + \lambda (t-n)c_0\nonumber\\
&\leq (1 - \lambda(t-n)\bar{a}\kappa)(1-\lambda\bar{a}\kappa)^n\E\lsrs{|\bar{\theta}_0^\lambda|^2} + c_0\left(1+\frac{1}{\bar{a}\kappa}\right),\nonumber
\end{align}
completing the proof. \qed

\subsection*{\textit{Proof of Lemma \ref{lemma:4.2}}(ii)}
For any $p\in[2,\infty)\cap\N$, $\lambda\in (0,\tilde{\lambda}_{\max})$, $t\in(n,n+1]$, $n\in\N_0$, we obtain an upper bound estimate for the $2p$-th moment $\E\lsrs{|\bar{\theta}_t^\lambda|^{2p}}$ of the continuous-time interpolated mTULA algorithm (\ref{eqn:tula_interpolated}). Recall the definitions of $\Delta_{n,t}^\lambda , \Xi_{n,t}^\lambda$ from (\ref{eqn:lemma4.2_delta_xi}). For all $x,y\in\R^d$, we have the inequality
\begin{align}
|x+y|^{2p}\leq |x|^{2p}+2p|x|^{2p-2}\lara{x,y}+\sum^{2p}_{k=2}\binom{2p}{k}|x|^{2p-k}|y|^k\label{eqn:lemma4.2_binom_inequality}
\end{align}
see, for example, Lemma A.3 of \cite{chau2021stochastic}. An application of this inequality yields, for every $p\in \N$, $t\in (n, n+1]$,
\begin{align}
&\E\lsrs{\ce{|\bar{\theta}_t^\lambda|^{2p}} \bar{\theta}_n^\lambda}\nonumber\\
&\leq |\Delta_{n,t}^\lambda|^{2p}  + 2p|\Delta_{n,t}^\lambda|^{2p-2}\E\lsrs{\ce{\lara{\Delta_{n,t}^\lambda, \Xi_{n,t}^\lambda}}\bar{\theta}_n^\lambda} + \sum^{2p}_{k=2}\binom{2p}{k}\E\lsrs{\ce{|\Delta_{n,t}^\lambda|^{2p-k}|\Xi_{n,t}^\lambda|^{k}} \bar{\theta}_n^\lambda}\nonumber\\
&=|\Delta_{n,t}^\lambda|^{2p}  + \sum^{2p-2}_{k=0}\binom{2p}{k+2}\E\lsrs{\ce{|\Delta_{n,t}^\lambda|^{2p-2-k}|\Xi_{n,t}^\lambda|^{k}|\Xi_{n,t}^\lambda|^{2}} \bar{\theta}_n^\lambda}\nonumber\\
&=|\Delta_{n,t}^\lambda|^{2p}  + \sum^{2p-2}_{k=0}\frac{2p(2p-1)}{(k+2)(k+1)}\binom{2p-2}{k}\E\lsrs{\ce{|\Delta_{n,t}^\lambda|^{2p-2-k}|\Xi_{n,t}^\lambda|^{k}|\Xi_{n,t}^\lambda|^{2}} \bar{\theta}_n^\lambda}\nonumber\\
&\leq |\Delta_{n,t}^\lambda|^{2p}  + p(2p-1)\sum^{2p-2}_{k=0}\binom{2p-2}{k}\E\lsrs{\ce{|\Delta_{n,t}^\lambda|^{2p-2-k}|\Xi_{n,t}^\lambda|^{k}|\Xi_{n,t}^\lambda|^{2}} \bar{\theta}_n^\lambda}\nonumber\\
&= |\Delta_{n,t}^\lambda|^{2p}  + p(2p-1)\E\lsrs{\ce{(|\Delta_{n,t}^\lambda|+|\Xi_{n,t}^\lambda|)^{2p-2}|\Xi_{n,t}^\lambda|^{2}} \bar{\theta}_n^\lambda}\nonumber\\
&\leq |\Delta_{n,t}^\lambda|^{2p}  + p(2p-1)2^{2p-3}\E\lsrs{\ce{(|\Delta_{n,t}^\lambda|^{2p-2}+|\Xi_{n,t}^\lambda|^{2p-2})|\Xi_{n,t}^\lambda|^{2}} \bar{\theta}_n^\lambda}\nonumber\\
&= |\Delta_{n,t}^\lambda|^{2p}  + p(2p-1)2^{2p-3}|\Delta_{n,t}^\lambda|^{2p-2}\E\lsrs{|\Xi_{n,t}^\lambda|^{2}} + p(2p-1)2^{2p-3}\E\lsrs{|\Xi_{n,t}^\lambda|^{2p}}\nonumber\\
&=|\Delta_{n,t}^\lambda|^{2p} + p(2p-1)2^{2p-2}\lambda\beta^{-1}(t-n)d\  |\Delta_{n,t}^\lambda|^{2p-2}\nonumber\\
&\quad + p(2p-1)2^{4p-3}(\lambda \beta^{-1} (t-n))^{p}p!\binom{\frac{d}{2}+p-1}{p}.\label{eqn:lemma4.2_tula_interpolated_moment2p_condexp_ub}
\end{align}
We have previously established in the proof of Lemma \ref{lemma:4.2}(i) that, for $\lambda \in \left(0,\min\{1,\frac{\bar{a}^2}{8K^4}\}\right)$, 
\begin{align}
|\Delta_{n,t}^\lambda|^2 \leq \left(1-\lambda(t-n)\frac{\bar{a}|\bar{\theta}_n^\lambda|^{r}}{(1+\lambda |\bar{\theta}_n^\lambda|^{2r})^{1/2}}\right)|\bar{\theta}_n^\lambda|^2 + \lambda(t-n)(2\bar{b}+2K^2),\label{eqn:lemma4.2_delta_moment2_condexp_ub3}
\end{align}
see (\ref{eqn:lemma4.2_delta_moment2_condexp_ub}) and (\ref{eqn:lemma4.2_T2_lb}). Furthermore, for every $\lambda \in (0, \frac{1}{\bar{a}^2}]$, 
\begin{align}
1\geq 1-\lambda(t-n)\frac{\bar{a}|\bar{\theta}_n^\lambda|^{r}}{(1+\lambda |\bar{\theta}_n^\lambda|^{2r})^{1/2}}
&> 1 - \lambda\frac{\bar{a}|\bar{\theta}_n^\lambda|^{r}}{(\lambda |\bar{\theta}_n^\lambda|^{2r})^{1/2}}\nonumber \\
&= 1 - \sqrt{\lambda}\bar{a}\nonumber\\
& \geq 0.\nonumber
\end{align}
Hence, raising (\ref{eqn:lemma4.2_delta_moment2_condexp_ub3}) to the $p$-th power yields, for $\lambda< \tilde{\lambda}_{\max}= \min\{1,\frac{\bar{a}^2}{8K^4}, \frac{1}{\bar{a}^2}\}$,
\begin{align}
|\Delta_{n,t}^\lambda|^{2p}
&=\sum^p_{k=0}\binom{p}{k}\lsrs{\left(1-\lambda(t-n)\frac{\bar{a}|\bar{\theta}_n^\lambda|^{r}}{(1+\lambda |\bar{\theta}_n^\lambda|^{2r})^{1/2}}\right)|\bar{\theta}_n^\lambda|^2}^{p-k}\lambda^k(t-n)^k(2\bar{b}+2K^2)^k\nonumber\\
&\leq \left(1-\lambda(t-n)\frac{\bar{a}|\bar{\theta}_n^\lambda|^{r}}{(1+\lambda |\bar{\theta}_n^\lambda|^{2r})^{1/2}}\right)|\bar{\theta}_n^\lambda|^{2p} + \sum^p_{k=1}\binom{p}{k}\lambda^k(t-n)^k(2\bar{b}+2K^2)^k|\bar{\theta}_n^\lambda|^{2p-2k}\nonumber\\
& =\left(1-\lambda(t-n)\frac{\bar{a}|\bar{\theta}_n^\lambda|^{r}}{2(1+\lambda |\bar{\theta}_n^\lambda|^{2r})^{1/2}}\right)|\bar{\theta}_n^\lambda|^{2p} - J_1^\lambda(\bar{\theta}_n^\lambda),\label{eqn:lemma4.2_delta_moment2p_condexp_ub}
\end{align}
where we define, for every $\theta\in\R^d$,
\begin{align}
J_1^\lambda(\theta):=\lambda(t-n)\frac{\bar{a}|\theta|^{r+2p}}{2(1+\lambda |\theta|^{2r})^{1/2}}-\sum^p_{k=1}\binom{p}{k}\lambda^k(t-n)^k(2\bar{b}+2K^2)^k|\theta|^{2p-2k}.\label{eqn:lemma4.2_J1}
\end{align}
For any $\theta\in\R^d$ and any $\lambda \in (0,\tilde{\lambda}_{\max})$, $J_1^\lambda(\theta)$ can be bounded below by
\begin{align}
&J_1^\lambda(\theta)\nonumber\\
&\geq \lambda(t-n)\frac{\bar{a}|\theta|^{r+2p}}{2(1+\lambda |\theta|^{2r})^{1/2}}-\sum^p_{k=1}\binom{p}{k}\lambda^k(t-n)^k(1+2\bar{b}+2K^2)^k|\theta|^{2p-2k}\nonumber\\
&\geq \lambda(t-n)\frac{\bar{a}|\theta|^{r+2p}}{2(1+\lambda |\theta|^{2r})^{1/2}}-\binom{p}{\lcrc{p/2}}\lambda(t-n)(1+2\bar{b}+2K^2)^p\sum^p_{k=1}|\theta|^{2p-2k}\cdot\frac{(1+\lambda |\theta|^{2r})^{1/2}}{(1+\lambda |\theta|^{2r})^{1/2}}\nonumber\\
&\geq \lambda(t-n)\frac{\bar{a}|\theta|^{r+2p}}{2(1+\lambda |\theta|^{2r})^{1/2}}-\binom{p}{\lcrc{p/2}}\lambda(t-n)(1+2\bar{b}+2K^2)^p\sum^p_{k=1}|\theta|^{2p-2k}\cdot\frac{1+|\theta|^r}{(1+\lambda |\theta|^{2r})^{1/2}}\nonumber\\
&=\frac{\lambda (t-n)}{(1+\lambda |\theta|^{2r})^{1/2}}\sum^p_{k=1}\lsrs{\frac{\bar{a}}{4p}|\theta|^{r+2p} - \binom{p}{\lcrc{p/2}}(1+2\bar{b}+2K^2)^p |\theta|^{2p-2k}}\nonumber\\
&\quad + \frac{\lambda (t-n)}{(1+\lambda |\theta|^{2r})^{1/2}}\sum^p_{k=1}\lsrs{\frac{\bar{a}}{4p}|\theta|^{r+2p} - \binom{p}{\lcrc{p/2}}(1+2\bar{b}+2K^2)^p |\theta|^{r+2p-2k}}.\label{eqn:lemma4.2_J1_lb}
\end{align}
Observe that, for $1\leq k \leq p$,
\begin{align}
\frac{\bar{a}}{4p}|\theta|^{r+2p} - \binom{p}{\lcrc{p/2}}(1+2\bar{b}+2K^2)^p |\theta|^{2p-2k} > 0 \iff |\theta|>\lsrs{\frac{4p\binom{p}{\lcrc{p/2}}(1+2\bar{b}+2K^2)^p}{\bar{a}}}^{\frac{1}{r+2k}}\label{eqn:lemma4.2_J1_piece1_lb}
\end{align}
and
\begin{align}
\frac{\bar{a}}{4p}|\theta|^{r+2p} - \binom{p}{\lcrc{p/2}}(1+2\bar{b}+2K^2)^p |\theta|^{r+2p-2k} > 0  \iff |\theta|>\lsrs{\frac{4p\binom{p}{\lcrc{p/2}}(1+2\bar{b}+2K^2)^p}{\bar{a}}}^{\frac{1}{2k}}.\label{eqn:lemma4.2_J1_piece2_lb}
\end{align}
Therefore, setting
\begin{align}
M_1(p):= \frac{4p\binom{p}{\lcrc{p/2}}(1+2\bar{b}+2K^2)^p}{\min\{1,\bar{a}\}},\qquad \forall p\in\N,\label{eqn:lemma4.2_constant_M1}
\end{align}
we have
\begin{align}
M_1(p)\geq \max_{1\leq k\leq p}\max\left\{\lsrs{\frac{4p\binom{p}{\lcrc{p/2}}(1+2\bar{b}+2K^2)^p}{\bar{a}}}^{\frac{1}{r+2k}}, \lsrs{\frac{4p\binom{p}{\lcrc{p/2}}(1+2\bar{b}+2K^2)^p}{\bar{a}}}^{\frac{1}{2k}}\right\}.\label{eqn:lemma4.2_M1_lb}
\end{align}
It follows from (\ref{eqn:lemma4.2_J1_lb}), (\ref{eqn:lemma4.2_J1_piece1_lb}), (\ref{eqn:lemma4.2_J1_piece2_lb}) and (\ref{eqn:lemma4.2_M1_lb}) that, for every $\lambda \in (0,\tilde{\lambda}_{\max})$ and $\theta\in\R^d$ such that $|\theta| > M_1(p)$,
\begin{align}
J_1^\lambda(\theta)>0.\label{eqn:lemma4.2_J1_lb2}
\end{align}
Therefore, denoting $S_{n,M_1(p)}:=\{\omega\in\Omega: |\bar{\theta}_n^\lambda(\omega)| > M_1(p)\}$, we have from (\ref{eqn:lemma4.2_delta_moment2p_condexp_ub}) and (\ref{eqn:lemma4.2_J1_lb2}), the fact that $s\mapsto \frac{s}{(1+\lambda s^2)^{1/2}}$ is non-decreasing on $[0,\infty)$ and $0<\lambda < \tilde{\lambda}_{\max} \leq 1$  that
\begin{align}
|\Delta_{n,t}^\lambda|^{2p}\1_{S_{n,M_1(p)}} 
&\leq\left(1-\lambda(t-n)\frac{\bar{a}|\bar{\theta}_n^\lambda|^{r}}{2(1+\lambda |\bar{\theta}_n^\lambda|^{2r})^{1/2}}\right)|\bar{\theta}_n^\lambda|^{2p} \1_{S_{n,M_1(p)}}\nonumber\\
&\leq \left(1-\lambda(t-n)\bar{a}\tilde{\kappa}(p)\right)|\bar{\theta}_n^\lambda|^{2p} \1_{S_{n,M_1(p)}},\label{eqn:lemma4.2_delta_moment2p_condexp_set_ub}
\end{align}
where
\begin{align}
\tilde{\kappa}(p):= \frac{[M_1(p)]^r}{2(1+[M_1(p)]^{2r})^{1/2}},\qquad \forall p\in\N.\label{eqn:lemma4.2_constant_kappatilde}
\end{align}
Furthermore, for every $\lambda\in (0,\tilde{\lambda}_{\max})$, we have
\begin{align}
&|\Delta_{n,t}^\lambda|^{2p}\1_{S_{n,M_1(p)}^c} \nonumber\\
&\leq|\bar{\theta}_n^\lambda|^{2p} \1_{S_{n,M_1(p)}^c}-J_1^\lambda(\bar{\theta}_n^\lambda)\1_{S_{n,M_1(p)}^c} \nonumber\\
&= \left(1-\lambda(t-n)\bar{a}\tilde{\kappa}(p)\right)|\bar{\theta}_n^\lambda|^{2p} \1_{S_{n,M_1(p)}^c} + \lambda (t-n)\bar{a}\tilde{\kappa}(p)|\bar{\theta}_n^\lambda|^{2p} \1_{S_{n,M_1(p)}^c} -J_1^\lambda(\bar{\theta}_n^\lambda)\1_{S_{n,M_1(p)}^c}\nonumber\\
&\leq \left(1-\lambda(t-n)\bar{a}\tilde{\kappa}(p)\right)|\bar{\theta}_n^\lambda|^{2p} \1_{S_{n,M_1(p)}^c} + \lambda (t-n)\bar{a}\tilde{\kappa}(p)[M_1(p)]^{2p} \1_{S_{n,M_1(p)}^c}\nonumber\\
&\quad + \lambda(t-n)\sum^p_{k=1}\binom{p}{k}(2\bar{b}+2K^2)^k [M_1(p)]^{2p-2k}\1_{S_{n,M_1(p)}^c}\nonumber\\
&=\left(1-\lambda(t-n)\bar{a}\tilde{\kappa}(p)\right)|\bar{\theta}_n^\lambda|^{2p} \1_{S_{n,M_1(p)}^c} + \lambda(t-n)c_1(p)\1_{S_{n,M_1(p)}^c},\label{eqn:lemma4.2_delta_moment2p_condexp_setcomp_ub}
\end{align}
where, for every $p\in\N$,
\begin{align}
c_1(p):=\bar{a}\tilde{\kappa}(p)[M_1(p)]^{2p}+\sum^p_{k=1}\binom{p}{k}(2\bar{b}+2K^2)^k [M_1(p)]^{2p-2k}.\label{eqn:lemma4.2_constant_c1}
\end{align}
Combining the inequalities (\ref{eqn:lemma4.2_delta_moment2p_condexp_set_ub}) and (\ref{eqn:lemma4.2_delta_moment2p_condexp_setcomp_ub}) yields, for $p\in \N$, $\lambda\in (0,\tilde{\lambda}_{\max})$,
\begin{align}
|\Delta_{n,t}^\lambda|^{2p}\leq \left(1-\lambda(t-n)\bar{a}\tilde{\kappa}(p)\right)|\bar{\theta}_n^\lambda|^{2p}  + \lambda(t-n)c_1(p).\label{eqn:lemma4.2_delta_moment2p_condexp_ub2}
\end{align}
By (\ref{eqn:lemma4.2_delta_moment2p_condexp_ub2}) and the fact that $\bar{a}\kappa + 2\bar{b}+2K^2<c_1(1)$, we also have, for every $p\in [2,\infty)\cap\N$, $\lambda\in (0,\tilde{\lambda}_{\max})$,
\begin{align}
|\Delta_{n,t}^\lambda|^{2p-2}\leq |\bar{\theta}_n^\lambda|^{2p-2}  + \lambda(t-n)c_1(p-1).\label{eqn:lemma4.2_delta_moment2pminus2_condexp_ub}
\end{align}
Substituting (\ref{eqn:lemma4.2_delta_moment2p_condexp_ub2}) and (\ref{eqn:lemma4.2_delta_moment2pminus2_condexp_ub}) into (\ref{eqn:lemma4.2_tula_interpolated_moment2p_condexp_ub}) yields, for every $p\in [2,\infty)\cap\N$, $\lambda\in (0,\tilde{\lambda}_{\max})$, $t\in (n,n+1]$, $n\in\N_0$,
\begin{align}
&\E\lsrs{\ce{|\bar{\theta}_t^\lambda|^{2p}} \bar{\theta}_n^\lambda}\nonumber\\
&\leq \left(1-\lambda(t-n)\bar{a}\tilde{\kappa}(p)\right)|\bar{\theta}_n^\lambda|^{2p} + p(2p-1)2^{2p-2}\lambda\beta^{-1}(t-n)d\  |\bar{\theta}_n^\lambda|^{2p-2} \\
&\quad +  p(2p-1)2^{2p-2}\lambda^2\beta^{-1}(t-n)^2d\ c_1(p-1) + \lambda(t-n)c_1(p)\nonumber\\
&\quad + p(2p-1)2^{4p-3}(\lambda \beta^{-1} (t-n))^{p}p!\binom{\frac{d}{2}+p-1}{p}\nonumber\\
&\leq \left(1-\lambda(t-n)\bar{a}\tilde{\kappa}(p)\right)|\bar{\theta}_n^\lambda|^{2p} + p(2p-1)2^{2p-2}\lambda\beta^{-1}(t-n)d\  |\bar{\theta}_n^\lambda|^{2p-2} + \lambda(t-n)c_2(p),\label{eqn:lemma4.2_tula_interpolated_moment2p_condexp_ub2}
\end{align}
where 
\begin{align}
c_2(p):=c_1(p) + p(2p-1)2^{2p-2}\beta^{-1}d\ c_1(p-1)+ p(2p-1)2^{4p-3}\beta^{-p}p!\binom{\frac{d}{2}+p-1}{p}.\label{eqn:lemma4.2_constant_c2}
\end{align}
Observe that, for every $\theta\in\R^d$,
\begin{align}
\frac{\bar{a}\tilde{\kappa}(p)}{2}|\theta|^{2p}-p(2p-1)2^{2p-2}\beta^{-1}d|\theta|^{2p-2} > 0 \iff |\theta| > M_2(p),\label{eqn:lemma4.2_tula_interpolated_moment2p_condexp_piece1_lb}
\end{align}
where
\begin{align}
M_2(p):=\lsrs{\frac{p(2p-1)2^{2p-1}\beta^{-1}d}{\bar{a}\tilde{\kappa}(p)}}^{1/2}.\label{eqn:lemma4.2_constant_M2}
\end{align}
Hence, by denoting $S_{n,M_2(p)}:=\{\omega\in\Omega: |\bar{\theta}_n^\lambda(\omega)| > M_2(p)\}$, we have by (\ref{eqn:lemma4.2_tula_interpolated_moment2p_condexp_ub2}) and (\ref{eqn:lemma4.2_tula_interpolated_moment2p_condexp_piece1_lb}),
\begin{align}
\E\lsrs{\ce{|\bar{\theta}_t^\lambda|^{2p}\1_{S_{n,M_2(p)}}} \bar{\theta}_n^\lambda}
&\leq \left(1-\lambda(t-n)\frac{\bar{a}\tilde{\kappa}(p)}{2}\right)|\bar{\theta}_n^\lambda|^{2p}\1_{S_{n,M_2(p)}}+ \lambda(t-n)c_2(p)\1_{S_{n,M_2(p)}}.\label{eqn:lemma4.2_tula_interpolated_moment2p_condexp_set_ub}
\end{align}
Furthermore,
\begin{align}
\E\lsrs{\ce{|\bar{\theta}_t^\lambda|^{2p}\1_{S_{n,M_2(p)}^c}} \bar{\theta}_n^\lambda}
&\leq \left(1-\lambda(t-n)\frac{\bar{a}\tilde{\kappa}(p)}{2}\right)|\bar{\theta}_n^\lambda|^{2p}\1_{S_{n,M_2(p)}^c}+ \lambda(t-n)c_3(p)\1_{S_{n,M_2(p)}^c},\label{eqn:lemma4.2_tula_interpolated_moment2p_condexp_setcomp_ub}
\end{align}
where
\begin{align}
c_3(p):=c_2(p)+p(2p-1)2^{2p-2}\beta^{-1}d\ [M_2(p)]^{2p-2}.\label{eqn:lemma4.2_constant_c3}
\end{align}
Combining (\ref{eqn:lemma4.2_tula_interpolated_moment2p_condexp_set_ub}) and (\ref{eqn:lemma4.2_tula_interpolated_moment2p_condexp_setcomp_ub}) yields, for $p\in [2,\infty)\cap\N$, $\lambda\in (0,\tilde{\lambda}_{\max})$,
\begin{align}
\E\lsrs{\ce{|\bar{\theta}_t^\lambda|^{2p}} \bar{\theta}_n^\lambda}
&\leq \left(1-\lambda(t-n)\frac{\bar{a}\tilde{\kappa}(p)}{2}\right)|\bar{\theta}_n^\lambda|^{2p}+ \lambda(t-n)c_3(p),\label{eqn:lemma4.2_tula_interpolated_moment2p_condexp_ub3}
\end{align}
where, for $p\in [2,\infty)\cap\N$, we define the constants
\begin{align}
\tilde{\kappa}(p) &:= \frac{[M_1(p)]^r}{2(1+[M_1(p)]^{2r})^{1/2}},\nonumber\\
M_1(p)&:= \frac{4p\binom{p}{\lcrc{p/2}}(1+2\bar{b}+2K^2)^p}{\min\{1,\bar{a}\}},\nonumber\\
c_3(p)&:=c_2(p)+p(2p-1)2^{2p-2}\beta^{-1}d\ [M_2(p)]^{2p-2},\nonumber\\
M_2(p)&:=\lsrs{\frac{p(2p-1)2^{2p-1}\beta^{-1}d}{\bar{a}\tilde{\kappa}(p)}}^{1/2},\nonumber\\
c_2(p)&:=c_1(p) + p(2p-1)2^{2p-2}\beta^{-1}d\ c_1(p-1)+ p(2p-1)2^{4p-3}\beta^{-p}p!\binom{\frac{d}{2}+p-1}{p},\nonumber\\
c_1(p)&:=\bar{a}\tilde{\kappa}(p)[M_1(p)]^{2p}+\sum^p_{k=1}\binom{p}{k}(2\bar{b}+2K^2)^k [M_1(p)]^{2p-2k}.\label{eqn:lemma4.2_constant_kappatilde_M1_c3_M2_c2_c1}
\end{align}
Observing that $\tilde{\kappa}(p) \geq \tilde{\kappa}(2)$ for $p\in [2,\infty)\cap\N$ and applying the same argument as in the proof of Lemma \ref{lemma:4.2}(i), we obtain from (\ref{eqn:lemma4.2_tula_interpolated_moment2p_condexp_ub3}), for $p\in [2,\infty)\cap\N$, $\lambda\in (0,\tilde{\lambda}_{\max})$, $t\in(n,n+1]$, $n\in \N_0$, the result
\begin{align}
\E\lsrs{|\bar{\theta}_t^\lambda|^{2p}}\leq \left(1-\lambda(t-n)\frac{\bar{a}\tilde{\kappa}(2)}{2}\right)\left(1-\lambda\frac{\bar{a}\tilde{\kappa}(2)}{2}\right)^n\E\lsrs{|\bar{\theta}_0^\lambda|^{2p}}+c_3(p)\left(1+\frac{2}{\bar{a}\tilde{\kappa}(2)}\right).\nonumber
\end{align}
This completes the proof. \qed

\subsection*{\textit{Proof of Lemma \ref{lemma:4.4}}} 
\label{proof:lemma4.4}
For every $p\in [2,\infty)\cap \N$, $\theta\in\R^d$, one obtains, by direct computation,
\begin{align}
\nabla V_p(\theta) = p V_{p-2}(\theta)\cdot \theta^T\label{eqn:lemma4.4_Vp_gradient}
\end{align}
and
\begin{align}
\Delta V_p(\theta) 
&= \tr\lsrs{\nabla\nabla^T V_p(\theta)}\nonumber\\
&= \tr\lsrs{\nabla\left(p V_{p-2}(\theta)\cdot \theta\right)}\nonumber\\
&= \tr\lsrs{p(p-2)(1+|\theta|^2)^{p/2-2}\theta\theta^T + p V_{p-2}(\theta)I_d}\nonumber\\
&= p(p-2)(1+|\theta|^2)^{p/2-2}|\theta|^2 + dp V_{p-2}(\theta).\label{eqn:lemma4.4_Vp_laplacian}
\end{align}

Combining (\ref{eqn:lemma4.4_Vp_gradient}) and (\ref{eqn:lemma4.4_Vp_laplacian}) then applying Remark \ref{remark:2.4} yields
\begin{align}
&\beta^{-1}\Delta V_p(\theta) - \lara{\nabla V_p(\theta), h(\theta)}\nonumber\\
=\ & \beta^{-1}p(p-2)(1+|\theta|^2)^{p/2-2}|\theta|^2 + \beta^{-1}dp V_{p-2}(\theta) - pV_{p-2}(\theta)\lara{\theta, h(\theta)}\nonumber \\
\leq\ & \beta^{-1}p(p-2)(1+|\theta|^2)^{p/2-2}|\theta|^2 + \beta^{-1}dp V_{p-2}(\theta) - \bar{a}pV_{p-2}(\theta)|\theta|^2 + \bar{b}' pV_{p-2}(\theta)\nonumber\\
=\ & - \bar{a}pV_p(\theta) + (\bar{a} + \bar{b}' + \beta^{-1} d)pV_{p-2}(\theta) + \beta^{-1}p(p-2)(1+|\theta|^2)^{p/2-2}|\theta|^2\nonumber\\
\leq\ & - \bar{a}pV_p(\theta) + (\bar{a} + \bar{b}' + \beta^{-1}(d+p-2))pV_{p-2}(\theta)\nonumber\\
=\ & -\frac{\bar{a} p}{2}V_p(\theta) - p\lsrs{\frac{\bar{a}}{2}V_p(\theta) - (\bar{a} + \bar{b}' + \beta^{-1}(d+p-2))V_{p-2}(\theta)}.\label{eqn:lemma4.4_drift_ub}
\end{align}

Define, for every $p\in [2,\infty)\cap\N$, the constant
\begin{align}
M_V(p):=\left(1+\frac{2\bar{b}'+2\beta^{-1}(d+p-2)}{\bar{a}}\right)^{1/2}.\label{eqn:lemma4.4_constant_MV}
\end{align}

Observe that
\begin{align}
\lsrs{\frac{\bar{a}}{2}V_p(\theta) - (\bar{a} + \bar{b}' + \beta^{-1}(d+p-2))V_{p-2}(\theta)} > 0 \iff |\theta| > M_V(p).\label{eqn:lemma4.4_drift_piece1_lb}
\end{align}

Hence, from (\ref{eqn:lemma4.4_drift_ub}) and (\ref{eqn:lemma4.4_drift_piece1_lb}), one obtains
\begin{align}
&\beta^{-1}\Delta V_p(\theta) - \lara{\nabla V_p(\theta), h(\theta)}\nonumber\\
\leq\ &-\frac{\bar{a} p}{2}V_p(\theta) - p\lsrs{\frac{\bar{a}}{2}V_p(\theta) - (\bar{a} + \bar{b}' + \beta^{-1}(d+p-2))V_{p-2}(\theta)}\1_{\{|\theta| \leq M_V(p)\}}\nonumber\\
\leq\ &-\frac{\bar{a} p}{2}V_p(\theta) + p (\bar{a} + \bar{b}' + \beta^{-1}(d+p-2))V_{p-2}(\theta)\1_{\{|\theta| \leq M_V(p)\}}\nonumber\\
\leq\ &-\frac{\bar{a} p}{2}V_p(\theta) + p (\bar{a} + \bar{b}' + \beta^{-1}(d+p-2))v_{p-2}(M_V(p))\nonumber\\
=\ &-\frac{\bar{a} p}{2}V_p(\theta) + \frac{\bar{a}p}{2}v_p(M_V(p))\nonumber\\
=\ &-c_{V,1}(p)V_p(\theta) + c_{V,2}(p),\qquad \forall \theta\in\R^d,\nonumber
\end{align}
with $c_{V,1}(p):= \frac{\bar{a} p}{2}$ and $c_{V,2}(p):= \frac{\bar{a}p}{2}v_p(M_V(p))$. This completes the proof. \qed

\subsection*{\textit{Proof of Lemma \ref{lemma:4.5}}} 
\label{proof:lemma4.5}
Let $p\in[2,\infty)\cap\N$, $\lambda\in (0,\tilde{\lambda}_{\max})$, $n\in \N$ and $t\in (nT, (n+1)T]$. By Ito's formula, one obtains
\begin{align}
\E\lsrs{V_p(\bar{\zeta}_t^{\lambda,n})}
=&\ \E\lsrs{V_p(\bar{\zeta}_{nT}^{\lambda,n})} + \lambda\int^t_{nT}\E\lsrs{\beta^{-1}\Delta V_p(\bar{\zeta}_s^{\lambda,n}) - \lara{\nabla V_p(\bar{\zeta}_s^{\lambda,n}), h(\bar{\zeta}_s^{\lambda,n})}}\dee s\nonumber\\
& + \E\lsrs{\int^t_{nT}\lara{\nabla V_p(\bar{\zeta}_s^{\lambda,n}), \sqrt{2\lambda\beta^{-1}}\dee B_s^\lambda}}\label{eqn:lemma4.5_aux_momentVp_ito}\\
=&\ \E\lsrs{V_p(\bar{\theta}_{nT}^\lambda)} + \lambda\int^t_{nT}\E\lsrs{\beta^{-1}\Delta V_p(\bar{\zeta}_s^{\lambda,n}) - \lara{\nabla V_p(\bar{\zeta}_s^{\lambda,n}), h(\bar{\zeta}_s^{\lambda,n})}}\dee s.\label{eqn:lemma4.5_aux_momentVp_eq}
\end{align}

To see that the last term of (\ref{eqn:lemma4.5_aux_momentVp_ito}) vanishes, it suffices to show that $\E\lsrs{\int^t_{nT}|\nabla V_p(\bar{\zeta}_s^{\lambda, n})|^2\ \dee s} < \infty$. To this end, define the stopping time $\tau_k := \inf\{s\geq nT: |\bar{\zeta}_s^{\lambda,n}| > k\}$. By Ito's formula applied to the stopped process $V_{2p}(\bar{\zeta}_{t\wedge \tau_k}^{\lambda,n})$ and Lemma \ref{lemma:4.4},
\begin{align}
\E\lsrs{V_{2p}(\bar{\zeta}_{t\wedge \tau_k}^{\lambda,n})}
=&\ \E\lsrs{V_{2p}(\bar{\theta}_{nT}^\lambda)} + \lambda\int^t_{nT}\E\lsrs{\beta^{-1}\Delta V_{2p}(\bar{\zeta}_{s\wedge \tau_k}^{\lambda,n}) - \lara{\nabla V_{2p}(\bar{\zeta}_{s\wedge \tau_k}^{\lambda,n}), h(\bar{\zeta}_{s\wedge\tau_k}^{\lambda,n})}}\dee s\nonumber\\
\leq &\ \E\lsrs{V_{2p}(\bar{\theta}_{nT}^\lambda)} + \lambda\int^t_{nT} \lsrs{-c_{V,1}(2p)\E\lsrs{V_{2p}(\bar{\zeta}_{s\wedge \tau_k}^{\lambda,n})} + c_{V,2}(2p)}\dee s\nonumber\\
\leq &\ C_1^* + C_2^*\int^t_{nT}\E\lsrs{V_{2p}(\bar{\zeta}_{s\wedge \tau_k}^{\lambda,n})}\ \dee s,\label{eqn:lemma4.5_aux_stopped_momentV2p_ub}
\end{align}

where $C_1^*, C_2^*$ are some non-negative constants which do not depend on $\tau_k$. Applying Fatou's Lemma, then Gronwall's Lemma to (\ref{eqn:lemma4.5_aux_stopped_momentV2p_ub}) yields
\begin{align}
\E\lsrs{V_{2p}(\bar{\zeta}_t^{\lambda,n})}\leq \liminf_{k\to\infty}\E\lsrs{V_{2p}(\bar{\zeta}_{t\wedge \tau_k}^{\lambda,n})} = C_1^* e^{C_2^*(t-nT)}.\label{eqn:lemma4.5_aux_momentV2p_ub}
\end{align}

From (\ref{eqn:lemma4.4_Vp_gradient}), one sees that $|\nabla V_p(\theta)|^2 \leq p^2 V_{2p}(\theta)$ for all $\theta\in \R^d$. It follows that 
\begin{align}
\E\lsrs{\int^t_{nT}|\nabla V_p(\bar{\zeta}_s^{\lambda, n})|^2\ \dee s}\leq p^2 \int^t_{nT} \E\lsrs{V_{2p}(\bar{\zeta}_s^{\lambda,n})}\ \dee s \leq p^2 \int^t_{nT} C_1^* e^{C_2^*(s-nT)}\ \dee s <\infty,\nonumber
\end{align}

thus establishing the veracity of (\ref{eqn:lemma4.5_aux_momentVp_eq}). By differentiating both sides of (\ref{eqn:lemma4.5_aux_momentVp_eq}) and then applying Lemma \ref{lemma:4.4}, one obtains
\begin{align}
\frac{\dee}{\dee t} \E\lsrs{V_p(\bar{\zeta}_t^{\lambda,n})} 
&= \lambda \E\lsrs{\beta^{-1}\Delta V_p(\bar{\zeta}_t^{\lambda,n}) - \lara{\nabla V_p(\bar{\zeta}_t^{\lambda,n}), h(\bar{\zeta}_t^{\lambda,n})}}\nonumber\\
&\leq -\lambda c_{V,1}(p)\E\lsrs{V_{p}(\bar{\zeta}_t^{\lambda, n})} + \lambda c_{V,2}(p)  = -\lambda c_{V,1}(p)\left(\E\lsrs{V_{p}(\bar{\zeta}_t^{\lambda, n})}-v_p(M_V(p))\right).\label{eqn:lemma4.5_aux_momentVp_derivative_ub}
\end{align}

The differential form of Gronwall's inequality applied to (\ref{eqn:lemma4.5_aux_momentVp_derivative_ub}) then yields
\begin{align}
\E\lsrs{V_p(\bar{\zeta}_t^{\lambda,n})}
&\leq v_p(M_V(p)) + \left(\E\lsrs{V_p(\bar{\theta}_{nT}^\lambda)}-v_p(M_V(p))\right)\exp\left\{\int^t_{nT}-\lambda c_{V,1}(p)\ \dee s\right\}\nonumber\\
&=e^{-\lambda \bar{a}p(t-nT)/2}\E\lsrs{V_p(\bar{\theta}_{nT}^\lambda)}+v_p(M_V(p))(1-e^{-\lambda \bar{a}p(t-nT)/2}),\nonumber
\end{align}
which completes the proof. \qed

\subsection*{\textit{Proof of Corollary \ref{corollary:4.6}}}
\label{proof:corollary4.6}
Let $p\in\N$, $\lambda \in (0, \tilde{\lambda}_{\max})$, $n\in\N$, and $t\in (nT, (n+1)T]$. From Corollary \ref{corollary:4.3} and Lemma \ref{lemma:4.5}, one obtains
\begin{align}
\E\lsrs{V_{2p}(\bar{\zeta}_t^{\lambda,n})}
\leq &\ e^{-\lambda \bar{a}p(t-nT)}\E\lsrs{\left(1+|\bar{\theta}_{nT}^\lambda|^2\right)^p}+v_{2p}(M_V(2p))\nonumber\\
\leq &\ 2^{p-1}e^{-\lambda \bar{a}p(t-nT)}\left(1+\E\lsrs{|\bar{\theta}_{nT}^\lambda|^{2p}}\right)+v_{2p}(M_V(2p))\nonumber\\
\leq &\ 2^{p-1}e^{-\lambda \bar{a}p(t-nT)}\left(1+e^{-\lambda\bar{a}\kappa_*nT}\E\lsrs{|\bar{\theta}_0^\lambda|^{2p}} + c_*(p)\left(1+\frac{1}{\bar{a}\kappa_*}\right)\right) +v_{2p}(M_V(2p))\nonumber\\
\leq &\ 2^{p-1} + 2^{p-1}e^{-\lambda \bar{a}\min\{p,\kappa_*\}t}\E\lsrs{|\bar{\theta}_0^{\lambda}|^{2p}} + 2^{p-1} c_*(p)\left(1+\frac{1}{\bar{a}\kappa_*}\right) + v_{2p}(M_V(2p)),\nonumber
\end{align}
completing the proof. \qed

\subsection*{\textit{Proof of Lemma \ref{lemma:4.8}}}
\label{proof:lemma4.8}
Let $\lambda \in (0, \tilde{\lambda}_{\max})$, $n\in\N$, and $t\in (nT, (n+1)T]$. From Assumption \ref{assumption:2} and Corollary \ref{corollary:4.3}, one obtains
\begin{align}
\E\lsrs{\lbrb{\bar{\theta}_t^\lambda - \bar{\theta}_{\lfrf{t}}^\lambda}^8}
&=\E\lsrs{\lbrb{\int^t_{\lfrf{t}}-\lambda h_\lambda(\bar{\theta}_{\lfrf{s}}^\lambda)\ \dee s + \sqrt{2\lambda\beta^{-1}}(B_t^\lambda - B_{\lfrf{t}}^\lambda)}^8}\nonumber\\
&\leq 128\lambda^8\E\lsrs{\lbrb{\int^t_{\lfrf{t}} h_\lambda(\bar{\theta}_{\lfrf{s}}^\lambda)\ \dee s}^8} + 128\ \E\lsrs{\lbrb{\sqrt{2\lambda \beta^{-1}}(B_t^\lambda-B_{\lfrf{t}}^\lambda)}^8}\nonumber\\
&= 128\lambda^8\E\lsrs{\lbrb{h_\lambda(\bar{\theta}_{\lfrf{t}}^\lambda)\cdot (t-\lfrf{t})}^8} + 2048\lambda^4\beta^{-4}d(d+2)(d+4)(d+6)(t-\lfrf{t})^4\nonumber\\
&\leq 128\lambda^8 \E\lsrs{\lbrb{h(\bar{\theta}_{\lfrf{t}}^\lambda)}^8} + 2048\lambda^4\beta^{-4}d(d+2)(d+4)(d+6)\nonumber\\
&\leq 128\lambda^8 K^8\E\lsrs{\left(1+|\bar{\theta}_{\lfrf{t}}^\lambda|^{r+1}\right)^8} + 2048\lambda^4\beta^{-4}d(d+2)(d+4)(d+6)\nonumber\\
&\leq 16384\lambda^8 K^8\left(1+\E\lsrs{|\bar{\theta}_{\lfrf{t}}^\lambda|^{8r+8}}\right) + 2048\lambda^4\beta^{-4}d(d+2)(d+4)(d+6)\nonumber\\
&\leq \lambda^4\left(e^{-\bar{a}\kappa_* n/2}\bar{C}_{1,1}\E\lsrs{|\bar{\theta}_0^{\lambda}|^{8r+8}}+\bar{C}_{2,1}\right),\nonumber
\end{align}
where
\begin{align}
\bar{C}_{1,1}&:=16384K^8,\nonumber\\
\bar{C}_{2,1}&:=16384K^8\left(1+c_*(4r+4)\left(1+\frac{1}{\bar{a}\kappa_*}\right)\right) + 2048\beta^{-4}d(d+2)(d+4)(d+6).\label{eqn:lemma4.8_constant_C11bar_C21bar}
\end{align}
Similarly, from Assumption \ref{assumption:2} and Corollary \ref{corollary:4.6}, one obtains
\begin{align}
\E\lsrs{\lbrb{\bar{\zeta}_t^{\lambda,n}-\bar{\zeta}_{\lfrf{t}}^{\lambda,n}}^4}
&=\E\lsrs{\lbrb{\int^t_{\lfrf{t}}-\lambda h(\bar{\zeta}_s^{\lambda, n})\ \dee s + \sqrt{2\lambda \beta^{-1}}(B_t^\lambda - B_{\lfrf{t}}^\lambda)}^4}\nonumber\\
&\leq 8\lambda^4\E\lsrs{\lbrb{\int^t_{\lfrf{t}}h(\bar{\zeta}_s^{\lambda, n})\ \dee s}^4} +  8\ \E\lsrs{\lbrb{\sqrt{2\lambda \beta^{-1}}(B_t^\lambda-B_{\lfrf{t}}^\lambda)}^4}\nonumber\nonumber\\
&\leq 8\lambda^4\E\lsrs{\left(\int^t_{\lfrf{t}}|h(\bar{\zeta}_s^{\lambda, n})|^4\ \dee s\right)\cdot\left(\int^t_{\lfrf{t}}1^{4/3}\ \dee s\right)^3} + 32\lambda^2\beta^{-2}d(d+2)\nonumber\\
&\leq 8\lambda^4\E\lsrs{\int^t_{\lfrf{t}}|h(\bar{\zeta}_s^{\lambda, n})|^4\ \dee s} + 32\lambda^2\beta^{-2}d(d+2)\nonumber\\
&\leq 8\lambda^4K^4\int^t_{\lfrf{t}}\E\lsrs{\left(1+|\bar{\zeta}_s^{\lambda,n}|^{r+1}\right)^4}\dee s + 32\lambda^2\beta^{-2}d(d+2)\nonumber\\
&\leq 64\lambda^4K^4\int^t_{\lfrf{t}}\E\lsrs{V_{4r+4}(\bar{\zeta}_s^{\lambda, n})}\ \dee s + 32\lambda^2\beta^{-2}d(d+2)\nonumber\\
&\leq \lambda^2\left(e^{- \bar{a}\min\{r+1, \kappa_*/2\}n}\bar{C}_{1,2}\E\lsrs{|\bar{\theta}_0^{\lambda}|^{4r+4}}+\bar{C}_{2,2}\right),\nonumber
\end{align}
where the second inequality is an application of H\"{o}lder's inequality, and
\begin{align}
\bar{C}_{1,2}:=&\ 2^{2r+7}K^4,\nonumber\\
\bar{C}_{2,2}:=&\ 64K^4\lsrs{2^{2r+1}+2^{2r+1}c_*(2r+2)\left(1+\frac{1}{\bar{a}\kappa_*}\right) + v_{4r+4}(M_V(4r+4))}  + 32\beta^{-2}d(d+2).\label{eqn:lemma4.8_constant_C12bar_C22bar}
\end{align}
This completes the proof. \qed

\begin{lemma}
\label{lemma:app211}
Let Assumptions \ref{assumption:1}, \ref{assumption:2}, \ref{assumption:3}, \ref{assumption:4} hold. Then, for any $\lambda \in (0,\tilde{\lambda}_{\max})$, $n\in\N_0$, and $t\in (nT, (n+1)T]$, the following inequality holds:
\begin{align} 
&-2\lambda\int^t_{nT}\E\lsrs{\lara{\bar{\zeta}_s^{\lambda,n}-\bar{\theta}_s^\lambda, h(\bar{\theta}_s^\lambda)-h(\bar{\theta}_{\lfrf{s}}^\lambda)-\nabla h(\bar{\theta}_{\lfrf{s}}^\lambda)(\bar{\theta}_s^\lambda-\bar{\theta}_{\lfrf{s}}^\lambda)}}\dee s \nonumber\\
 &\leq \lambda\bar{L}\int^t_{nT}\E\lsrs{\lbrb{\bar{\zeta}_s^{\lambda, n} - \bar{\theta}_s^\lambda}^2}\dee s\nonumber\\ 
&\quad +\lambda^2\Bigg{(}e^{-\bar{a}\kappa_* n/2}\cdot 27\bar{L}^{-1}\bar{L}_\nabla^2\cdot \E\lsrs{|\bar{\theta}_0^\lambda|^{4\nu}} + e^{-\bar{a}\kappa_* n/2}\cdot \frac{\bar{L}^{-1}\bar{L}_\nabla^2\bar{C}_{1,1}}{2}\cdot\E\lsrs{|\bar{\theta}_0^{\lambda}|^{8r+8}}\nonumber\\
&\quad +\bar{L}^{-1}\bar{L}_\nabla^2\left(\frac{27}{2} + 27c_*(2\nu)\left(1+\frac{1}{\bar{a}\kappa_*}\right) + \frac{\bar{C}_{2,1}}{2}\right)\Bigg{)},\nonumber
\end{align}
where the constant $\kappa_*$ is given in Corollary \ref{corollary:4.3} and the constants $\bar{C}_{1,1}, \bar{C}_{2,1}$ are given in (\ref{eqn:lemma4.8_constant_C11bar_C21bar}).
\end{lemma}

\textit{Proof.}  By Young's inequality, Remark \ref{remark:2.6}, and Lemma \ref{lemma:4.8}, one obtains
\begin{align}
&-2\lambda\int^t_{nT}\E\lsrs{\lara{\bar{\zeta}_s^{\lambda,n}-\bar{\theta}_s^\lambda, h(\bar{\theta}_s^\lambda)-h(\bar{\theta}_{\lfrf{s}}^\lambda)-\nabla h(\bar{\theta}_{\lfrf{s}}^\lambda)(\bar{\theta}_s^\lambda-\bar{\theta}_{\lfrf{s}}^\lambda)}}\dee s \nonumber\\
\leq &\ \lambda\bar{L}\int^t_{nT}\E\lsrs{\lbrb{\bar{\zeta}_s^{\lambda, n} - \bar{\theta}_s^\lambda}^2}\dee s + \lambda\bar{L}^{-1}\int^t_{nT}\E\lsrs{\lbrb{h(\bar{\theta}_s^\lambda)-h(\bar{\theta}_{\lfrf{s}}^\lambda)-\nabla h(\bar{\theta}_{\lfrf{s}}^\lambda)(\bar{\theta}_s^\lambda-\bar{\theta}_{\lfrf{s}}^\lambda)}^2} \dee s \nonumber\\
\leq &\ \lambda\bar{L}\int^t_{nT}\E\lsrs{\lbrb{\bar{\zeta}_s^{\lambda, n} - \bar{\theta}_s^\lambda}^2}\dee s + \lambda\bar{L}^{-1}L_\nabla^2\int^t_{nT}\E\lsrs{\left(1+|\bar{\theta}_s^\lambda|^\nu+|\bar{\theta}_{\lfrf{s}}^\lambda|^\nu\right)^2\cdot \lbrb{\bar{\theta}_s^\lambda-\bar{\theta}_{\lfrf{s}}^\lambda}^4} \dee s\nonumber\\
\leq &\ \lambda\bar{L}\int^t_{nT}\E\lsrs{\lbrb{\bar{\zeta}_s^{\lambda, n} - \bar{\theta}_s^\lambda}^2}\dee s + \lambda\bar{L}^{-1}\bar{L}_\nabla^2\int^t_{nT}\left\{\E\lsrs{\lbrb{1+|\bar{\theta}_s^\lambda|^\nu+|\bar{\theta}_{\lfrf{s}}^\lambda|^\nu}^4}\right\}^{1/2}\cdot\left\{\E\lsrs{\lbrb{\bar{\theta}_s^\lambda-\bar{\theta}_{\lfrf{s}}^\lambda}^8}\right\}^{1/2}\dee s\nonumber\\
\leq &\ \lambda\bar{L}\int^t_{nT}\E\lsrs{\lbrb{\bar{\zeta}_s^{\lambda, n} - \bar{\theta}_s^\lambda}^2}\dee s \nonumber\\
&\ + \lambda\bar{L}^{-1}\bar{L}_\nabla^2\int^t_{nT}\left\{27\E\lsrs{\left(1+|\bar{\theta}_s^\lambda|^{4\nu}+|\bar{\theta}_{\lfrf{s}}^\lambda|^{4\nu}\right)}\right\}^{1/2}\cdot\left\{\E\lsrs{\lbrb{\bar{\theta}_s^\lambda-\bar{\theta}_{\lfrf{s}}^\lambda}^8} \right\}^{1/2}\dee s\nonumber\\
\leq &\ \lambda\bar{L}\int^t_{nT}\E\lsrs{\lbrb{\bar{\zeta}_s^{\lambda, n} - \bar{\theta}_s^\lambda}^2}\dee s \nonumber\\
&\ + \lambda\bar{L}^{-1}\bar{L}_\nabla^2\int^t_{nT} \left\{27\left(1+2e^{-\bar{a}\kappa_* n/2}\E\lsrs{|\bar{\theta}_0^\lambda|^{4\nu}}+2c_*(2\nu)\left(1+\frac{1}{\bar{a}\kappa_*}\right)\right)\right\}^{1/2}\nonumber\\
&\ \times\left\{\lambda^4\left(e^{-\bar{a}\kappa_* n/2}\bar{C}_{1,1}\E\lsrs{|\bar{\theta}_0^{\lambda}|^{8r+8}}+\bar{C}_{2,1}\right) \right\}^{1/2}\dee s\nonumber\\
\leq &\ \lambda\bar{L}\int^t_{nT}\E\lsrs{\lbrb{\bar{\zeta}_s^{\lambda, n} - \bar{\theta}_s^\lambda}^2}\dee s \nonumber\\
&\ + \lambda^3\bar{L}^{-1}\bar{L}_\nabla^2\int^t_{nT} \left\{27\left(1+2e^{-\bar{a}\kappa_* n/2}\E\lsrs{|\bar{\theta}_0^\lambda|^{4\nu}}+2c_*(2\nu)\left(1+\frac{1}{\bar{a}\kappa_*}\right)\right)\right\}^{1/2}\nonumber\\
&\ \times\left\{e^{-\bar{a}\kappa_* n/2}\bar{C}_{1,1}\E\lsrs{|\bar{\theta}_0^{\lambda}|^{8r+8}}+\bar{C}_{2,1}\right\}^{1/2}\dee s\nonumber\\
\leq &\ \lambda\bar{L}\int^t_{nT}\E\lsrs{\lbrb{\bar{\zeta}_s^{\lambda, n} - \bar{\theta}_s^\lambda}^2}\dee s \nonumber\\
&\ + \lambda^2\bar{L}^{-1}\bar{L}_\nabla^2\left\{27\left(1+2e^{-\bar{a}\kappa_* n/2}\E\lsrs{|\bar{\theta}_0^\lambda|^{4\nu}}+2c_*(2\nu)\left(1+\frac{1}{\bar{a}\kappa_*}\right)\right)\right\}^{1/2}\nonumber\\
&\ \times\left\{e^{-\bar{a}\kappa_* n/2}\bar{C}_{1,1}\E\lsrs{|\bar{\theta}_0^{\lambda}|^{8r+8}}+\bar{C}_{2,1}\right\}^{1/2}\nonumber\\
\leq &\ \lambda\bar{L}\int^t_{nT}\E\lsrs{\lbrb{\bar{\zeta}_s^{\lambda, n} - \bar{\theta}_s^\lambda}^2}\dee s\nonumber\\ 
&\ +\lambda^2\Bigg{(}e^{-\bar{a}\kappa_* n/2}\cdot 27\bar{L}^{-1}\bar{L}_\nabla^2\cdot \E\lsrs{|\bar{\theta}_0^\lambda|^{4\nu}} + e^{-\bar{a}\kappa_* n/2}\cdot \frac{\bar{L}^{-1}\bar{L}_\nabla^2\bar{C}_{1,1}}{2}\cdot\E\lsrs{|\bar{\theta}_0^{\lambda}|^{8r+8}}\nonumber\\
&\ +\bar{L}^{-1}\bar{L}_\nabla^2\left(\frac{27}{2} + 27c_*(2\nu)\left(1+\frac{1}{\bar{a}\kappa_*}\right) + \frac{\bar{C}_{2,1}}{2}\right)\Bigg{)}\nonumber,
\end{align}
where the third inequality is an application of the Cauchy-Schwarz inequality.\qed

\begin{lemma}
\label{lemma:app2121}
Let Assumptions \ref{assumption:1}, \ref{assumption:2}, \ref{assumption:3}, \ref{assumption:4} hold. Then, for any $\lambda \in (0,\tilde{\lambda}_{\max})$, $n\in\N_0$, and $t\in (nT, (n+1)T]$, the following inequality holds:
\begin{align}
& -2\lambda\int^t_{nT}\E\lsrs{\lara{\bar{\zeta}_s^{\lambda,n}-\bar{\theta}_s^\lambda, \nabla h(\bar{\theta}_{\lfrf{s}}^\lambda)\left\{\int^s_{\lfrf{s}}-\lambda h_\lambda(\bar{\theta}_{\lfrf{u}}^\lambda)\ \dee u \right\}}}\dee s,\nonumber\\
\leq &\ \lambda\bar{L}\int^t_{nT}\E\lsrs{\lbrb{\bar{\zeta}_s^{\lambda, n} - \bar{\theta}_s^\lambda}^2}\dee s +\lambda^2\Bigg{(}e^{-\bar{a}\kappa_* n/2}\cdot 12\bar{L}^{-1}C_\nabla^2 K^2\cdot \E\lsrs{|\bar{\theta}_0^\lambda|^{2\nu+2r+4}}\nonumber\\
&\ +12\bar{L}^{-1}C_\nabla^2 K^2\left(1+c_*(\nu+r+2)\left(1+\frac{1}{\bar{a}\kappa_*}\right)\right)\Bigg{)},\nonumber
\end{align}
where the constants $\kappa_*$, $c_*(\nu+r+2)$ are given in Corollary \ref{corollary:4.3}.
\end{lemma}

\textit{Proof.} 
Applying in succession Young's inequality, Remark \ref{remark:2.6}, Assumption \ref{assumption:2}, and Corollary \ref{corollary:4.3} yields
\begin{align}
& -2\lambda\int^t_{nT}\E\lsrs{\lara{\bar{\zeta}_s^{\lambda,n}-\bar{\theta}_s^\lambda, \nabla h(\bar{\theta}_{\lfrf{s}}^\lambda)\left\{\int^s_{\lfrf{s}}-\lambda h_\lambda(\bar{\theta}_{\lfrf{u}}^\lambda)\ \dee u \right\}}}\dee s,\nonumber\\
\leq &\ \lambda\bar{L}\int^t_{nT}\E\lsrs{\lbrb{\bar{\zeta}_s^{\lambda, n} - \bar{\theta}_s^\lambda}^2}\dee s + \lambda\bar{L}^{-1}\int^t_{nT}\E\lsrs{\lbrb{\nabla h(\bar{\theta}_{\lfrf{s}}^\lambda)\left\{\int^s_{\lfrf{s}}-\lambda h_\lambda(\bar{\theta}_{\lfrf{u}}^\lambda)\ \dee u \right\}}^2}\dee s\nonumber\\
\leq &\ \lambda\bar{L}\int^t_{nT}\E\lsrs{\lbrb{\bar{\zeta}_s^{\lambda, n} - \bar{\theta}_s^\lambda}^2}\dee s + \lambda^3\bar{L}^{-1}\int^t_{nT}\E\lsrs{\lbrb{\nabla h(\bar{\theta}_{\lfrf{s}}^\lambda)h_\lambda(\bar{\theta}_{\lfrf{s}}^\lambda)}^2}\dee s\nonumber\\
\leq &\ \lambda\bar{L}\int^t_{nT}\E\lsrs{\lbrb{\bar{\zeta}_s^{\lambda, n} - \bar{\theta}_s^\lambda}^2}\dee s + \lambda^3\bar{L}^{-1}\int^t_{nT}\E\lsrs{|\nabla h(\bar{\theta}_{\lfrf{s}}^\lambda)|^2\cdot |h_\lambda(\bar{\theta}_{\lfrf{s}}^\lambda)|^2}\dee s\nonumber\\
\leq &\ \lambda\bar{L}\int^t_{nT}\E\lsrs{\lbrb{\bar{\zeta}_s^{\lambda, n} - \bar{\theta}_s^\lambda}^2}\dee s + \lambda^3\bar{L}^{-1}\int^t_{nT}\E\lsrs{|\nabla h(\bar{\theta}_{\lfrf{s}}^\lambda)|^2\cdot |h(\bar{\theta}_{\lfrf{s}}^\lambda)|^2}\dee s\nonumber\\
\leq &\ \lambda\bar{L}\int^t_{nT}\E\lsrs{\lbrb{\bar{\zeta}_s^{\lambda, n} - \bar{\theta}_s^\lambda}^2}\dee s + \lambda^3\bar{L}^{-1}C_\nabla^2 K^2\int^t_{nT}\E\lsrs{(1+|\bar{\theta}_{\lfrf{s}}^\lambda|^{\nu+1})^2\cdot (1+|\bar{\theta}_{\lfrf{s}}^\lambda|^{r+1})^2}\dee s\nonumber\\
\leq &\ \lambda\bar{L}\int^t_{nT}\E\lsrs{\lbrb{\bar{\zeta}_s^{\lambda, n} - \bar{\theta}_s^\lambda}^2}\dee s + 4\lambda^3\bar{L}^{-1}C_\nabla^2 K^2\int^t_{nT}\E\lsrs{(1+|\bar{\theta}_{\lfrf{s}}^\lambda|^{2\nu+2})(1+|\bar{\theta}_{\lfrf{s}}^\lambda|^{2r+2})}\dee s\nonumber\\
\leq &\ \lambda\bar{L}\int^t_{nT}\E\lsrs{\lbrb{\bar{\zeta}_s^{\lambda, n} - \bar{\theta}_s^\lambda}^2}\dee s + 12\lambda^3\bar{L}^{-1}C_\nabla^2 K^2\int^t_{nT}\left(1+\E\lsrs{|\bar{\theta}_{\lfrf{s}}^\lambda|^{2\nu+2r+4}}\right)\dee s\nonumber\\
\leq &\ \lambda\bar{L}\int^t_{nT}\E\lsrs{\lbrb{\bar{\zeta}_s^{\lambda, n} - \bar{\theta}_s^\lambda}^2}\dee s \nonumber\\
&\ + 12\lambda^3\bar{L}^{-1}C_\nabla^2 K^2\int^t_{nT}\lsrs{1+e^{-\bar{a}\kappa_* n/2}\E\lsrs{|\bar{\theta}_0^\lambda|^{2\nu+2r+4}}+c_*(\nu+r+2)\left(1+\frac{1}{\bar{a}\kappa_*}\right)}\dee s\nonumber\\
\leq &\ \lambda\bar{L}\int^t_{nT}\E\lsrs{\lbrb{\bar{\zeta}_s^{\lambda, n} - \bar{\theta}_s^\lambda}^2}\dee s +\lambda^2\Bigg{(}e^{-\bar{a}\kappa_* n/2}\cdot 12\bar{L}^{-1}C_\nabla^2 K^2\cdot \E\lsrs{|\bar{\theta}_0^\lambda|^{2\nu+2r+4}}\nonumber\\
&\ +12\bar{L}^{-1}C_\nabla^2 K^2\left(1+c_*(\nu+r+2)\left(1+\frac{1}{\bar{a}\kappa_*}\right)\right)\Bigg{)}.\nonumber
\end{align}\qed

\begin{lemma}
\label{lemma:app2122}
Let Assumptions \ref{assumption:1}, \ref{assumption:2}, \ref{assumption:3}, \ref{assumption:4} hold. Then, for any $\lambda \in (0,\tilde{\lambda}_{\max})$, $n\in\N_0$, and $t\in (nT, (n+1)T]$, the following inequality holds:
\begin{align}
&\left\{\E\lsrs{\lbrb{\int^s_{\lfrf{s}}(h(\bar{\zeta}_u^{\lambda,n})-h(\bar{\zeta}_{\lfrf{u}}^{\lambda, n}))\ \dee u}^2}\right\}^{1/2}\nonumber\\
\leq &\ \left\{\E\lsrs{\left(\int^s_{\lfrf{s}} 1^2\ \dee u\right)\cdot\left(\int^s_{\lfrf{s}}\lbrb{h(\bar{\zeta}_u^{\lambda,n})-h(\bar{\zeta}_{\lfrf{u}}^{\lambda, n})}^2\ \dee u\right)}\right\}^{1/2}\nonumber\\
\leq &\ \lambda^{1/2} \Bigg{(}e^{-\bar{a}\min\{r, \kappa_*/2\}n}\cdot L^2 6^r 2^{2r-2}\cdot\E\lsrs{|\bar{\theta}_0^{\lambda}|^{4r}} + e^{- \bar{a}\min\{r+1, \kappa_*/2\}n}\cdot\frac{L^26^r\bar{C}_{1,2}}{2}\cdot\E\lsrs{|\bar{\theta}_0^{\lambda}|^{4r+4}}\nonumber\\
&\ +  L^2 6^r\left(2^{2r-2}\left(1+c_*(2r)\left(1+\frac{1}{\bar{a}\kappa_*}\right)\right)+\frac{v_{4r}(M_V(4r))}{2}+\frac{\bar{C}_{2,2}}{2}\right)\Bigg{)}^{1/2},\label{eqn:lemma4.9_I2122_comp1_ub}\\
&\left\{\E\lsrs{\lbrb{\nabla h(\bar{\theta}_{\lfrf{s}}^\lambda)(B_s^\lambda - B_{\lfrf{s}}^\lambda)}^2}\right\}^{1/2}\nonumber\\
\leq&\ \left(e^{-\bar{a}\kappa_* n/2}\cdot\frac{C_\nabla^2}{2}\cdot \E\lsrs{|\bar{\theta}_0^\lambda|^{4\nu+4}}+\frac{C_\nabla^2}{2}\left(8d(d+2)+1+c_*(2\nu+2)\left(1+\frac{1}{\bar{a}\kappa_*}\right)\right)\right)^{1/2},\label{eqn:lemma4.9_I2122_comp2_ub}
\end{align}
where the constants $\kappa_*$, $c_*(2r)$, $c_*(2\nu+2)$ are given in Corollary \ref{corollary:4.3}, the constant $M_V(4r)$ is given in Lemma \ref{lemma:4.4}, and the constants $ \bar{C}_{1,2}, \bar{C}_{2,2}$ are given in (\ref{eqn:lemma4.8_constant_C12bar_C22bar}).
\end{lemma}

\textit{Proof.} We first show that \eqref{eqn:lemma4.9_I2122_comp1_ub} holds. Indeed, by the Cauchy-Schwarz inequality, Assumption \ref{assumption:2}, Corollary \ref{corollary:4.6}, and Lemma \ref{lemma:4.8}, one obtains
\begin{align}
&\left\{\E\lsrs{\lbrb{\int^s_{\lfrf{s}}(h(\bar{\zeta}_u^{\lambda,n})-h(\bar{\zeta}_{\lfrf{u}}^{\lambda, n}))\ \dee u}^2}\right\}^{1/2}\nonumber\\
\leq &\ \left\{\E\lsrs{\left(\int^s_{\lfrf{s}} 1^2\ \dee u\right)\cdot\left(\int^s_{\lfrf{s}}\lbrb{h(\bar{\zeta}_u^{\lambda,n})-h(\bar{\zeta}_{\lfrf{u}}^{\lambda, n})}^2\ \dee u\right)}\right\}^{1/2}\nonumber\\
\leq &\ \left\{\int^s_{\lfrf{s}}\E\lsrs{\lbrb{h(\bar{\zeta}_u^{\lambda,n})-h(\bar{\zeta}_{\lfrf{u}}^{\lambda, n})}^2}\ \dee u\right\}^{1/2}\nonumber\\
\leq &\ L\left\{\int^s_{\lfrf{s}}\E\lsrs{\left(1+|\bar{\zeta}_u^{\lambda,n}|+|\bar{\zeta}_{\lfrf{u}}^{\lambda, n}|\right)^{2r}\cdot \lbrb{\bar{\zeta}_u^{\lambda,n}-\bar{\zeta}_{\lfrf{u}}^{\lambda, n}}^2}\ \dee u\right\}^{1/2}\nonumber\\
\leq &\ L\left\{\int^s_{\lfrf{s}}\left(\E\lsrs{\left(1+|\bar{\zeta}_u^{\lambda,n}|+|\bar{\zeta}_{\lfrf{u}}^{\lambda, n}|\right)^{4r}}\right)^{1/2}\cdot \left(\E\lsrs{\lbrb{\bar{\zeta}_u^{\lambda,n}-\bar{\zeta}_{\lfrf{u}}^{\lambda, n}}^4}\right)^{1/2}\ \dee u\right\}^{1/2}\nonumber\\
\leq &\ L\left\{\int^s_{\lfrf{s}}\left(3^{2r}\E\lsrs{\left(1+|\bar{\zeta}_u^{\lambda,n}|^2+|\bar{\zeta}_{\lfrf{u}}^{\lambda, n}|^2\right)^{2r}}\right)^{1/2}\cdot \left(\E\lsrs{\lbrb{\bar{\zeta}_u^{\lambda,n}-\bar{\zeta}_{\lfrf{u}}^{\lambda, n}}^4}\right)^{1/2}\ \dee u\right\}^{1/2}\nonumber\\
\leq &\ L\left\{\int^s_{\lfrf{s}}\left(3^{2r}\cdot 2^{2r-1}\E\lsrs{V_{4r}(\bar{\zeta}_u^{\lambda,n})+V_{4r}(\bar{\zeta}_{\lfrf{u}}^{\lambda,n})}\right)^{1/2}\cdot \left(\E\lsrs{\lbrb{\bar{\zeta}_u^{\lambda,n}-\bar{\zeta}_{\lfrf{u}}^{\lambda, n}}^4}\right)^{1/2}\ \dee u\right\}^{1/2}\nonumber\\
\leq &\ L\ \Bigg{\{}\int^s_{\lfrf{s}}6^r\times\left(e^{-\bar{a}\min\{r, \kappa_*/2\}n}\cdot 2^{2r-1} \cdot\E\lsrs{|\bar{\theta}_0^{\lambda}|^{4r}} + 2^{2r-1}\left(1+c_*(2r)\left(1+\frac{1}{\bar{a}\kappa_*}\right)\right) + v_{4r}(M_V(4r))\right)^{1/2}\nonumber\\
&\ \times \left(\lambda^2\left(e^{- \bar{a}\min\{r+1, \kappa_*/2\}n}\bar{C}_{1,2}\E\lsrs{|\bar{\theta}_0^{\lambda}|^{4r+4}}+\bar{C}_{2,2}\right)\right)^{1/2}\ \dee u\Bigg{\}}^{1/2}\nonumber\\
\leq &\ \lambda^{1/2} \Bigg{(}e^{-\bar{a}\min\{r, \kappa_*/2\}n}\cdot L^2 6^r 2^{2r-2}\cdot\E\lsrs{|\bar{\theta}_0^{\lambda}|^{4r}} + e^{- \bar{a}\min\{r+1, \kappa_*/2\}n}\cdot\frac{L^26^r\bar{C}_{1,2}}{2}\cdot\E\lsrs{|\bar{\theta}_0^{\lambda}|^{4r+4}}\nonumber\\
&\ +  L^2 6^r\left(2^{2r-2}\left(1+c_*(2r)\left(1+\frac{1}{\bar{a}\kappa_*}\right)\right)+\frac{v_{4r}(M_V(4r))}{2}+\frac{\bar{C}_{2,2}}{2}\right)\Bigg{)}^{1/2}.\nonumber
\end{align}
Moreover, we note that (\ref{eqn:lemma4.9_I2122_comp1_ub}) still holds even in the case $r=0$ in view of Remark \ref{remark:4.7}.\\

To show that \eqref{eqn:lemma4.9_I2122_comp2_ub} holds, we note that
\begin{align}
&\left\{\E\lsrs{\lbrb{\nabla h(\bar{\theta}_{\lfrf{s}}^\lambda)(B_s^\lambda - B_{\lfrf{s}}^\lambda)}^2}\right\}^{1/2}\nonumber\\
\leq&\ \left\{\E\lsrs{|\nabla h(\bar{\theta}_{\lfrf{s}}^\lambda)|^2\cdot |B_s^\lambda - B_{\lfrf{s}}^\lambda|^2}\right\}^{1/2}\nonumber\\
\leq&\ \left\{\E\lsrs{\lbrb{\nabla h(\bar{\theta}_{\lfrf{s}}^\lambda)}^4}\right\}^{1/4}\cdot\left\{\E\lsrs{\lbrb{B_s^\lambda-B_{\lfrf{s}}^\lambda}^4}\right\}^{1/4}\nonumber\\
\leq&\ C_\nabla \cdot (d(d+2))^{1/4} \left\{\E\lsrs{\left(1+|\bar{\theta}_{\lfrf{s}}^\lambda|^{\nu + 1}\right)^4}\right\}^{1/4}\nonumber\\
\leq&\ C_\nabla \cdot (8d(d+2))^{1/4} \left(1+\E\lsrs{|\bar{\theta}_{\lfrf{s}}^\lambda|^{4\nu + 4}}\right)^{1/4}\nonumber\\
\leq&\ C_\nabla \cdot (8d(d+2))^{1/4} \left(1+e^{-\bar{a}\kappa_* n/2}\E\lsrs{|\bar{\theta}_0^\lambda|^{4\nu+4}}+c_*(2\nu+2)\left(1+\frac{1}{\bar{a}\kappa_*}\right)\right)^{1/4}\nonumber\\
\leq&\ \left(e^{-\bar{a}\kappa_* n/2}\cdot\frac{C_\nabla^2}{2}\cdot \E\lsrs{|\bar{\theta}_0^\lambda|^{4\nu+4}}+\frac{C_\nabla^2}{2}\left(8d(d+2)+1+c_*(2\nu+2)\left(1+\frac{1}{\bar{a}\kappa_*}\right)\right)\right)^{1/2}.\nonumber
\end{align}\qed

\begin{lemma}
\label{lemma:app31}
Let Assumptions \ref{assumption:1}, \ref{assumption:2}, \ref{assumption:3}, \ref{assumption:4} hold. Then, for any $\lambda \in (0,\tilde{\lambda}_{\max})$, $n\in\N_0$, and $t\in (nT, (n+1)T]$, the following inequality holds:
\begin{align}
& -2\lambda\int^t_{nT}\E\lsrs{\lara{\bar{\zeta}_s^{\lambda,n}-\bar{\theta}_s^\lambda, h(\bar{\theta}_{\lfrf{s}}^\lambda)-h_\lambda(\bar{\theta}_{\lfrf{s}}^\lambda)}}\dee s\nonumber\\
\leq&\ \lambda\bar{L}\int^t_{nT}\E\lsrs{\lbrb{\bar{\zeta}_s^{\lambda, n} - \bar{\theta}_s^\lambda}^2}\dee s + \lambda^2\Bigg{(}e^{-\bar{a}\kappa_* n/2}\cdot 4\bar{L}^{-1}K^2\cdot \E\lsrs{|\bar{\theta}_0^\lambda|^{6r+2}}\nonumber\\
&\ +2\bar{L}^{-1}K^2\left(1+2c_*(3r+1)\left(1+\frac{1}{\bar{a}\kappa_*}\right)\right)\Bigg{)},\nonumber
\end{align}
where the constants $\kappa_*$, $c_*(3r+1)$ are given in Corollary \ref{corollary:4.3}.
\end{lemma}

\textit{Proof.} By using Young's inequality, Assumption \ref{assumption:2}, and Corollary \ref{corollary:4.3}, we obtain that 
\begin{align}
& -2\lambda\int^t_{nT}\E\lsrs{\lara{\bar{\zeta}_s^{\lambda,n}-\bar{\theta}_s^\lambda, h(\bar{\theta}_{\lfrf{s}}^\lambda)-h_\lambda(\bar{\theta}_{\lfrf{s}}^\lambda)}}\dee s\nonumber\\
\leq&\ \lambda\bar{L}\int^t_{nT}\E\lsrs{\lbrb{\bar{\zeta}_s^{\lambda, n} - \bar{\theta}_s^\lambda}^2}\dee s + \lambda\bar{L}^{-1}\int^t_{nT}\E\lsrs{\lbrb{h(\bar{\theta}_{\lfrf{s}}^\lambda)-h_\lambda(\bar{\theta}_{\lfrf{s}}^\lambda)}^2}\dee s\nonumber\\
=&\ \lambda\bar{L}\int^t_{nT}\E\lsrs{\lbrb{\bar{\zeta}_s^{\lambda, n} - \bar{\theta}_s^\lambda}^2}\dee s + \lambda\bar{L}^{-1}\int^t_{nT}\E\lsrs{\lbrb{\frac{h(\bar{\theta}_{\lfrf{s}}^\lambda)\left((1+\lambda|\bar{\theta}_{\lfrf{s}}^\lambda|^{2r})^{1/2}-1\right)}{(1+\lambda|\bar{\theta}_{\lfrf{s}}^\lambda|^{2r})^{1/2}}}^2}\dee s\nonumber\\
\leq&\ \lambda\bar{L}\int^t_{nT}\E\lsrs{\lbrb{\bar{\zeta}_s^{\lambda, n} - \bar{\theta}_s^\lambda}^2}\dee s + \lambda\bar{L}^{-1}\int^t_{nT}\E\lsrs{\frac{|h(\bar{\theta}_{\lfrf{s}}^\lambda)|^2\left(1+\lambda|\bar{\theta}_{\lfrf{s}}^\lambda|^{2r}-1\right)^2}{1+\lambda|\bar{\theta}_{\lfrf{s}}^\lambda|^{2r}}}\dee s\nonumber\\
\leq&\ \lambda\bar{L}\int^t_{nT}\E\lsrs{\lbrb{\bar{\zeta}_s^{\lambda, n} - \bar{\theta}_s^\lambda}^2}\dee s + \lambda^3\bar{L}^{-1}\int^t_{nT}\E\lsrs{|h(\bar{\theta}_{\lfrf{s}}^\lambda)|^2\cdot |\bar{\theta}_{\lfrf{s}}^\lambda|^{4r}}\dee s\nonumber\\
\leq&\ \lambda\bar{L}\int^t_{nT}\E\lsrs{\lbrb{\bar{\zeta}_s^{\lambda, n} - \bar{\theta}_s^\lambda}^2}\dee s + \lambda^3\bar{L}^{-1}K^2\int^t_{nT}\E\lsrs{(1+|\bar{\theta}_{\lfrf{s}}^\lambda|^{r+1})^2\cdot |\bar{\theta}_{\lfrf{s}}^\lambda|^{4r}}\dee s\nonumber\\
\leq&\ \lambda\bar{L}\int^t_{nT}\E\lsrs{\lbrb{\bar{\zeta}_s^{\lambda, n} - \bar{\theta}_s^\lambda}^2}\dee s + 2\lambda^3\bar{L}^{-1}K^2\int^t_{nT}\E\lsrs{(1+|\bar{\theta}_{\lfrf{s}}^\lambda|^{2r+2})\cdot |\bar{\theta}_{\lfrf{s}}^\lambda|^{4r}}\dee s\nonumber\\
\leq&\ \lambda\bar{L}\int^t_{nT}\E\lsrs{\lbrb{\bar{\zeta}_s^{\lambda, n} - \bar{\theta}_s^\lambda}^2}\dee s + 2\lambda^3\bar{L}^{-1}K^2\int^t_{nT}\left(1+2\E\lsrs{|\bar{\theta}_{\lfrf{s}}^\lambda|^{6r+2}}\right)\dee s\nonumber\\
\leq&\ \lambda\bar{L}\int^t_{nT}\E\lsrs{\lbrb{\bar{\zeta}_s^{\lambda, n} - \bar{\theta}_s^\lambda}^2}\dee s\nonumber\\
&\ + 2\lambda^3\bar{L}^{-1}K^2\int^t_{nT}\left(1+2e^{-\bar{a}\kappa_*n/2}\E\lsrs{|\bar{\theta}_0^\lambda|^{6r+2}}+2c_*(3r+1)\left(1+\frac{1}{\bar{a}\kappa_*}\right)\right)\dee s\nonumber\\
\leq&\ \lambda\bar{L}\int^t_{nT}\E\lsrs{\lbrb{\bar{\zeta}_s^{\lambda, n} - \bar{\theta}_s^\lambda}^2}\dee s + \lambda^2\Bigg{(}e^{-\bar{a}\kappa_* n/2}\cdot 4\bar{L}^{-1}K^2\cdot \E\lsrs{|\bar{\theta}_0^\lambda|^{6r+2}}\nonumber\\
&\ +2\bar{L}^{-1}K^2\left(1+2c_*(3r+1)\left(1+\frac{1}{\bar{a}\kappa_*}\right)\right)\Bigg{)}.\nonumber
\end{align}\qed

\subsection*{\textit{Proof of Lemma \ref{lemma:4.9}}}
\label{proof:lemma4.9}
By Ito's formula, one obtains, for every $n\in\N_0$, $t\in (nT, (n+1)T]$, 
\begin{align}
W_2^2(\mc{L}(\bar{\theta}_t^\lambda), \mc{L}(\bar{\zeta}_t^{\lambda,n}))
&\leq \E\lsrs{\lbrb{\bar{\zeta}_t^{\lambda, n} - \bar{\theta}_t^\lambda}^2}\nonumber\\
&=-2\lambda\int^t_{nT}\E\lsrs{\lara{\bar{\zeta}_s^{\lambda,n}-\bar{\theta}_s^\lambda, h(\bar{\zeta}_s^{\lambda,n})-h_\lambda(\bar{\theta}_{\lfrf{s}}^\lambda)}}\dee s\nonumber\\
&\leq \bar{I}_{1,1} + \bar{I}_{2,1} + \bar{I}_{3,1},\label{eqn:lemma4.9_wass2_ub}
\end{align}
where
\begin{align}
\bar{I}_{1,1} &:= -2\lambda\int^t_{nT}\E\lsrs{\lara{\bar{\zeta}_s^{\lambda,n}-\bar{\theta}_s^\lambda, h(\bar{\zeta}_s^{\lambda,n})-h(\bar{\theta}_s^\lambda)}}\dee s,\nonumber\\
\bar{I}_{2,1} &:= -2\lambda\int^t_{nT}\E\lsrs{\lara{\bar{\zeta}_s^{\lambda,n}-\bar{\theta}_s^\lambda, h(\bar{\theta}_s^\lambda)-h(\bar{\theta}_{\lfrf{s}}^\lambda)}}\dee s,\label{eqn:lemma4.9_auxdef_I11_I21_I31}\\
\bar{I}_{3,1} &:= -2\lambda\int^t_{nT}\E\lsrs{\lara{\bar{\zeta}_s^{\lambda,n}-\bar{\theta}_s^\lambda, h(\bar{\theta}_{\lfrf{s}}^\lambda)-h_\lambda(\bar{\theta}_{\lfrf{s}}^\lambda)}}\dee s.\nonumber
\end{align}
We now bound each of these integrals. One obtains immediately from Remark \ref{remark:2.4} that
\begin{align}
\bar{I}_{1,1} \leq 2\lambda\bar{L}\int^t_{nT}\E\lsrs{\lbrb{\bar{\zeta}_s^{\lambda, n} - \bar{\theta}_s^\lambda}^2}\dee s.\label{eqn:lemma4.9_I11_ub}
\end{align}
The second integral $\bar{I}_{2,1}$ of (\ref{eqn:lemma4.9_wass2_ub}) can be further expanded as $\bar{I}_{2,1}=\bar{I}_{2,1}^{(1)} + \bar{I}_{2,1}^{(2)}$, where
\begin{align}
\bar{I}_{2,1}^{(1)}&:= -2\lambda\int^t_{nT}\E\lsrs{\lara{\bar{\zeta}_s^{\lambda,n}-\bar{\theta}_s^\lambda, h(\bar{\theta}_s^\lambda)-h(\bar{\theta}_{\lfrf{s}}^\lambda)-\nabla h(\bar{\theta}_{\lfrf{s}}^\lambda)(\bar{\theta}_s^\lambda-\bar{\theta}_{\lfrf{s}}^\lambda)}}\dee s,\nonumber\\
\bar{I}_{2,1}^{(2)}&:= -2\lambda\int^t_{nT}\E\lsrs{\lara{\bar{\zeta}_s^{\lambda,n}-\bar{\theta}_s^\lambda, \nabla h(\bar{\theta}_{\lfrf{s}}^\lambda)(\bar{\theta}_s^\lambda-\bar{\theta}_{\lfrf{s}}^\lambda)}}\dee s.\label{eqn:lemma4.9_auxdef_I211_I212}
\end{align}

By using Lemma \ref{lemma:app211}, we have that
\begin{align} 
\bar{I}_{2,1}^{(1)} &\leq \lambda\bar{L}\int^t_{nT}\E\lsrs{\lbrb{\bar{\zeta}_s^{\lambda, n} - \bar{\theta}_s^\lambda}^2}\dee s\nonumber\\ 
&\quad +\lambda^2\Bigg{(}e^{-\bar{a}\kappa_* n/2}\cdot 27\bar{L}^{-1}\bar{L}_\nabla^2\cdot \E\lsrs{|\bar{\theta}_0^\lambda|^{4\nu}} + e^{-\bar{a}\kappa_* n/2}\cdot \frac{\bar{L}^{-1}\bar{L}_\nabla^2\bar{C}_{1,1}}{2}\cdot\E\lsrs{|\bar{\theta}_0^{\lambda}|^{8r+8}}\nonumber\\
&\quad +\bar{L}^{-1}\bar{L}_\nabla^2\left(\frac{27}{2} + 27c_*(2\nu)\left(1+\frac{1}{\bar{a}\kappa_*}\right) + \frac{\bar{C}_{2,1}}{2}\right)\Bigg{)}.\label{eqn:lemma4.9_I211_ub}
\end{align}

The integral $\bar{I}_{2,1}^{(2)}$ can be further expanded as
\begin{align}
\bar{I}_{2,1}^{(2)}
&=-2\lambda\int^t_{nT}\E\lsrs{\lara{\bar{\zeta}_s^{\lambda,n}-\bar{\theta}_s^\lambda, \nabla h(\bar{\theta}_{\lfrf{s}}^\lambda)\left\{\int^s_{\lfrf{s}}-\lambda h_\lambda(\bar{\theta}_{\lfrf{u}}^\lambda)\ \dee u + \sqrt{2\lambda\beta^{-1}}(B_s^\lambda - B_{\lfrf{s}}^\lambda)\right\}}}\dee s\nonumber\\
&=\bar{I}_{2,1}^{(2,1)}+\bar{I}_{2,1}^{(2,2)},\label{eqn:lemma4.9_I212_eq}
\end{align}
where
\begin{align}
\bar{I}_{2,1}^{(2,1)}&:= -2\lambda\int^t_{nT}\E\lsrs{\lara{\bar{\zeta}_s^{\lambda,n}-\bar{\theta}_s^\lambda, \nabla h(\bar{\theta}_{\lfrf{s}}^\lambda)\left\{\int^s_{\lfrf{s}}-\lambda h_\lambda(\bar{\theta}_{\lfrf{u}}^\lambda)\ \dee u \right\}}}\dee s,\nonumber\\
\bar{I}_{2,1}^{(2,2)}&:=-2\lambda\int^t_{nT}\E\lsrs{\lara{\bar{\zeta}_s^{\lambda,n}-\bar{\theta}_s^\lambda, \nabla h(\bar{\theta}_{\lfrf{s}}^\lambda)\sqrt{2\lambda\beta^{-1}}(B_s^\lambda - B_{\lfrf{s}}^\lambda)}}\dee s.\label{eqn:lemma4.9_auxdef_I2121_I2122}
\end{align}
By using Lemma \ref{lemma:app2121}, we obtain that
\begin{align}
\bar{I}_{2,1}^{(2,1)}
\leq &\ \lambda\bar{L}\int^t_{nT}\E\lsrs{\lbrb{\bar{\zeta}_s^{\lambda, n} - \bar{\theta}_s^\lambda}^2}\dee s +\lambda^2\Bigg{(}e^{-\bar{a}\kappa_* n/2}\cdot 12\bar{L}^{-1}C_\nabla^2 K^2\cdot \E\lsrs{|\bar{\theta}_0^\lambda|^{2\nu+2r+4}}\nonumber\\
&\ +12\bar{L}^{-1}C_\nabla^2 K^2\left(1+c_*(\nu+r+2)\left(1+\frac{1}{\bar{a}\kappa_*}\right)\right)\Bigg{)}.\label{eqn:lemma4.9_I2121_ub}
\end{align}

It remains to bound $\bar{I}_{2,1}^{(2,2)}$.  Observe that
\begin{align}
&\bar{I}_{2,1}^{(2,2)}\nonumber\\
=&\ -2\lambda\int^t_{nT}\E\lsrs{\lara{\left\{\bar{\zeta}_{\lfrf{s}}^{\lambda,n}-\bar{\theta}_{\lfrf{s}}^\lambda-\lambda\int^s_{\lfrf{s}}(h(\bar{\zeta}_u^{\lambda,n})-h_\lambda(\bar{\theta}_{\lfrf{u}}^\lambda))\ \dee u\right\}, \nabla h(\bar{\theta}_{\lfrf{s}}^\lambda)\sqrt{2\lambda\beta^{-1}}(B_s^\lambda - B_{\lfrf{s}}^\lambda)}}\dee s\nonumber\\
=&\ -2\lambda\int^t_{nT}\E\lsrs{\lara{\bar{\zeta}_{\lfrf{s}}^{\lambda,n}-\bar{\theta}_{\lfrf{s}}^\lambda, \nabla h(\bar{\theta}_{\lfrf{s}}^\lambda)\sqrt{2\lambda\beta^{-1}}(B_s^\lambda - B_{\lfrf{s}}^\lambda)}}\dee s\nonumber\\
&+ 2\lambda^2\int^t_{nT}\E\lsrs{\lara{\int^s_{\lfrf{s}}(h(\bar{\zeta}_u^{\lambda,n})-h_\lambda(\bar{\theta}_{\lfrf{u}}^\lambda))\ \dee u, \nabla h(\bar{\theta}_{\lfrf{s}}^\lambda)\sqrt{2\lambda\beta^{-1}}(B_s^\lambda - B_{\lfrf{s}}^\lambda)}}\dee s.\label{eqn:lemma4.9_I2122_eq}
\end{align}
Note that for each $s\in (nT, t]$, both $\bar{\zeta}_{\lfrf{s}}^{\lambda,n}$ and $\bar{\theta}_{\lfrf{s}}^\lambda$ are $\mc{F}_{\lfrf{s}}$-measurable and $(B_s^\lambda - B_{\lfrf{s}}^\lambda)$ is independent of $\mc{F}_{\lfrf{s}}$. Since $\E\lsrs{(B_s^\lambda - B_{\lfrf{s}}^\lambda)}=0$, one obtains as a result of conditioning on $\mc{F}_{\lfrf{s}}$, the Tower Property of conditional expectation, and the fact that $\lfrf{u}=\lfrf{s}$ for $u\in [\lfrf{s}, s]$ that
\begin{align}
\E\lsrs{\lara{\bar{\zeta}_{\lfrf{s}}^{\lambda,n}-\bar{\theta}_{\lfrf{s}}^\lambda, \nabla h(\bar{\theta}_{\lfrf{s}}^\lambda)\sqrt{2\lambda\beta^{-1}}(B_s^\lambda - B_{\lfrf{s}}^\lambda)}} &= 0,\label{eqn:lemma4.9_I2122_ip1}\\
\E\lsrs{\lara{\int^s_{\lfrf{s}}h_\lambda(\bar{\theta}_{\lfrf{u}}^\lambda)\ \dee u, \nabla h(\bar{\theta}_{\lfrf{s}}^\lambda)\sqrt{2\lambda\beta^{-1}}(B_s^\lambda - B_{\lfrf{s}}^\lambda)}} &= 0,\label{eqn:lemma4.9_I2122_ip2}\\
\E\lsrs{\lara{\int^s_{\lfrf{s}}h(\bar{\zeta}_{\lfrf{u}}^{\lambda,n})\ \dee u, \nabla h(\bar{\theta}_{\lfrf{s}}^\lambda)\sqrt{2\lambda\beta^{-1}}(B_s^\lambda - B_{\lfrf{s}}^\lambda)}} &= 0.\label{eqn:lemma4.9_I2122_ip3}
\end{align} 
Therefore, substituting (\ref{eqn:lemma4.9_I2122_ip1}), (\ref{eqn:lemma4.9_I2122_ip2}), and (\ref{eqn:lemma4.9_I2122_ip3}) into (\ref{eqn:lemma4.9_I2122_eq}) yields
\begin{align}
\bar{I}_{2,1}^{(2,2)}
=&\ 2\lambda^{5/2}\sqrt{2\beta^{-1}}\int^t_{nT}\E\lsrs{\lara{\int^s_{\lfrf{s}}(h(\bar{\zeta}_u^{\lambda,n})-h(\bar{\zeta}_{\lfrf{u}}^{\lambda, n}))\ \dee u, \nabla h(\bar{\theta}_{\lfrf{s}}^\lambda)(B_s^\lambda - B_{\lfrf{s}}^\lambda)}}\dee s\nonumber\\
\leq& 2\lambda^{5/2}\sqrt{2\beta^{-1}}\int^t_{nT}\E\lsrs{\lbrb{\int^s_{\lfrf{s}}(h(\bar{\zeta}_u^{\lambda,n})-h(\bar{\zeta}_{\lfrf{u}}^{\lambda, n}))\ \dee u}\cdot \lbrb{\nabla h(\bar{\theta}_{\lfrf{s}}^\lambda)(B_s^\lambda - B_{\lfrf{s}}^\lambda)}}\dee s\nonumber\\
\leq& 2\lambda^{5/2}\sqrt{2\beta^{-1}}\int^t_{nT}\left\{\E\lsrs{\lbrb{\int^s_{\lfrf{s}}(h(\bar{\zeta}_u^{\lambda,n})-h(\bar{\zeta}_{\lfrf{u}}^{\lambda, n}))\ \dee u}^2}\right\}^{1/2}\cdot \left\{\E\lsrs{\lbrb{\nabla h(\bar{\theta}_{\lfrf{s}}^\lambda)(B_s^\lambda - B_{\lfrf{s}}^\lambda)}^2}\right\}^{1/2}\dee s\label{eqn:lemma4.9_I2122_ub},
\end{align}
where the Cauchy-Schwarz inequality was applied twice successively. Using Lemma \ref{lemma:app2122} and substituting (\ref{eqn:lemma4.9_I2122_comp1_ub}) and (\ref{eqn:lemma4.9_I2122_comp2_ub}) into (\ref{eqn:lemma4.9_I2122_ub}) then yield
\begin{align}
\bar{I}_{2,1}^{(2,2)}
\leq&\ 2\lambda^3\sqrt{2\beta^{-1}}\int^t_{nT}\Bigg{(}e^{-\bar{a}\min\{r, \kappa_*/2\}n}\cdot L^2 6^r 2^{2r-2}\cdot\E\lsrs{|\bar{\theta}_0^{\lambda}|^{4r}}\nonumber\\
&\ + e^{- \bar{a}\min\{r+1, \kappa_*/2\}n}\cdot\frac{L^26^r\bar{C}_{1,2}}{2}\cdot\E\lsrs{|\bar{\theta}_0^{\lambda}|^{4r+4}}\nonumber\\
&\ +  L^2 6^r\left(2^{2r-2}\left(1+c_*(2r)\left(1+\frac{1}{\bar{a}\kappa_*}\right)\right)+\frac{v_{4r}(M_V(4r))}{2}+\frac{\bar{C}_{2,2}}{2}\right)\Bigg{)}^{1/2}\nonumber\\
&\ \times \Bigg{(}e^{-\bar{a}\kappa_* n/2}\cdot\frac{C_\nabla^2}{2}\cdot \E\lsrs{|\bar{\theta}_0^\lambda|^{4\nu+4}}  +\frac{C_\nabla^2}{2}\left(8d(d+2)+1+c_*(2\nu+2)\left(1+\frac{1}{\bar{a}\kappa_*}\right)\right)\Bigg{)}^{1/2}\ \dee s\nonumber\\
\leq &\ \lambda^2\Bigg{(}e^{-\bar{a}\min\{r, \kappa_*/2\}n}\cdot \sqrt{2\beta^{-1}}L^2 6^r 2^{2r-2}\cdot\E\lsrs{|\bar{\theta}_0^{\lambda}|^{4r}}\nonumber\\
&\ + e^{- \bar{a}\min\{r+1, \kappa_*/2\}n}\cdot\frac{\sqrt{2\beta^{-1}}L^26^r\bar{C}_{1,2}}{2}\cdot\E\lsrs{|\bar{\theta}_0^{\lambda}|^{4r+4}}\nonumber\\
&\ +e^{-\bar{a}\kappa_* n/2}\cdot\frac{\sqrt{2\beta^{-1}}C_\nabla^2}{2}\cdot \E\lsrs{|\bar{\theta}_0^\lambda|^{4\nu+4}}\nonumber\\
&\ + \sqrt{2\beta^{-1}}L^2 6^r\left(2^{2r-2}\left(1+c_*(2r)\left(1+\frac{1}{\bar{a}\kappa_*}\right)\right)+\frac{v_{4r}(M_V(4r))}{2}+\frac{\bar{C}_{2,2}}{2}\right)\nonumber\\
&\ + \frac{\sqrt{2\beta^{-1}}C_\nabla^2}{2}\left(8d(d+2)+1+c_*(2\nu+2)\left(1+\frac{1}{\bar{a}\kappa_*}\right)\right)\Bigg{).}\label{eqn:lemma4.9_I2122_ub2}
\end{align}
Furthermore, by using Lemma \ref{lemma:app31}, we have that
\begin{align}
\bar{I}_{3,1}\leq&\ \lambda\bar{L}\int^t_{nT}\E\lsrs{\lbrb{\bar{\zeta}_s^{\lambda, n} - \bar{\theta}_s^\lambda}^2}\dee s + \lambda^2\Bigg{(}e^{-\bar{a}\kappa_* n/2}\cdot 4\bar{L}^{-1}K^2\cdot \E\lsrs{|\bar{\theta}_0^\lambda|^{6r+2}}\nonumber\\
&\ +2\bar{L}^{-1}K^2\left(1+2c_*(3r+1)\left(1+\frac{1}{\bar{a}\kappa_*}\right)\right)\Bigg{)}.\label{eqn:lemma4.9_I31_ub}
\end{align}
Combining the upper bounds in (\ref{eqn:lemma4.9_I11_ub}), (\ref{eqn:lemma4.9_I211_ub}), (\ref{eqn:lemma4.9_I2121_ub}), (\ref{eqn:lemma4.9_I2122_ub2}), and (\ref{eqn:lemma4.9_I31_ub}), we obtain
\begin{align}
&\E\lsrs{\lbrb{\bar{\zeta}_t^{\lambda, n} - \bar{\theta}_t^\lambda}^2}\nonumber\\
\leq &\ \bar{I}_{1,1} + \bar{I}_{2,1}^{(1)} + \bar{I}_{2,1}^{(2,1)}+ \bar{I}_{2,1}^{(2,2)} +\bar{I}_{3,1}\nonumber\\
\leq&\ 5\lambda\bar{L}\int^t_{nT}\E\lsrs{\lbrb{\bar{\zeta}_s^{\lambda, n} - \bar{\theta}_s^\lambda}^2}\dee s\nonumber\\
&\ +\lambda^2\Bigg{(}e^{-\bar{a}\kappa_* n/2}\cdot 27\bar{L}^{-1}\bar{L}_\nabla^2\cdot \E\lsrs{|\bar{\theta}_0^\lambda|^{4\nu}} + e^{-\bar{a}\kappa_* n/2}\cdot \frac{\bar{L}^{-1}\bar{L}_\nabla^2\bar{C}_{1,1}}{2}\cdot\E\lsrs{|\bar{\theta}_0^{\lambda}|^{8r+8}}\nonumber\\
&\ \qquad +\bar{L}^{-1}\bar{L}_\nabla^2\left(\frac{27}{2} + 27c_*(2\nu)\left(1+\frac{1}{\bar{a}\kappa_*}\right) + \frac{\bar{C}_{2,1}}{2}\right)\Bigg{)}\nonumber\\
&\ +\lambda^2\Bigg{(}e^{-\bar{a}\kappa_* n/2}\cdot 12\bar{L}^{-1}C_\nabla^2 K^2\cdot \E\lsrs{|\bar{\theta}_0^\lambda|^{2\nu+2r+4}} +12\bar{L}^{-1}C_\nabla^2 K^2\left(1+c_*(\nu+r+2)\left(1+\frac{1}{\bar{a}\kappa_*}\right)\right)\Bigg{)}\nonumber\\
&\ +\lambda^2\Bigg{(}e^{-\bar{a}\min\{r, \kappa_*/2\}n}\cdot \sqrt{2\beta^{-1}}L^2 6^r 2^{2r-2}\cdot\E\lsrs{|\bar{\theta}_0^{\lambda}|^{4r}}\nonumber\\
&\ \qquad + e^{- \bar{a}\min\{r+1, \kappa_*/2\}n}\cdot\frac{\sqrt{2\beta^{-1}}L^26^r\bar{C}_{1,2}}{2}\cdot\E\lsrs{|\bar{\theta}_0^{\lambda}|^{4r+4}}\nonumber\\
&\ \qquad +e^{-\bar{a}\kappa_* n/2}\cdot\frac{\sqrt{2\beta^{-1}}C_\nabla^2}{2}\cdot \E\lsrs{|\bar{\theta}_0^\lambda|^{4\nu+4}}\nonumber\\
&\ \qquad + \sqrt{2\beta^{-1}}L^2 6^r\left(2^{2r-2}\left(1+c_*(2r)\left(1+\frac{1}{\bar{a}\kappa_*}\right)\right)+\frac{v_{4r}(M_V(4r))}{2}+\frac{\bar{C}_{2,2}}{2}\right)\nonumber\\
&\ \qquad + \frac{\sqrt{2\beta^{-1}}C_\nabla^2}{2}\left(8d(d+2)+1+c_*(2\nu+2)\left(1+\frac{1}{\bar{a}\kappa_*}\right)\right)\Bigg{)}\nonumber\\
&\ + \lambda^2\Bigg{(}e^{-\bar{a}\kappa_* n/2}\cdot 4\bar{L}^{-1}K^2\cdot \E\lsrs{|\bar{\theta}_0^\lambda|^{6r+2}}+2\bar{L}^{-1}K^2\left(1+2c_*(3r+1)\left(1+\frac{1}{\bar{a}\kappa_*}\right)\right)\Bigg{)}\nonumber\\
\leq&\ 5\lambda\bar{L}\int^t_{nT}\E\lsrs{\lbrb{\bar{\zeta}_s^{\lambda, n} - \bar{\theta}_s^\lambda}^2}\dee s + \lambda^2 e^{-5\bar{L}}\left(e^{-\bar{a}\min\{r,\kappa_*/2\}n}\bar{C}_0\E\lsrs{|\bar{\theta}_0^\lambda|^{r_*}} +\bar{C}_1\right),\label{eqn:lemma4.9_wass2_ub2}
\end{align}
with 
\begin{align}
\bar{C}_0 :=&\ e^{5\bar{L}}\Bigg{(} 27\bar{L}^{-1}\bar{L}_\nabla^2 + \frac{\bar{L}^{-1}\bar{L}_\nabla^2\bar{C}_{1,1}}{2} + 12\bar{L}^{-1}C_\nabla^2 K^2 + \sqrt{2\beta^{-1}}L^2 6^r 2^{2r-2} + \frac{\sqrt{2\beta^{-1}}L^26^r\bar{C}_{1,2}}{2} \nonumber\\
&\ + \frac{\sqrt{2\beta^{-1}}C_\nabla^2}{2} + 4\bar{L}^{-1}K^2\Bigg{)},\nonumber\\
\bar{C}_1 := &\ \bar{C}_0 + e^{5\bar{L}}\Bigg{(} \bar{L}^{-1}\bar{L}_\nabla^2\left(\frac{27}{2} + 27c_*(2\nu)\left(1+\frac{1}{\bar{a}\kappa_*}\right) + \frac{\bar{C}_{2,1}}{2}\right)\nonumber\\
&\ + 12\bar{L}^{-1}C_\nabla^2 K^2\left(1+c_*(\nu+r+2)\left(1+\frac{1}{\bar{a}\kappa_*}\right)\right)\nonumber\\
&\ + \sqrt{2\beta^{-1}}L^2 6^r\left(2^{2r-2}\left(1+c_*(2r)\left(1+\frac{1}{\bar{a}\kappa_*}\right)\right)+\frac{v_{4r}(M_V(4r))}{2}+\frac{\bar{C}_{2,2}}{2}\right)\nonumber\\
&\ + \frac{\sqrt{2\beta^{-1}}C_\nabla^2}{2}\left(8d(d+2)+1+c_*(2\nu+2)\left(1+\frac{1}{\bar{a}\kappa_*}\right)\right)\nonumber\\
&\ +2\bar{L}^{-1}K^2\left(1+2c_*(3r+1)\left(1+\frac{1}{\bar{a}\kappa_*}\right)\right)\Bigg{)},\label{eqn:lemma4.9_constant_C0_C1}
\end{align}
where the constants $\kappa_*$ and $c_*(p)$ for $p\in\N$ are given in Corollary \ref{corollary:4.3}, $M_V(p)$ for $p\in\N_0$ is given in Lemma \ref{lemma:4.4}, $c_*(0):=1$ (see Remark \ref{remark:4.7}), and the constants $\bar{C}_{1,1}, \bar{C}_{2,1}, \bar{C}_{1,2}, \bar{C}_{2,2}$ are given in Lemma \ref{lemma:4.8}. Finally, by applying Gronwall's lemma to (\ref{eqn:lemma4.9_wass2_ub2}), one obtains
\begin{align}
\E\lsrs{\lbrb{\bar{\zeta}_t^{\lambda, n} - \bar{\theta}_t^\lambda}^2}\leq \lambda^2\left(e^{-\bar{a}\min\{r,\kappa_*/2\}n}\bar{C}_0\E\lsrs{|\bar{\theta}_0^\lambda|^{r_*}} +\bar{C}_1\right),\nonumber
\end{align}
which completes the proof. \qed
\subsection*{\textit{Proof of Lemma \ref{lemma:4.10}}}
\label{proof:lemma4.10}
We refer the reader to the proof of Lemma A.3 in \cite{lim2023non}. \qed

\subsection*{\textit{Proof of Proposition \ref{proposition:4:11}}} 
\label{proof:proposition4.11}
We verify that the assumptions of Theorem 2.2 of \cite{eberle2019quantitative} are satisfied. To this end, assign $\kappa\leftarrow\bar{L}$, $V\leftarrow V_2$, $C\leftarrow c_{V,2}(2)$, and $\lambda\leftarrow c_{V,1}(2)$, where here $\kappa$, $V$, $C$, and $\lambda$ are in the notation of \cite{eberle2019quantitative}. With these choices of constants, one sees that Assumptions 2.1 and 2.2 of \cite{eberle2019quantitative} hold due to \ref{remark:2.5} and Lemma \ref{lemma:4.4}, respectively. Moreover, it follows from (\ref{eqn:lemma4.4_Vp_gradient}) that $\frac{\nabla V_2(\theta)}{V_2(\theta)}\to 0$ as $|\theta|\to\infty$, which implies that Assumptions 2.4 and 2.5 of \cite{eberle2019quantitative} hold.\\

Now let 
\begin{align}
R_1 :=&\ \operatorname{diam}\{(\theta,\theta')\in \R^d\times\R^d: V_2(\theta)+V_2(\theta')\leq 2c_{V,2}(2)/c_{V,1}(2)\},\nonumber\\
R_2 :=&\ \operatorname{diam}\{(\theta,\theta')\in \R^d\times\R^d: V_2(\theta)+V_2(\theta')\leq 4c_{V,2}(2)(1+c_{V,1}(2))/c_{V,1}(2)\},\nonumber\\
Q(\epsilon) :=&\ \sup_{\theta\in\R^d}\frac{|\nabla V_2(\theta)|}{\max\{V_2(\theta),\epsilon^{-1}\}},\nonumber\\
\varphi(r;\epsilon) :=&\ \exp\{-\beta \bar{L}r^2/8-2Q(\epsilon)r\},\nonumber\\
\Phi(r;\epsilon) :=&\ \int^r_0\varphi(s;\epsilon)\ \dee s,\nonumber\\
g(r;\epsilon) :=& 1 - \frac{1}{4}\frac{\int^{r\wedge R_2}_0\beta\Phi(s;\epsilon)\varphi(s;\epsilon)^{-1}\ \dee s}{\int ^{R_2}_0\beta\Phi(s;\epsilon)\varphi(s;\epsilon)^{-1}\ \dee s} - \frac{1}{4}\frac{\int^{r\wedge R_1}_0\Phi(s;\epsilon)\varphi(s;\epsilon)^{-1}\ \dee s}{\int ^{R_1}_0\Phi(s;\epsilon)\varphi(s;\epsilon)^{-1}\ \dee s},\nonumber\\
f(r; \epsilon) :=& \int^{r\wedge R_2}_0 \varphi(s;\epsilon)g(s; \epsilon)\ \dee s,\nonumber\\
\phi_0:=&\ \frac{1}{\int^{R_2}_0\beta\Phi(s;\epsilon)\varphi(s;\epsilon)^{-1}\ \dee s},\nonumber\\
\dot{c}_0 :=&\ \min\{\phi_0, c_{V,1}(2)/2, 2c_{V,2}(2)\epsilon c_{V,1}(2)\}.\nonumber
\end{align}

Since the assumptions of Theorem 2.2 of \cite{eberle2019quantitative} hold, it follows, after a suitable rescaling, that for solutions $Z_t, Z_t'$, to the Langevin SDE (\ref{eqn:langevin}) with square-integrable initial conditions $Z_0=\theta_0$, $Z_0'=\theta_0'$, the inequality
\begin{align}
\mc{W}^\epsilon_\rho(\mc{L}(Z_t), \mc{L}(Z_t'))\leq e^{-\dot{c}_0 t}\mc{W}^\epsilon_\rho(\mc{L}(\theta_0),\mc{L}(\theta_0'))\label{eqn:proposition4.11_contraction}
\end{align}
holds for all $t\geq 0$ and any $\epsilon>0$ satisfying
\begin{align}
\frac{1}{2\beta c_{V,2}(2)\epsilon}\geq \int^{R_1}_0\Phi(s; \epsilon)\varphi(s; \epsilon)^{-1}\ \dee s,\label{eqn:propostion4.11_epsilon_constraint}
\end{align}

where $\mc{W}_\rho$ is the multiplicative semi-metric defined as
\begin{align}
\mc{W}^\epsilon_\rho(\mu,\mu'):=\inf_{\zeta\in\mc{C}(\mu,\mu')}\int_{\R^d\times\R^d} f(|\theta-\theta'|;\epsilon)(1+\epsilon V_2(\theta)+\epsilon V_2(\theta'))\ \dee\zeta(\theta,\theta').\label{eqn:proposition4.11_semimetric}
\end{align}

The remainder of the proof follows that of Lemma 3.24 of \cite{chau2021stochastic}. Using the same arguments as in \cite{chau2021stochastic}, one can derive the inequality $Q(\epsilon) \leq 1$ for all $\epsilon>0$. Furthermore, from the definitions of $R_1$ and $R_2$, one deduces that $R_1\leq \overline{R}_1$ and $\underline{R}_2\leq R_2\leq \overline{R}_2$, where
\begin{align}
\overline{R}_1 :=&\ 2(2c_{V,2}(2)/c_{V,1}(2)-1)^{1/2},\nonumber\\
\underline{R}_2 :=&\ (4c_{V,2}(2)(1+c_{V,1}(2))/c_{V,1}(2)-2)^{1/2},\nonumber\\
\overline{R}_2 :=&\ 2(4c_{V,2}(2)(1+c_{V,1}(2))/c_{V,1}(2)-1)^{1/2}.\nonumber
\end{align}
Using these estimates, one can bound the RHS of the inequality (\ref{eqn:propostion4.11_epsilon_constraint}) above by
\begin{align}
\int^{R_1}_0\Phi(s;\epsilon)\varphi(s;\epsilon)^{-1}\ \dee s
=&\ \int^{R_1}_0\int^s_0\exp\{\beta \bar{L}(s^2-r^2)/8+2Q(\epsilon)(s-r)\}\ \dee r\ \dee s\nonumber\\
\leq&\ \int^{R_1}_0\int^s_0\exp\{\beta \bar{L}(s^2-r^2)/8+2(s-r)\}\ \dee r\ \dee s\nonumber\\
=&\ \int^{R_1}_0\exp\left\{\frac{\beta \bar{L}}{8}\left(s+\frac{8}{\beta\bar{L}}\right)^2\right\}\int^s_0\exp\left\{-\frac{\beta \bar{L}}{8}\left(r+\frac{8}{\beta\bar{L}}\right)^2\right\}\ \dee r\ \dee s\nonumber\\
\leq&\ \int^{R_1}_0\exp\left\{\frac{\beta \bar{L}}{8}\left(s+\frac{8}{\beta\bar{L}}\right)^2\right\}\int_\R\exp\left\{-\frac{\beta \bar{L}}{8}r^2\right\}\ \dee r\ \dee s\nonumber\\
=&\ \sqrt{\frac{8\pi}{\beta \bar{L}}}\int^{R_1}_0\exp\left\{\left(s\sqrt{\beta\bar{L}/8} +\sqrt{8/(\beta\bar{L})}\right)^2\right\}\ \dee s\nonumber\\
\leq&\ \sqrt{\frac{8\pi}{\beta \bar{L}}}\int^{\overline{R}_1}_0\exp\left\{\left(s\sqrt{\beta\bar{L}/8} +\sqrt{8/(\beta\bar{L})}\right)^2\right\}\ \dee s.\label{eqn:proposition4.11_epsilon_constraint_relaxed}
\end{align}
By substituting (\ref{eqn:proposition4.11_epsilon_constraint_relaxed}) into (\ref{eqn:propostion4.11_epsilon_constraint}) and rearranging the terms, one sees that (\ref{eqn:proposition4.11_contraction}) holds for any choice of $\epsilon > 0$ satisfying
\begin{align}
\epsilon \leq 1\wedge \left(4c_{V,2}(2)\sqrt{2\pi\beta/\bar{L}}\int^{\overline{R}_1}_0\exp\left\{\left(s\sqrt{\beta\bar{L}/8} +\sqrt{8/(\beta\bar{L})}\right)^2\right\}\ \dee s\right)^{-1}.\label{eqn:proposition4.11_epsilon_constraint_simplified}
\end{align}

Similarly, replacing $R_1$ and $\overline{R}_1$ in the above argument with $R_2$ and $\overline{R}_2$ respectively yields the estimate
\begin{align}
\phi_0^{-1}
=&\ \beta\int^{R_2}_0 \Phi(s;\epsilon)\varphi(s;\epsilon)^{-1}\ \dee s\nonumber\\
\leq&\ \sqrt{\frac{8\pi\beta}{\bar{L}}}\int^{\overline{R}_2}_0\exp\left\{\left(s\sqrt{\beta\bar{L}/8} +\sqrt{8/(\beta\bar{L})}\right)^2\right\}\ \dee s\nonumber\\
\leq&\ \sqrt{\frac{8\pi\beta}{\bar{L}}}\overline{R}_2\exp\left\{\left(\overline{R}_2\sqrt{\beta\bar{L}/8} +\sqrt{8/(\beta\bar{L})}\right)^2\right\}=:\phi^{-1}. \label{eqn:proposition4.11_phi_lb}
\end{align}

Therefore, (\ref{eqn:proposition4.11_contraction}) and (\ref{eqn:proposition4.11_phi_lb}) together imply that for all $t\geq 0$ and $\epsilon >0$ satisfying (\ref{eqn:proposition4.11_epsilon_constraint_simplified}), 
\begin{align}
\mc{W}^\epsilon_\rho(\mc{L}(Z_t), \mc{L}(Z_t'))\leq e^{-\dot{c} t}\mc{W}^\epsilon_\rho(\mc{L}(\theta_0),\mc{L}(\theta_0')),\label{eqn:proposition4.11_contraction2}
\end{align}
where $\dot{c} := \min\{\phi, c_{V,1}(2)/2, 2c_{V,2}(2)\epsilon c_{V,1}(2)\}$.\\

It remains to relate the semi-metric $\mc{W}_\rho^\epsilon$ to our functional $w_{1,2}$. Fix an $\epsilon>0$ satisfying (\ref{eqn:proposition4.11_epsilon_constraint_simplified}). From now, we omit the explicit dependence of our functions of interest on $\epsilon$ for brevity. From Equation 5.4 of \cite{eberle2019quantitative} and $Q(\epsilon)\leq 1$, the inequalities
\begin{align}
\frac{1}{2}r\exp\{-\beta\bar{L}R_2^2/8-2R_2\}\leq \frac{1}{2}\Phi(r)\leq f(r) \leq \Phi(r) \leq r\label{eqn:proposition4.11_inequalities}
\end{align}
hold true for $r\in [0,R_2]$. Furthermore, it is clear from the definition of $f$ that $f(r)=f(R_2)$ for all $r\geq R_2$. Let $\theta,\theta'\in \R^d$ and $r=|\theta-\theta'|$. Suppose $r\leq R_2$. One then obtains, from the inequalities in (\ref{eqn:proposition4.11_inequalities}),
\begin{align}
(1\wedge |\theta - \theta'|)(1+ V_2(\theta)+ V_2(\theta'))
\leq &\ \epsilon^{-1}|\theta-\theta'|(\epsilon+\epsilon V_2(\theta)+\epsilon V_2(\theta'))\nonumber\\
\leq &\ \epsilon^{-1}|\theta-\theta'|(1+\epsilon V_2(\theta)+\epsilon V_2(\theta'))\nonumber\\
\leq &\ 2\epsilon^{-1}\exp\{\beta\bar{L}R_2^2/8+2R_2\}f(|\theta-\theta'|)(1+\epsilon V_2(\theta)+\epsilon V_2(\theta'))\nonumber\\
\leq &\ 2\epsilon^{-1}\exp\{\beta\bar{L}\overline{R}_2^2/8+2\overline{R}_2\}f(|\theta-\theta'|)(1+\epsilon V_2(\theta)+\epsilon V_2(\theta')),\label{eqn:proposition4.11_tightbound1}
\end{align}
where the second and third inequalities are due to $\epsilon \leq 1$ and $r\leq 2f(r)\exp\{\beta\bar{L}R_2^2/8+R_2\}$ for every $r\in[0,R_2]$, respectively, as well as
\begin{align}
f(|\theta-\theta'|)(1+\epsilon V_2(\theta)+\epsilon V_2(\theta'))
\leq&\ |\theta-\theta'|(1+\epsilon V_2(\theta)+\epsilon V_2(\theta'))\nonumber\\
\leq&\ |\theta-\theta'|(1+ V_2(\theta)+ V_2(\theta'))\nonumber\\
=&\ (|\theta-\theta'|\1_{\{r\leq 1\}}+|\theta-\theta'|\1_{\{r> 1\}})(1+ V_2(\theta)+ V_2(\theta'))\nonumber\\
=&\  (\1_{\{r\leq 1\}}+|\theta-\theta'|\1_{\{r> 1\}})(1\wedge|\theta-\theta'|)(1+ V_2(\theta)+ V_2(\theta'))\nonumber\\
\leq&\  (\1_{\{r\leq 1\}}+R_2\1_{\{r> 1\}})(1\wedge|\theta-\theta'|)(1+ V_2(\theta)+ V_2(\theta'))\nonumber\\
\leq&\ (1+\overline{R}_2)(1\wedge|\theta-\theta'|)(1+ V_2(\theta)+ V_2(\theta')).\label{eqn:proposition4.11_tightbound2}
\end{align}

In the case $r>R_2$, one obtains, similarly,
\begin{align}
(1\wedge |\theta - \theta'|)(1+ V_2(\theta)+ V_2(\theta'))
\leq &\ \epsilon^{-1}(1+\epsilon V_2(\theta)+\epsilon V_2(\theta'))\nonumber\\
\leq &\ \frac{2}{\epsilon R_2}\exp\{\beta\bar{L}R_2^2/8+2R_2\}f(R_2)(1+\epsilon V_2(\theta)+\epsilon V_2(\theta'))\nonumber\\
= &\ \frac{2}{\epsilon R_2}\exp\{\beta\bar{L}R_2^2/8+2R_2\}f(|\theta-\theta'|)(1+\epsilon V_2(\theta)+\epsilon V_2(\theta'))\nonumber\\
\leq &\ 2(\epsilon \underline{R}_2)^{-1}\exp\{\beta\bar{L}\overline{R}_2^2/8+2\overline{R}_2\}f(|\theta-\theta'|)(1+\epsilon V_2(\theta)+\epsilon V_2(\theta'))\label{eqn:proposition4.11_tightbound3}
\end{align}
and
\begin{align}
f(|\theta-\theta'|)(1+\epsilon V_2(\theta)+\epsilon V_2(\theta'))
=&\ f(R_2)(1+\epsilon V_2(\theta)+\epsilon V_2(\theta'))\nonumber\\
\leq&\ R_2(1+\epsilon V_2(\theta)+\epsilon V_2(\theta'))\nonumber\\
\leq&\ R_2(1+ V_2(\theta)+ V_2(\theta'))\nonumber\\
=&\ (R_2\1_{\{R_2\leq 1\}}+R_2\1_{\{R_2> 1\}})(1+ V_2(\theta)+ V_2(\theta'))\nonumber\\
=&\  (\1_{\{R_2\leq 1\}}+R_2\1_{\{R_2> 1\}})(1\wedge R_2)(1+ V_2(\theta)+ V_2(\theta'))\nonumber\\
\leq&\ (1+\overline{R}_2)(1\wedge|\theta-\theta'|)(1+ V_2(\theta)+ V_2(\theta')).\label{eqn:proposition4.11_tightbound4}
\end{align}

Therefore, (\ref{eqn:proposition4.11_tightbound1}), (\ref{eqn:proposition4.11_tightbound2}), (\ref{eqn:proposition4.11_tightbound3}), and (\ref{eqn:proposition4.11_tightbound4}) together imply that
\begin{align}
(\epsilon/2)\exp\{-\beta\bar{L}\overline{R}_2^2/8-2\overline{R}_2\}w_{1,2}\leq\mc{W}_\rho^\epsilon \leq (1+\overline{R}_2)w_{1,2}.\label{eqn:proposition4.11_tightbound}
\end{align}

Substituting (\ref{eqn:proposition4.11_tightbound}) into (\ref{eqn:proposition4.11_contraction2}), one obtains
\begin{align}
w_{1,2}(\mc{L}(Z_t),\mc{L}(Z_t'))
\leq &\ 2\epsilon^{-1}\exp\{\beta\bar{L}\overline{R}_2^2/8+2\overline{R}_2\}\mc{W}_\rho^\epsilon(\mc{L}(Z_t),\mc{L}(Z_t'))\nonumber\\
\leq &\ 2\epsilon^{-1}\exp\{\beta\bar{L}\overline{R}_2^2/8+2\overline{R}_2\}e^{-\dot{c}t}\mc{W}_\rho^\epsilon(\mc{L}(\theta_0),\mc{L}(\theta_0'))\nonumber\\
\leq &\ \hat{c}e^{-ct}w_{1,2}(\mc{L}(\theta_0),\mc{L}(\theta_0')),\nonumber
\end{align}
where $\hat{c}:= 2(1+\overline{R}_2)\exp\{\beta\bar{L}\overline{R}_2^2/8+2\overline{R}_2\}/\epsilon$. This completes the proof. \qed

\subsection*{\textit{Proof of Lemma \ref{lemma:4.12}}}
\label{proof:lemma4.12}
We follow the idea of the proof of Lemma 4.7 of \cite{lim2023non}. Let $\lambda\in(0,\tilde{\lambda}_{\max})$, $n\in\N_0$, and let $t\in(nT, (n+1)T]$. By applying successively the triangle inequality, Definition \ref{definition:4.1}, Lemma \ref{lemma:4.10}, and the Cauchy-Schwarz and Minkowski inequalities, we have
\begin{align}
&W_1(\mc{L}(\bar{\zeta}_t^{\lambda,n}),\mc{L}(Z_t^\lambda))\nonumber\\
\leq&\ \sum^n_{k=1} W_1(\mc{L}(\bar{\zeta}_t^{\lambda,k}),\mc{L}(\bar{\zeta}_t^{\lambda,k-1}))\nonumber\\
=&\ \sum^n_{k=1} W_1(\mc{L}(\zeta_t^{kT, \bar{\theta}_{kT}^\lambda, \lambda}),\mc{L}(\zeta_t^{kT, \bar{\zeta}_{kT}^{\lambda, k-1}, \lambda}))\nonumber\\
\leq&\ \sum^n_{k=1}w_{1,2}(\mc{L}(\zeta_t^{kT, \bar{\theta}_{kT}^\lambda, \lambda}),\mc{L}(\zeta_t^{kT, \bar{\zeta}_{kT}^{\lambda, k-1}, \lambda}))\nonumber\\
\leq&\ \hat{c}\sum^n_{k=1}e^{-\dot{c}\lambda (t-kT)}w_{1,2}(\mc{L}(\bar{\theta}_{kT}^\lambda),\mc{L}(\bar{\zeta}_{kT}^{\lambda, k-1}))\nonumber\\
\leq&\ \hat{c}\sum^n_{k=1}e^{-\dot{c}(n-k)/2}w_{1,2}(\mc{L}(\bar{\theta}_{kT}^\lambda),\mc{L}(\bar{\zeta}_{kT}^{\lambda, k-1}))\nonumber\\
=&\ \hat{c}\sum^n_{k=1}e^{-\dot{c}(n-k)/2}\inf_{\psi\in\mc{C}(\mc{L}(\bar{\theta}_{kT}^\lambda),\mc{L}(\bar{\zeta}_{kT}^{\lambda, k-1}))}\int_{\R^d\times\R^d}(1\wedge |\theta-\theta'|)(1+V_2(\theta)+V_2(\theta'))\ \dee\psi(\theta,\theta')\nonumber\\
\leq& \ \hat{c}\sum^n_{k=1}e^{-\dot{c}(n-k)/2}\inf_{\psi\in\mc{C}(\mc{L}(\bar{\theta}_{kT}^\lambda),\mc{L}(\bar{\zeta}_{kT}^{\lambda, k-1}))}\Bigg{(}\left\{\int_{\R^d\times\R^d}|\theta-\theta'|^2\ \dee\psi(\theta,\theta')\right\}^{1/2}\nonumber\\
&\ \times \left\{\int_{\R^d\times\R^d}(1+V_2(\theta)+V_2(\theta'))^2\ \dee\psi(\theta,\theta')\right\}^{1/2}\Bigg{)}\nonumber\\
\leq&\ \hat{c}\sum^n_{k=1}e^{-\dot{c}(n-k)/2}W_2(\mc{L}(\bar{\theta}_{kT}^\lambda),\mc{L}(\bar{\zeta}_{kT}^{\lambda, k-1}))\left\lVert1+V_2(\bar{\theta}_{kT}^\lambda)+V_2(\bar{\zeta}_{kT}^{\lambda, k-1})\right\rVert_{L^2(\Omega)}\nonumber\\
\leq&\ \hat{c}\sum^n_{k=1}e^{-\dot{c}(n-k)/2}W_2(\mc{L}(\bar{\theta}_{kT}^\lambda),\mc{L}(\bar{\zeta}_{kT}^{\lambda, k-1}))\lsrs{1+\left\{\E[V_4(\bar{\theta}_{kT}^\lambda)]\right\}^{1/2} + \left\{\E[V_4(\bar{\zeta}_{kT}^{\lambda, k-1})]\right\}^{1/2}}.\label{eqn:lemma4.12_wass1_ub}
\end{align}
In addition, by Young's inequality, Lemma \ref{lemma:4.9} and Corollaries \ref{corollary:4.6} and \ref{corollary:4.3}, one further obtains
\begin{align}
&W_1(\mc{L}(\bar{\zeta}_t^{\lambda,n}),\mc{L}(Z_t^\lambda))\nonumber\\
\leq&\ \lambda^{-1}\hat{c}\sum^n_{k=1}e^{-\dot{c}(n-k)/2}W_2^2(\mc{L}(\bar{\theta}_{kT}^\lambda),\mc{L}(\bar{\zeta}_{kT}^{\lambda, k-1})) + \frac{3\lambda }{4}\hat{c}\sum^n_{k=1}e^{-\dot{c}(n-k)/2}\lsrs{1+\E[V_4(\bar{\theta}_{kT}^\lambda)]+ \E[V_4(\bar{\zeta}_{kT}^{\lambda, k-1})]}\nonumber\\
\leq&\ \lambda\hat{c}\sum^n_{k=1}e^{-\dot{c}(n-k)/2}\lsrs{e^{-\bar{a}\min\{r,\kappa_*/2\}(k-1)}\bar{C}_0\E\lsrs{|\bar{\theta}_0^\lambda|^{r_*}}+\bar{C}_1}\nonumber\\
&\ +\frac{3\lambda}{4}\hat{c}\sum^n_{k=1}e^{-\dot{c}(n-k)/2}\Bigg{[}1+2+2e^{-\bar{a}\kappa_*(k-1)/2}\E\lsrs{|\bar{\theta}_0^\lambda|^4}+2c_*(2)\left(1+\frac{1}{\bar{a}\kappa_*}\right)\nonumber\\
&\ +2+2e^{-\bar{a}\min\{2,\kappa_*\}(k-1)/2}\E\lsrs{|\bar{\theta}_0^\lambda|^4}+2c_*(2)\left(1+\frac{1}{\bar{a}\kappa_*}\right)+v_4(M_V(4))\Bigg{]}\nonumber\\
\leq&\ \lambda\hat{c}\sum^n_{k=1}e^{-\dot{c}(n-k)/2}\lsrs{e^{-\bar{a}\min\{r,\kappa_*/2\}(k-1)}\bar{C}_0\E\lsrs{|\bar{\theta}_0^\lambda|^{r_*}}+\bar{C}_1}\nonumber\\
&\ +\frac{3\lambda}{4}\hat{c}\sum^n_{k=1}e^{-\dot{c}(n-k)/2}\left[9+4e^{-\bar{a}\min\{1,\kappa_*/2\}(k-1)}\E\lsrs{|\bar{\theta}_0^\lambda|^{r_*}}+4c_*(2)\left(1+\frac{1}{\bar{a}\kappa_*}\right)+v_4(M_V(4))\right]\nonumber\\
=&\ \lambda\hat{c}\sum^n_{k=1}e^{-\dot{c}(n-k)/2}\lsrs{e^{-\bar{a}\min\{1, r,\kappa_*/2\}(k-1)}(\bar{C}_0+3)\E\lsrs{|\bar{\theta}_0^\lambda|^{r_*}}}\nonumber\\
&\ +\lambda\hat{c}\sum^n_{k=1}e^{-\dot{c}(n-k)/2}\left[\bar{C}_1+\frac{27}{4}+3c_*(2)\left(1+\frac{1}{\bar{a}\kappa_*}\right)+\frac{3}{4}v_4(M_V(4))\right]\nonumber\\
\leq &\ \lambda\hat{c}\sum^n_{k=1}e^{-\min\{\dot{c}/2, \bar{a}, \bar{a}r,\bar{a}\kappa_*/2\}(n-1)}(\bar{C}_0+3)\E\lsrs{|\bar{\theta}_0^\lambda|^{r_*}}\nonumber\\
&\ +\lambda\hat{c}\sum^\infty_{k=0}e^{-\dot{c}k/2}\left[\bar{C}_1+\frac{27}{4}+3c_*(2)\left(1+\frac{1}{\bar{a}\kappa_*}\right)+\frac{3}{4}v_4(M_V(4))\right]\nonumber\\
=&\ \lambda \hat{c}n e^{-\min\{\dot{c}/2, \bar{a}, \bar{a}r,\bar{a}\kappa_*/2\}(n-1)}(\bar{C}_0+3)\E\lsrs{|\bar{\theta}_0^\lambda|^{r_*}}\nonumber\\
&\ + \lambda\cdot \frac{\hat{c}}{1-e^{-\dot{c}/2}}\left[\bar{C}_1+\frac{27}{4}+3c_*(2)\left(1+\frac{1}{\bar{a}\kappa_*}\right)+\frac{3}{4}v_4(M_V(4))\right]\nonumber\\
\leq&\ \lambda \hat{c}\cdot \frac{1+(n-1)}{1+\min\{\dot{c}/2, \bar{a}, \bar{a}r,\bar{a}\kappa_*/2\}(n-1)/2}e^{-\min\{\dot{c}/2, \bar{a}, \bar{a}r,\bar{a}\kappa_*/2\}(n-1)/2} (\bar{C}_0+3)\E\lsrs{|\bar{\theta}_0^\lambda|^{r_*}}\nonumber\\
&\ +\lambda\cdot \frac{2\hat{c}e^{\dot{c}/2}}{\dot{c}}\left[\bar{C}_1+\frac{27}{4}+3c_*(2)\left(1+\frac{1}{\bar{a}\kappa_*}\right)+\frac{3}{4}v_4(M_V(4))\right]\nonumber\\
\leq &\ \lambda\left(e^{-\min\{\dot{c}/4, \bar{a}/2, \bar{a}r/2,\bar{a}\kappa_*/4\}n}\bar{C}_2\E\lsrs{|\bar{\theta}_0^\lambda|^{r_*}} + \bar{C}_3\right),\nonumber
\end{align}
where we used $e^s\geq 1+s$, for all $s\in \R$ in the second last inequality, with
\begin{align}
\bar{C}_2 &:=\hat{c}e^{\min\{\dot{c}/4, \bar{a}/2, \bar{a}r/2,\bar{a}\kappa_*/4\}}\left(1+\frac{1}{\min\{\dot{c}/4, \bar{a}/2, \bar{a}r/2,\bar{a}\kappa_*/4\}}\right)(\bar{C}_0+3),\nonumber\\
\bar{C}_3 &:= \frac{2\hat{c}e^{\dot{c}/2}}{\dot{c}}\left[\bar{C}_1+\frac{27}{4}+3c_*(2)\left(1+\frac{1}{\bar{a}\kappa_*}\right)+\frac{3}{4}v_4(M_V(4))\right],\label{eqn:lemma4.12_constant_C2_C3}
\end{align}
where the constants $\kappa_*$ and $c_*(2)$ are given in Corollary \ref{corollary:4.3}, $M_V(4)$ is given in Lemma \ref{lemma:4.4}, $c_*(0):=1$ (see Remark \ref{remark:4.7}), $\bar{C}_0$ and $\bar{C}_1$ are given in Lemma \ref{lemma:4.9}, and $\dot{c}$ and $\hat{c}$ are given in Proposition \ref{proposition:4:11}. This completes the proof. \qed

\section*{Acknowledgement}
Financial support by the MOE AcRF Tier 2 Grant \textit{MOE-T2EP20222-0013}
and the Guangzhou-HKUST(GZ) Joint Funding Program (No. 2024A03J0630) is gratefully acknowledged.
\newpage
\section*{Appendix A. Analytic Expression of Constants}
\label{appendix:A}
\begin{tabularx}{\textwidth}[H]
{ 
	>{\hsize=0.5\hsize\linewidth=\hsize}X
	>{\hsize=0.3\hsize\linewidth=\hsize}X|
	>{\hsize=2.2\hsize\linewidth=\hsize}X
} 
{CONSTANT} & & FULL EXPRESSION\\
\hline
Remark \ref{remark:2.4} & $R$ & $\max\left\{\left(\tfrac{4b}{a}\right)^{1/(r-\bar{r})},2^{1/r}\right\}$ \\

 & $\bar{a}$ & $\frac{a}{2}\1_{\{r>0\}}+\tilde{a}\1_{\{r=0\}}$\\

 & $\bar{b}$ & $\left(\left(b+\tfrac{a}{2}\right)R^{\bar{r}+2}+\tfrac{K^2}{2a}\right)\1_{\{r>0\}}+\tilde{b}\1_{\{r=0\}}$\\

 & $\bar{b}'$ & $(\bar{b}+2^{2/r}\bar{a})\1_{\{r>0\}} + \tilde{b}\1_{\{r=0\}}$ \\
\hline
Remark \ref{remark:2.5} & $\bar{R}$ & $(b/a)^{1/(r-\bar{r})}$ \\
 & $\bar{L}$ & $L(1+2\bar{R})^{r}$ \\
\hline

Remark \ref{remark:2.6} & $C_\nabla$ &  $2\max\{2^{\nu-1}L_\nabla, |\nabla h(0)|\}$\\
 & $\bar{L}_\nabla$ & $3^{\nu - 1}L_\nabla$\\
\hline

Lemma  \ref{lemma:4.2} & $\kappa$ &  $\tfrac{1}{\sqrt{2}}$\\

 & $c_0$ &  $\bar{a}\kappa+2\bar{b}+2d\beta^{-1}+2K^2$ \\
 
 & $M_1(p)$ & $\frac{4p\binom{p}{\lcrc{p/2}}(1+2\bar{b}+2K^2)^p}{\min\{1,\bar{a}\}}$\\
 
 & $\tilde{\kappa}(p)$ & $\frac{[M_1(p)]^r}{2(1+[M_1(p)]^{2r})^{1/2}}$\\
 & $M_2(p)$ & $\lsrs{\frac{p(2p-1)2^{2p-1}\beta^{-1}d}{\bar{a}\tilde{\kappa}(p)}}^{1/2}$\\
 
 & $c_1(p)$ & $\bar{a}\tilde{\kappa}(p)[M_1(p)]^{2p}+\sum^p_{k=1}\binom{p}{k}(2\bar{b}+2K^2)^k [M_1(p)]^{2p-2k}$\\
 
 & $c_2(p)$ & $c_1(p) + p(2p-1)2^{2p-2}\beta^{-1}d\ c_1(p-1)+ p(2p-1)2^{4p-3}\beta^{-p}p!\binom{\frac{d}{2}+p-1}{p}$\\
 
 & $c_3(p)$ & $c_2(p)+p(2p-1)2^{2p-2}\beta^{-1}d\ [M_2(p)]^{2p-2}$\\
\hline

Corollary \ref{corollary:4.3} & $\kappa_*$ & $\min\{\kappa, \frac{\tilde{\kappa}(2)}{2}\}$\\
 & $c_*(p)$ & $\1_{\{p=0\}}+\max\{c_0,c_3(p)\1_{\{p\geq 2\}}\}$\\
\hline

Lemma \ref{lemma:4.4} & $M_V(p)$ & $\frac{\bar{a} p}{2}$\\

 & $c_{V,1}(p)$ & $\frac{\bar{a} p}{2}$\\
 
 & $c_{V,2}(p)$ & $\left(1+\frac{2\bar{b}'+2\beta^{-1}(d+p-2)}{\bar{a}}\right)^{1/2}$\\
\hline

Lemma \ref{lemma:4.8} & $\bar{C}_{1,1}$ & $16384K^8$\\

 & $\bar{C}_{2,1}$ & $16384K^8\left(1+c_*(4r+4)\left(1+\frac{1}{\bar{a}\kappa_*}\right)\right) + 2048\beta^{-4}d(d+2)(d+4)(d+6)$ \\
 
 & $\bar{C}_{1,2}$ & $2^{2r+7}K^4$ \\
 
 & $\bar{C}_{2,2}$ & $64K^4\lsrs{2^{2r+1}+2^{2r+1}c_*(2r+2)\left(1+\frac{1}{\bar{a}\kappa_*}\right) + v_{4r+4}(M_V(4r+4))}$\\
\hline

Lemma \ref{lemma:4.9} & $\bar{C}_0$ & $e^{5\bar{L}}\Bigg{(} 27\bar{L}^{-1}\bar{L}_\nabla^2 + \frac{\bar{L}^{-1}\bar{L}_\nabla^2\bar{C}_{1,1}}{2} + 12\bar{L}^{-1}C_\nabla^2 K^2 + \sqrt{2\beta^{-1}}L^2 6^r 2^{2r-2} + \frac{\sqrt{2\beta^{-1}}L^26^r\bar{C}_{1,2}}{2} + \frac{\sqrt{2\beta^{-1}}C_\nabla^2}{2} + 4\bar{L}^{-1}K^2\Bigg{)}$ \\

 & $\bar{C}_1$ & $\bar{C}_0 + e^{5\bar{L}}\Bigg{(} \bar{L}^{-1}\bar{L}_\nabla^2\left(\frac{27}{2} + 27c_*(2\nu)\left(1+\frac{1}{\bar{a}\kappa_*}\right) + \frac{\bar{C}_{2,1}}{2}\right)$\newline
 
 $+ 12\bar{L}^{-1}C_\nabla^2 K^2\left(1+c_*(\nu+r+2)\left(1+\frac{1}{\bar{a}\kappa_*}\right)\right)$
 \newline
 
 $+ \sqrt{2\beta^{-1}}L^2 6^r\left(2^{2r-2}\left(1+c_*(2r)\left(1+\frac{1}{\bar{a}\kappa_*}\right)\right)+\frac{v_{4r}(M_V(4r))}{2}+\frac{\bar{C}_{2,2}}{2}\right)$\newline
 
 $+ \frac{\sqrt{2\beta^{-1}}C_\nabla^2}{2}\left(8d(d+2)+1+c_*(2\nu+2)\left(1+\frac{1}{\bar{a}\kappa_*}\right)\right)$\newline

 $ +2\bar{L}^{-1}K^2\left(1+2c_*(3r+1)\left(1+\frac{1}{\bar{a}\kappa_*}\right)\right)\Bigg{)}$ \\
\hline

Proposition \ref{proposition:4:11} & $\overline{R}_1$ & $2(2c_{V,2}(2)/c_{V,1}(2)-1)^{1/2}$\\

 & $\overline{R}_2$ & $2(4c_{V,2}(2)(1+c_{V,1}(2))/c_{V,1}(2)-1)^{1/2}$\\
 
 & $\epsilon$ & $\in\left(0,1\wedge \left(4c_{V,2}(2)\sqrt{2\pi\beta/\bar{L}}\int^{\overline{R}_1}_0\exp\left\{\left(s\sqrt{\beta\bar{L}/8} +\sqrt{8/(\beta\bar{L})}\right)^2\right\}\ \dee s\right)^{-1}\right]$ \\
 
 & $\hat{c}$ & $2(1+\overline{R}_2)\exp\{\beta\bar{L}\overline{R}_2^2/8+2\overline{R}_2\}/\epsilon$\\
 
 & $\dot{c}$ & $\min\Bigg{\{}\left(\overline{R}_2\sqrt{8\pi\beta/\bar{L}}\exp\left\{\left(\overline{R}_2\sqrt{\beta\bar{L}/8}+\sqrt{8/(\beta\bar{L})}\right)^2\right\}\right)^{-1},$ \newline 
 $ c_{V,1}(2)/2,\ 2c_{V,2}(2)\epsilon c_{V,1}(2)\Bigg{\}}$\\
\hline

Lemma \ref{lemma:4.12} & $\bar{C}_2$ & $\hat{c}e^{\min\{\dot{c}/4, \bar{a}/2, \bar{a}r/2,\bar{a}\kappa_*/4\}}\left(1+\frac{1}{\min\{\dot{c}/4, \bar{a}/2, \bar{a}r/2,\bar{a}\kappa_*/4\}}\right)(\bar{C}_0+3)$\\

 & $\bar{C}_3$ & $\frac{2\hat{c}e^{\dot{c}/2}}{\dot{c}}\left[\bar{C}_1+\frac{27}{4}+3c_*(2)\left(1+\frac{1}{\bar{a}\kappa_*}\right)+\frac{3}{4}v_4(M_V(4))\right]$\\

\hline 
Theorem \ref{theorem:2.7} & $C_0$ & $\min\{\dot{c}/4, \bar{a}/2, \bar{a}r/2,\bar{a}\kappa_*/4\}$ \\
 & $C_1$ & $e^{C_0}\lsrs{\bar{C}_0^{1/2}+\bar{C}_2+\hat{c}\left(3+\int_{\R^d} V_2(\theta)\ \dee\pi_\beta(\theta)\right)}$\\
 & $C_2$ & $\bar{C}_1^{1/2}+\bar{C}_3$\\

\hline

Theorem \ref{theorem:2.8} & $\bar{C}_4$ & $\sqrt{2\hat{c}}e^{\min\{\dot{c}/8,\bar{a}/4,\bar{a}r/4,\bar{a}\kappa_*/8\}}\left(1+\frac{1}{\min\{\dot{c}/8,\bar{a}/4,\bar{a}r/4,\bar{a}\kappa_*/8\}}\right)\left(\bar{C}_0^{1/2}+\frac{1}{\sqrt{2}}\right)$\\
 & $\bar{C}_5$ & $\frac{4\sqrt{2\hat{c}}e^{\dot{c}/4}}{\dot{c}}\left(\bar{C}_1^{1/2}+\frac{1+2\sqrt{2}}{4}+\sqrt{\frac{c_*(2)}{2}}\left(1+\frac{1}{\bar{a}\kappa_*}\right)^{1/2}+\frac{(v_4(M_V(4)))^{1/2}}{4}\right)$\\
 & $C_3$ & $\min\{\dot{c}/8, \bar{a}/4, \bar{a}r/4,\bar{a}\kappa_*/8\}$ \\
 & $C_4$ & $e^{C_3}\lsrs{\bar{C}_0^{1/2}+\bar{C}_4+\sqrt{2\hat{c}}\left(1+\left(2+\int_{\R^d} V_2(\theta)\ \dee\pi_\beta(\theta)\right)^{1/2}\right)}$ \\
 & $C_5$ & $\bar{C}_1^{1/2}+\bar{C}_5$ \\
\end{tabularx}
\newpage
%
%
%
%
\printbibliography
\end{document}